\documentclass[10pt,a4paper,twoside]{article}
\pdfoutput=1

\usepackage[margin=3cm,top=3cm,bottom=3cm]{geometry}

\usepackage{monstyle}
\usepackage{import}
\usepackage{ifthen}
\usetikzlibrary{decorations.pathmorphing}

\newtheorem{thmintrovariant}{Theorem}

\newcommand{\col}{\colon}

\newcommand{\ar}{\stackrel}

\newcommand{\LX}{\mathcal{L}X}
\newcommand{\tri}{\triangleleft}
\newcommand{\oact}{{\zO^\o_{\mathrm{act}}}}
\newcommand{\atom}{{\mathrm{Atom}_\zO}}

\newcommand{\Mark}{{\mathrm{Mark}}}
\newcommand{\dst}{{\mathrm{d}\mathscr{S}\mathrm{t}}}
\newcommand{\dgCat}{{\mathrm{dg}\mathscr{C}\mathrm{at}}}

\newcommand{\sspan}{{\mathsf{Span}_2}}

\newcommand{\bo}{\mathscr{BO}}
\newcommand{\lift}{\mathscr{L}}

\renewcommand{\Fin}{{\mathbb{F}}}

\tikzcdset{column sep/verytiny/.initial=.3em}
\tikzcdset{crossing over clearance=.8ex}
\tikzcdset{scale cd/.style={every label/.append style={scale=#1},
    cells={nodes={scale=#1}}}}

\begin{document}
    \title
    {Brane actions for coherent \texorpdfstring{$\infty$}{infinity}-operads}
    \author{Hugo Pourcelot\footnote{
    Université Sorbonne Paris Nord, LAGA, CNRS, UMR 7539, F-93430, Villetaneuse, France
    }}
	\date{February 23, 2023}
    \maketitle

\begin{abstract}
    We prove that Mann--Robalo's construction of the brane action \cite{robalo} extends to general coherent $\infty$-operads, with possibly multiple colors and non-contractible spaces of unary operations.
    This requires to establish two results regarding spaces of extensions that were left unproven in the aforementioned construction. 
    First, we show that Lurie's and Mann--Robalo's models for such spaces are equivalent.
    Second, we prove that the space of extensions in the sense of Lurie is not in general equivalent to the homotopy fiber 
    of the associated forgetful morphism, but rather to its homotopy quotient  by the $\infty$-groupoid of unary operations, correcting an oversight in existing literature.
    As an application, we obtain that the $\infty$-operads of $B$-framed little disks are coherent and therefore yield new operations on spaces of branes of perfect derived stacks.
\end{abstract}

\tableofcontents

\pagestyle{plain}

	\section{Introduction} 
\label{sec:introduction}

%


\paragraph*{The brane action.}
Introduced by Toën in \cite{toen_branes}, reformulated and further developed by Mann and Robalo in \cite{robalo,robalo_mann_survol_smf} and Kern in \cite{kern_these}, the brane action is an operadic construction that encodes in a unified way certain "configuration-type" invariants of geometric objects. The motivating examples are that of Gromov--Witten invariants from enumerative geometry and string topology operations; in both cases, the brane action reflects part of the structure of an associated topological conformal field theory, in the sense of Costello \cite{costello_2007_TCFT_CY}.

Formally, given a reduced\footnote{Recall that an $\infty$-operad $\zO$ is called \emph{reduced} if its underlying $\infty$-category is a trivial $\infty$-groupoid. In the monochromatic case, this means that $\zO(1)\simeq *$.} coherent monochromatic $\infty$-operad $\zO^\o$, the brane action consists in associating to every operation $\sigma \in \zO(k)$ its space of extensions $\Ext(\sigma)$, informally described as parametrizing those operations $\sigma^+$ of arity $k+1$ that restrict to $\sigma$ on the first $k$ inputs (see definition \ref{df:Ext(sigma)} for details). 
Each such space $\Ext(\sigma)$ then fits into a cospan diagram  
\begin{equation}
    \label{eqn:cospan_diagram}
    \begin{tikzcd}
        \Ext(\id)^{\amalg k} \ar[r] & \Ext(\sigma) & \ar[l] \Ext(\id)
\end{tikzcd}
\end{equation}
which altogether assembles into an $\zO$-algebra structure on $\Ext(\id)$, viewed as an object in the $\infty$-category of cospans of spaces.

We refer to section \ref{sec:the_brane_fibration} for precise definitions of spaces of extensions and {coherent} $\infty$-operads, as well as Mann--Robalo's construction of the associated brane action, which we will extend to the \emph{non-reduced coloured} setting.
Let us mention that a large part of our work is actually devoted to providing the first proofs in literature of several claims used in their work, thereby completing the construction of brane actions from \cite{robalo}.

\paragraph*{Application to Gromov--Witten invariants and string topology.}
The first application of this construction is the categorification of Gromov--Witten invariants, due to Mann and Robalo. In \cite{robalo}, the authors showed that the brane action, when instantiated with a variant of the $\infty$-operad of stable algebraic curves with marked points, recovers the structure of genus 0 Gromov--Witten invariants of a smooth projective complex variety\footnote{The restriction to genus $0$ is only necessary because the higher genus situation would require to develop the brane action for modular $\infty$-operads. Joint work of the author and David Kern aims at generalizing this construction.}. Their construction therefore lifts Gromov--Witten invariants from cohomology or K-theory to a purely geometric level, which then induce categorical invariants, by considering for example the derived categories of quasicoherent sheaves.

The main examples of coherent $\infty$-operads come from geometric structures of "configuration type": the prototypical example is the $\infty$-operad of little disks $\bE_n$. 
Similarly to the case of Gromov--Witten invariants, one expects that the brane action applied to the little disks $\infty$-operad $\bE_2$ can recover the classical construction of string topology operations on the chains of the free loop space $\Map(S^1,X)$ of a closed oriented manifold $X$. The connection is as follows: the cospans of spaces of extensions (of diagram \eqref{eqn:cospan_diagram}) constructed by the brane action for $\bE_2$ are equivalent to the classical cobordism-type diagrams that induce string topology operations upon taking pull-push on homology using Gysin morphisms.

\paragraph*{Brane topology.}
More generally, one may consider spaces of branes  $\Map(S^{n-1},X)$, which are higher dimensional analogs of free loop spaces and whose study is called \emph{brane topology}. A central conjecture\footnote{This conjecture appears explicitly in the introduction of \cite{ginot-tradler-zeinalian_higherHH} and is based on \cite[Section 5.4]{cohen-voronov}.} states that the singular chains on higher loop spaces should inherit an algebra structure over the $\infty$-operad $\bE_{n}^{\mathrm{fr}}$ of \emph{framed} little disks of dimension $n$. This is motivated by work of Sullivan and Voronov \cite{cohen-voronov} which considers the homology version of this conjecture.
The case of the underlying $\bE_{n}$-algebra structure has been proven by Ginot--Tradler--Zeinalian in \cite{ginot-tradler-zeinalian_higherHH}, under the assumption that $X$ is an $(n-1)$-connected Poincaré duality space, using Hochschild homology models (see also \cite{hu_2006} for related results).
Their algebraic methods can be refined to give an $\bE_{n-1}^\mathrm{fr}\o \bE_1$-algebra structure (see \cite[Corollary 5.25]{ginot_hodge_filtration}, but do not readily generalize to the fully-framed situation; one might therefore try to use the brane action instead. However, since the $\infty$-operad $\bE_{n}^\mathrm{fr}$ is not reduced, having the non-trivial group $\bE_n^\mathrm{fr}(1)\simeq \mathrm{SO}(n)$ as unary operations, one cannot use Toën's or Mann--Robalo's theorems.
This motivates the search for a construction of brane actions for non-reduced coherent $\infty$-operads.

\subsection*{Main results} 
\label{sub:statement_of_the_results}

Our contribution is three-fold: we generalize the brane action to {general} coherent $\infty$-operads, complete the proof of the previous approach and finally obtain operations on spaces of branes in more general settings, both topological and derived algebro-geometric. Further applications to brane topology, in particular concerning the conjectural $\bE_{n+1}^{\mathrm{fr}}$-algebra structure on singular chains on spaces of branes, will be developed in a sequel to this paper.

\paragraph*{Extension of the brane action to general coherent \texorpdfstring{$\infty$}{infinity}-operads.}
The first result of this paper generalizes the construction of brane actions to the general setting of coherent $\infty$-operads,  showing that  the reducedness and monochromaticity assumptions from Toën's and Mann--Robalo's results \cite{toen_branes,robalo} are actually unnecessary.

\begin{thmintro}
    \label{thm:main_thm}
    Let $\zO^\o$ be a coherent $\infty$-operad. Then the collection of spaces $\Ext(\id_X)$, for varying colors $X\in \zO$, carries a canonical $\zO$-algebra structure in $\Cospan(\zS)$, with structural maps generalizing the cospan diagrams \eqref{eqn:cospan_diagram}.
\end{thmintro}

Our approach is based on Mann--Robalo's construction of the brane action. It is encapsulated as an explicit cartesian fibration, which we call the \emph{brane fibration}, 
$$\pi\colon \bo\too \Tw(\Env(\zO))^\o$$ 
over the twisted arrow $\infty$-category of the symmetric monoidal envelope of $\zO^\o$,  whose classifying functor gives the desired $\zO$-algebra structure in cospans of spaces (see section \ref{sec:construction_of_the_brane_fibration} for details).

This approach relies on identifications between three possible models for the space of extensions of an operation $\sigma$ of arity $n$: first, the fibers $\bo_\sigma$ of the brane fibration, second, Lurie's space $\Ext(\sigma)$ defined in \cite{ha}\footnote{up to a very slight correction, see remark \ref{rk:difference_df_Ext_with_HA}.} and third, in the monochromatic case, the homotopy fiber at $\sigma$ of the forgetful map $\zO(n+1)\to \zO(n)$, as considered in Toën's original work \cite{toen_branes}. 
We call these spaces respectively Mann--Robalo, Lurie and Toën's model for spaces of extensions.
Identifications between these models are used in existing literature, but left unproven.\footnote{Note that another statement generalizing Toën and Mann--Robalo's results appears in Kern's thesis as \cite[Theorem 2.2.2.0.1]{kern_these}, under an assumption called \emph{hapaxunitality}, which allows for \emph{certain} coloured non-reduced $\infty$-operads. This is analogous to our theorem \ref{thm:main_thm}, but with slightly different generality: on the one hand, the hapaxunitality condition is more restrictive than our assumptions, on the other hand, Kern's results apply also to examples of non-coherent $\infty$-operads, in which case the brane action becomes a lax action. More importantly, as the proof of Kern's result is based on Mann--Robalo's construction, its validity relies on our theorems \ref{thm:comparaison_bo_ext_intro} and \ref{thm:quotient_par_O1_intro}.}

Our next two results fill this gap, by providing an equivalence between the first two models (for general unital $\infty$-operads, see theorem \ref{thm:comparaison_bo_ext_intro}) and the last two models (for \emph{reduced} unital $\infty$-operads, as the comparison has to be modified in the general case, see theorem \ref{thm:quotient_par_O1_intro}).

\paragraph*{Comparison of models of spaces of extensions.}
When trying to identify the two definitions $\bo_\sigma$ and $\Ext(\sigma)$ for the space of extensions of an operation $\sigma$, it seems that no straightforward comparison is available. 
For instance, there is no obvious  morphism of simplicial sets relating the two models that would furnish a candidate for the equivalence. 
We prove the following comparison result by providing an explicit, ad-hoc zigzag of homotopy equivalences between $\bo_\sigma$ and $\Ext(\sigma)$, 

\begin{thmintro}[Theorem \ref{thm:comparaison_bo_ext}]
    \label{thm:comparaison_bo_ext_intro}
    Let $\sigma$ be an active morphism in a unital $\infty$-operad $\zO^\o$. Then the fiber $\bo_\sigma$ of the brane fibration and the $\infty$-category of extensions $\Ext(\sigma)$ are equivalent.
\end{thmintro}

\paragraph*{Computation of homotopy type of spaces of extensions.}
To apply brane actions to particular examples of coherent $\infty$-operads, or to prove that a given $\infty$-operad is coherent, one needs to compute spaces of extensions. For that purpose, Mann--Robalo's model $\bo_\sigma$ and Lurie's $\Ext(\sigma)$ are essentially of no use: for example, even in the simple case of the $\infty$-operad of littke disks $\bE_2$, identifiying the homotopy type of the space of extensions $\Ext(\id)$ seems a difficult combinatorial problem.

Actually, in practice all known computations in the literature involving spaces of extensions rely on the equivalence with {Toën's model} described above.
Such an equivalence for Lurie's model $\Ext(\sigma)$ (and therefore also for Mann--Robalo's model, by theorem \ref{thm:comparaison_bo_ext_intro}) is claimed in \cite[Section 5.1.1]{ha}. More precisely, a comparison map is defined and asserted to be an equivalence. 
However, we find that Lurie's model $\Ext(\sigma)$ only agrees with Toën's for \emph{reduced} $\infty$-operads. In fact, we provide a counter-example (see proposition \ref{prop:example_computation_rmk_HA_fails}) when this assumption fails, thereby contradicting the corresponding statement in \cite{ha}. The reason comes from a non-trivial action of the $\infty$-group of unary operations on the last color of extensions (we refer to section \ref{sub:difference_existing_literature} for a more detailed discussion).

The general, unreduced situation is clarified by the following result, which exhibits $\Ext(\sigma)$ as a quotient of Toën's model by this $\zO(1)$-action. 

\begin{thmintro}[Theorem \ref{thm:quotient_par_O(1)}]
	\label{thm:quotient_par_O1_intro}
	Let $\zO^\o$ be a monochromatic unital $\infty$-operad whose underlying $\infty$-category $\zO$ is an $\infty$-groupoid and let $\sigma\in \zO(n)$ an operation of arity $n$. Choose a semi-inert morphism $i\colon \langle n \rangle \to \langle n+1\rangle$ in $\zO^\o$. Then the space $\Ext(\sigma)$ is equivalent to the homotopy quotient of $\zO(n+1)\mathop{\times}^{\mathrm{h}}_{\zO(n)}\{\sigma\}$ by the natural action of the $\infty$-group $\zO(1)$ of unary operations on the additional color of the extensions.
\end{thmintro}


\paragraph*{Operations on spaces of branes.}
Recall that the $\infty$-operad of little disks $\bE_{n+1}^\o$ is coherent for any $n\geq 0$, by Lurie's result \cite[Theorem 5.1.1.1]{ha} (whose proof relies on the validity of our theorem \ref{thm:quotient_par_O1_intro}).
Using the computation tool given by the previous theorem, we extend this coherence result to the variants $\bE_B^\o$ of $\bE_{n+1}^\o$  of $B$-framed little disks, in the sense of \cite{ha} (recalled in section \ref{sub:coherence_of_variants_of_the_little_disks_operad}). This family of examples include in particular the $\infty$-operad of framed little disks $\bE_n^{\mathrm{fr}}\simeq \bE_n \rtimes SO(n)$, and more generally semi-direct products $\bE_n\rtimes G$ with a topological group $G$ equipped with a morphism to the group of self-homeomorphisms of  $\bR^n$.

\begin{thmintro}[Theorem \ref{thm:B-framed_little_disks_coherent_operad}]
    \label{thmintro:B_framed_little_disks_coherent}
    The $\infty$-operads of $B$-framed little disks $\bE_B^\o$ are coherent.
\end{thmintro}
One can prove that the space of extensions $\Ext(\id_b)$ of any color $b\in B$ is homotopy equivalent to the sphere $S^{n}$.
As a consequence of theorems \ref{thmintro:B_framed_little_disks_coherent} and \ref{thm:main_thm}, we obtain an $\bE_B$-algebra structure on $S^{n}$ in cospans of spaces, and hence corresponding operations on spaces of branes, at the level of spans.

\begin{cor}
    \label{cor:Map_Sn_X_twisted_En_algebra}
    Let $X$ be a topological space. Then the space of branes $\Map(S^n,X)$ has an $\bE_B$-algebra structure in $\Span(\zS)$ given by the brane action.
\end{cor}

In the particular case of the $\infty$-operad $\bE_{n+1}^{\mathrm{fr}}$, this structure at the level of spans is a first step towards the conjectural extension of Sullivan--Voronov brane topology operations to the chain level of $\Map(S^n,X)$.

Finally, using the notion of perfect stacks from \cite{integral_transforms} and the universal property of higher categories of spans \cite{stefanich_higher_sheaf_theory}, we can invert the wrong-way map in spans and obtain brane topology structures in the setting of derived algebraic geometry.
\begin{cor}
    [Corollary \ref{cor:E_B_algebra_inverting_spans}]
    Let $X$ be a perfect stack. Then the $\infty$-category of quasicoherent sheaves on its space of branes $\Map((S^{n-1})_{\mathrm{cst}},X)$ carries a canonical $\bE_B$-algebra (in particular $\bE_n^{\mathrm{fr}}$-algebra) structure in dg-categories. 
\end{cor}

\paragraph*{Organization of the paper.} 
We start in section \ref{sec:the_brane_fibration} by recalling some important constructions: the $\infty$-categories of spans and that of twisted arrows, the precise definition of Lurie's model $\Ext(\sigma)$ for the space of extensions and the definition of coherent $\infty$-operads. We then define the brane fibration, following Mann--Robalo, and outline the proof of theorem \ref{thm:main_thm}.
This proof is then completed in section \ref{sec:cartesianity_of_brane_fibration}, by establishing that the functor  $\pi\colon \bo\to \Tw(\Env(\zO))^\o$ is indeed a cartesian fibration (theorem \ref{thm:pi_is_cartesian_fibration}).

Section \ref{sec:comparison_with_lurie_s_spaces_of_extensions} is devoted to the proof of theorem \ref{thm:comparaison_bo_ext_intro}, that is the comparison between Mann--Robalo's and Lurie's model of spaces of extensions, via the construction of an explicit zigzag of homotopy equivalences. 
Next, we deal in section \ref{sec:homotopy_type_of_spaces_of_extensions} with the problem of computing the homotopy type of spaces of extensions, by establishing an equivalence between Toën's and Lurie's models, thereby proving theorem \ref{thm:quotient_par_O1_intro}. 
Moreover, we discuss how our results differ from a claim in \cite{ha} and provide a counterexample to the latter statement. 

We end in section \ref{sec:applications} with applications to brane topology, starting with a proof of coherence of the $\infty$-operad of $B$-framed little disks (theorem \ref{thmintro:B_framed_little_disks_coherent}).
We conclude with a construction of new operations on spaces of branes, both in the topological context (at the span level) and for derived algebraic stacks (at the level of derived categories).


Finally, we provide an appendix, composed of two parts. The first contains definitions and results concerning marked simplicial sets that are crucial to the arguments of section \ref{sec:comparison_with_lurie_s_spaces_of_extensions}; they do not appear in the literature and might be of independent interest. The second part gathers some results from the literature, mainly on $\infty$-group actions, and is provided merely for the reader's convenience.





\paragraph*{Notations and conventions.} 
\label{sub:notations}
We work in the particular model of $\infty$-category theory given by quasicategories and use Lurie\rq s presentation of $\infty$-operads. Our notations generally follow those of \cite{htt} and \cite{ha}.

\begin{itemize}
	\item 
        {Particular arrows:} monomorphisms are denoted as $\begin{tikzcd}
			[cramped, sep=small] A\ar[r,hook] &B
		\end{tikzcd}$,
		cofibrations as $\begin{tikzcd}[cramped, sep = small]A \ar[r,tail] &B\end{tikzcd}$
		and atomic morphisms (see definition \ref{df:semi-inert_atomic_maps}) as $\begin{tikzcd}
			[cramped, sep = small] A\ar[r,hook,harpoon]& B
		\end{tikzcd}$.
		\todo{donner les conventions sur les flèches : cofibration, atomic, act, etc. et les appliquer partout}
	\item
		When considering a diagram $X\colon P\to \zC$ from a poset $P$ to an $\infty$-category $\zC$ and a sequence $i_0 \leq i_1 \leq \dots \leq i_n$ in $P$, we will write $X_{i_0}\dots X_{i_n}$ for the $n$-simplex of $X\circ \langle i_0\dots i_n\rangle \colon \Delta^n \to P \to \zC$.
		For instance, the notation $X_iX_j$ denotes the unique morphism $X_i\to X_j$ of the diagram.
	\item
		Given a finite linear order $I = \{i_0 < i_1 < \dots < i_n\}$, the full subsimplex of $\Delta^{I}$ on the objects $i_{j_0} < \dots < i_{j_k}$ will be denoted $\Delta^{i_{j_0}\dots i_{j_k}}$
		(unless $k=0$). Similarly, $\Lambda_{i_{j_p}}^{i_{j_0}\dots i_{j_k}}$ stands for the horn in $\Delta^{i_{j_0}\dots i_{j_k}}$ obtained by removing the face opposed to vertex $i_{j_p}$. 
		For instance, the horns $\Lambda^{12}_1$ and $\Lambda^{12}_2$ are respectively the simplicial subsets $\Delta^{\{2\}}$ and $\Delta^{\{1\}}$ of the $1$-simplex $\Delta^{12}$, while the notation $\Lambda^{12}_0$ does not make sense in our convention.
    \item 
        For simplicity, given an $\infty$-operad $\zO^\o$, we will often write $\zE$ for its symmetric monoidal envelope $\Env(\zO)^\o$ and $\zT$ for the associated twisted arrows $\infty$-category $\Tw(\Env(\zO))^\o$ (see notation \ref{nota:E_T}).
    \item We let $\bF_*$ denote the nerve of the category of pointed finite sets. We usually identify $\bF_*$ with its equivalent full subcategory on the pointed sets $\langle n \rangle = (\{0,\dots,n\},0)$.
\end{itemize}

\paragraph*{Acknowledgements.}
I would like to thank my advisor Grégory Ginot for his time, guidance and constant support during my PhD. 
I am also grateful to Marco Robalo for numerous useful discussions regarding this work.

This project has received funding from the European Union’s Horizon 2020 research and innovation programme under the Marie Skłodowska-Curie grant agreement No 754362 \includegraphics[height=0.8em]{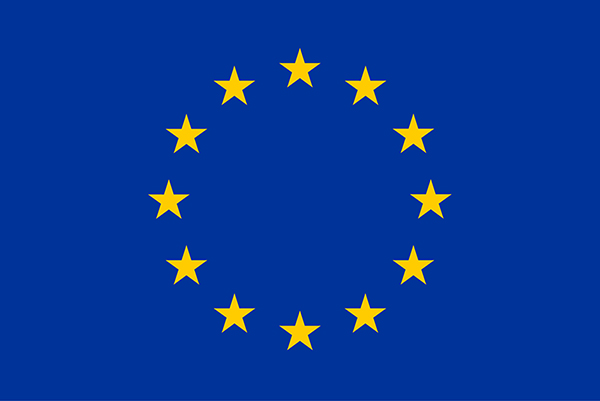}.


	\section{The brane fibration} 
\label{sec:the_brane_fibration}

In this section, following \cite{robalo}, we explain how the brane action of theorem \ref{thm:main_thm} arises from a certain fibration, which we call \emph{the brane fibration}. 
Before giving the precise construction, we recall the notions of $\infty$-categories of (co)spans, of twisted arrows, of spaces of extensions in the sense Lurie and the definition of coherent $\infty$-operads.


\subsection{Categories of spans and of twisted arrows} 
\label{sec:categories_of_spans_and_of_twisted_arrows}

Given an $\infty$-category $\zC$ with finite limits, we may form the $\infty$-category $\Span(\zC)$ of spans in $\zC$, whose objects are those of $\zC$, morphisms between two objects $X$ and $Y$ are given by span diagrams
$\begin{tikzcd}[cramped, sep = small]
	X & Z\ar[l]\ar[r] & Y
\end{tikzcd}$
and composition is given by taking pullback (see \cite{barwick_Q-construction_for_exact_quasicategories} or \cite{haugseng_spans} for a rigorous $\infty$-categorical definition). Dually, if $\zC$ has finite colimits, we may consider its $\infty$-category of cospans $\Cospan(\zC)$ defined as $\Span(\zC\op)$.

\sloppy
The $\infty$-category $\Span(\zC)$ has a canonical symmetric monoidal structure $\Span(\zC)^{\o_\times}$ induced from the cartesian monoidal structure on $\zC^\times$, although $\Span(\zC)^{\o_\times}$ is not itself cartesian.

\begin{df}[Twisted arrow $\infty$-category]
	\label{df:Tw}
	Let $s\colon \mathsf{\Delta}\to\mathsf{\Delta}$ be the functor given by $s[n] = [n]\ast[n]\op$. Precomposition with $s$ yields an endofunctor $s^*\colon \sSet\to\sSet$ that we shall denote $\Tw$. 
	Left Kan extension of $s$ along the Yoneda embedding of $\mathsf{\Delta}$ induces a functor $s_*\colon \sSet\to \sSet$ left adjoint to $\Tw$.
The image under $\Tw$ of an $\infty$-category $\zC$ is again an $\infty$-category $\Tw(\zC)$ called its \emph{twisted arrow} $\infty$-category, whose $n$-simplices are $(2n+1)$-simplices $s_*(\Delta^n)= \Delta^n\ast \Delta^{n,\mathrm{op}}\to \zC$, represented as
\[
\begin{tikzcd}
X_0 \arrow[r] \arrow[d] & X_1 \arrow[r] \arrow[d] & \dots \arrow[r] & X_n \arrow[d]       \\
\bar{X_0}               & \bar{X_1} \arrow[l]     & \dots \arrow[l] & \bar{X_n}. \arrow[l]
\end{tikzcd}
\]
To depict a morphism in $\Tw(\zC)$, that is, a twisted arrow between two arrows $f$ and $g$ of $\zC$, we will often write $f\rightsquigarrow g$.
\end{df}

\begin{rk}
	Given a $2$-simplex $\sigma$
	\[
	\begin{tikzcd}
	& Y \ar[rd,"g"] & \\
	X \ar[ru,"f"] \ar[rr,"h"] & & Z
	\end{tikzcd}
	\]
	in $\zC$ that exhibits $h$ as a composite of $g$ and $f$, we obtain twisted arrows $h\rightsquigarrow g$ and $h\rightsquigarrow f$ in $\Tw(\zC)$ given respectively by the following commutative squares:
	\[
	\begin{tikzcd}  
		X \ar[r,"h"] \ar[rd,"h"]\ar[d,"f" left]	&	Z\\
		Y \ar[r,"g"']					& 	Z \ar[u,"\id_Z"'] 
	\end{tikzcd}
	\qquad \text{and}\qquad
	\begin{tikzcd}  
		X \ar[r,"h"] \ar[rd,"f"]\ar[d,"\id_X" left]	&	Z \\
		X \ar[r,"f"']					& 	Y \ar[u,"g"'] 
	\end{tikzcd}
	\]
	in which the $2$-simplices are either degenerate, are equal to $\sigma$.
\end{rk}

By \cite[Example 5.2.2.23.]{ha}, any symmetric monoidal $\infty$-category $\zC^\o$ induces a symmetric monoidal structure $\Tw(\zC)^\o$ on the twisted arrow $\infty$-category $\Tw(\zC)$, in which the tensor product of two morphisms $f\colon x\to y$ and $g\colon z\to t$ is the obvious arrow of the form $f\otimes g \colon x\otimes z\to y\otimes t$.

An important feature of the construction of $\infty$-category of twisted arrows is the following universal property.

\begin{prop}[Universal property of $\Tw$ and $\Span$]
	\label{prop:universal_property_Tw_Span}
	\begin{enumerate}
		Let $\zC$ and $\zD$ be two $\infty$-categories and assume that $\zD$ has all finite limits. Then:
		\item 
			There is a natural equivalence between the space of functors $\zC \to \Span(\zD)$ and that of functors $F\colon \Tw(\zC)\to \zD$ satisfying the \emph{pullback condition}: namely that for every $2$-simplex $h\colon X\stackrel{f}{\to} Y\stackrel{g}{\to} Z$ exhibiting $h$ as a composite of $g$ and $f$, the induced square
			\[
			\begin{tikzcd}  
				h \ar[r,squiggly] \ar[d,squiggly]	&	g \ar[d,squiggly] \\
				f \ar[r,squiggly]					& 	\id_{Y}
			\end{tikzcd}
			\]
			in $\Tw(\zC)$ is sent by $F$ to a cartesian square in $\zD$.
		\item 
            \sloppy
			If $\zC^\o$ is a symmetric monoidal $\infty$-category with underlying $\infty$-category $\zC$, then there is a natural equivalence between the space of symmetric monoidal functors $\zC^\o \to \Span(\zD)^{\o_\times}$ and that of symmetric monoidal functors $\Tw(\zC)^\o \to \zD^\times$ satisfying the above pullback condition.
	\end{enumerate}
\end{prop}
\todo{[DONE] vérifier un peu ça et ajouter des ref}

A proof of the first part of this result can be found in the appendix of \cite[Theorem 2.18.]{haugseng-hebestreit-linskens-nuiten}: there, the statement takes the stronger form of an adjunction between $\Cat_\infty$ and a certain $\infty$-category $\Cat_\infty^\mathrm{dir}$ of small $\infty$-categories with directions, which fully-faithfully contains the $\infty$-category of small $\infty$-categories with finite limits and functors preserving them. In particular, this requires to enhance $\Tw$ to a functor $\Cat_\infty \to \Cat_\infty^\mathrm{dir}$. The extension to the symmetric monoidal case is explained in \cite[Corollary 2.1.3.]{robalo}.


\subsection{Extensions and coherent \texorpdfstring{$\infty$}{infinity}-operads} 
\label{sec:extensions_and_coherent_operads}

In this subsection, we recall the definition of the $\infty$-category of extensions of an operation in an $\infty$-operad and the closely related notion of coherence. We essentially follow \cite[Section 3.3.1]{ha}, except for a small difference in the definition of $\Ext(\sigma)$ (see remark \ref{rk:difference_df_Ext_with_HA}).

Let $p\colon \zO^\o \to \Fin_*$ be a unital $\infty$-operad.

\begin{df}[Semi-inert and atomic maps]
	\label{df:semi-inert_atomic_maps}
	Let  $f\colon X\to Y$ be a morphism in $\zO^\o$, corresponding to a morphism $\alpha = p(f) \colon \langle n\rangle \to \langle m\rangle$ in $\Fin_*$ together with a family of multimorphisms $f_j\colon \{X_i\}_{\alpha(i)=j} \to Y_j$ for $j\in \langle m\rangle^\circ$.
	We say that $f$ is \emph{semi-inert} if for every $j\in \langle m\rangle^\circ$
	\begin{itemize}
		\item either the set $\alpha\inv\{j\}$ is empty, or
		\item the set $\alpha\inv\{j\}$ is the singleton $\{i_j\}$ and the map $f_j\colon X_{i_j}\to Y_j$ is an equivalence.
	\end{itemize}
    Following the terminology of \cite{kern_these}, we say that $f$ is \emph{atomic} if it is semi-inert and lies over an inclusion $\alpha\colon \langle n\rangle \to \langle n+1\rangle$. In other words, $f$ is atomic if and only if it is semi-inert with no non-trivial factorization through another semi-inert morphism.
	Given a commutative diagram
	\[
	\begin{tikzcd}  
		X \ar[r,""] \ar[d,"f_X" left, harpoon, hook]	&	Y \ar[d,"f_Y", harpoon, hook] \\
		X' \ar[r,"f"]					& Y' 	
	\end{tikzcd}
	\qquad \text{or} \qquad
	\begin{tikzcd}[column sep = small]
		& X\ar[ld,"f_X" left, harpoon, hook]\ar[rd,"f_Y", harpoon, hook] & \\
		X' \ar[rr,"f"] & & Y'
	\end{tikzcd}
	\]
	with $f_X$ and $f_Y$ atomic, we say that $f$ is \emph{compatible with extension} if $f$ sends the unique color of $p(X')\setminus \im(p(f_X))$ to the unique color of $p(Y')\setminus \im(p(f_Y))$.
\end{df}

\begin{rk}
	\label{rk:comparison_df_atomic_m-semi-inert_HA}
	In \cite[Definition 3.3.2.3.]{ha}, the notion of a \emph{$m$-semi-inert morphism} is introduced, for $m\in \bN$. In terms of this definition, a morphism $f$ in $\zO^\o$ is atomic if and only if it is $1$-semi-inert but not $0$-semi-inert.
\end{rk}

\begin{df}[$\infty$-category of extensions]
	\label{df:Ext(sigma)}
	Let $\sigma\colon \Delta^n\to \oact$ be an $n$-simplex corresponding to a sequence of active morphisms $X_0\stackrel{f_1}{\to} \dots \stackrel{f_n}{\to} X_n $. Given a downward-closed subset $S = \{0,\dots,r\}\subseteq [n]$, let $\Ext(\sigma,S)$ be the (non-full) subcategory of $\Fun(\Delta^n,\zO^\o)_{\sigma/}$ whose
	\begin{itemize}
		\item \textbf{objects} are diagrams $\Delta^1\times\Delta^n\to \zO^\o$ represented as
		\[
		\begin{tikzcd}
		X_0 \ar[r,"f_1"]\ar[d,"g_0",hook,harpoon] & \dots \ar[r,"f_{r}"] & X_{r} \ar[r]\ar[d,"g_{r}",hook,harpoon]& X_{r+1} \ar[r]\ar[d,"g_{r+1}"right, "\wr"left] & \dots \ar[r,"f_n"] & X_n\ar[d,"g_n"right, "\wr"left]\\
		X'_0 \ar[r,"f'_1"] & \dots \ar[r,"f'_{r}"] & X'_{r} \ar[r] & X'_{r+1} \ar[r] & \dots \ar[r,"f'_n"] & X'_n
		\end{tikzcd}
		\]
		satisfying the following conditions:
		\begin{enumerate}
			\item if $i\notin S$, then $g_i$ is an equivalence,
			\item if $i \in S$, then $g_i$ is atomic,
			\item if $i \in S\setminus \{0\}$, then $f_i'$ is compatible with extension,
			\item each map $f_i'$ is active;
		\end{enumerate}
		\item \textbf{morphisms} are diagrams $\Delta^2\times\Delta^n\to \zO^\o$ represented as 
		\[
		\begin{tikzcd}
		X_0 \ar[r,"f_1"]\ar[d,"g_0",hook,harpoon]\ar[dd,hook,harpoon,bend right=40]  & \dots \ar[r,"f_n"] & X_n\ar[d,"g_n"left, "\wr"right]\ar[dd,bend left=40,"\wr"right]\\
		X'_0 \ar[r,"f'_1"]\ar[d,"h_0"] & \dots \ar[r,"f'_n"] & X'_n\ar[d,"h_n"left]\\
		X''_0 \ar[r,"f''_1"] & \dots \ar[r,"f''_n"] & X''_n
		\end{tikzcd}
		\]
		in which the morphisms $h_i\colon X_i'\to X''_i$ are compatible with extension for all $i\in S$.
	\end{itemize}
	Given an active morphism $\sigma\colon \Delta^1\to \oact$, we write $\Ext(\sigma)$ for $\Ext(\sigma,\{0\})$. We call $\Ext(\sigma)$ the \emph{$\infty$-category of extensions} of $\sigma$. When the underlying $\infty$-category $\zO$ of $\zO^\o$ is an $\infty$-groupoid, $\Ext(\sigma)$ is a Kan complex and therefore refered to as the \emph{space} of extensions of $\sigma$.
\end{df}

\begin{ex}[Description of $\Ext(\sigma)$ in the discrete case]
    \label{ex:description_simplices_Ext_discrete_case}
    Let $\zO_\Delta$ be an operad in sets, $\zO^\o$ its homotopy coherent nerve and $\sigma\colon \langle m\rangle \to \langle 1\rangle$ an active morphism in $\zO^\o$. 
    Then the $k$-simplices of $\Ext(\sigma)$ are those functors between $1$-categories $[1]\times [1+k]\to \zO_\Delta$ whose associated diagrams is of the form
	\begin{equation}
        \label{eqn:description_smplices_Ext_discrete_case}
        \begin{tikzcd}
            \langle m\rangle \arrow[r, "\sigma"] \arrow[d, "\text{atomic}", hook, harpoon]                   & \langle 1\rangle \arrow[d, "\sim"]               \\
            \langle m+1\rangle \arrow[d, "\vdots" description, no head, dotted] \arrow[r, "\text{act}"] & \langle 1\rangle \arrow[d, "\vdots" description] \\
            \langle m+1\rangle \arrow[d] \arrow[r, "\text{act}"]                                        & \langle 1\rangle \arrow[d, "\sim"]               \\
            \langle m+1\rangle \arrow[r, "\text{act}"]                                                  & \langle 1\rangle.
        \end{tikzcd}
    \end{equation}
    and such that all the left vertical morphism $\langle m+1\rangle \to \langle m+1\rangle$ are compatible with extensions.
\end{ex}

\begin{rk}[Difference with the existing definition]
	\label{rk:difference_df_Ext_with_HA}
	The previous definition is slightly different from the initial definition from \cite{ha} in that we impose a condition on the morphisms in $\Ext(\sigma)$, rather than defining it as a \emph{full} subcategory of $\Fun(\Delta^n,\zO^\o)_{\sigma/}$. The reason for this choice is that the space defined in \cite[Definition 3.3.1.4.]{ha}, that we shall denote $\Ext^\mathrm{HA}(\sigma)$ here, does not have the expected homotopy type. To see this, consider the example of the commutative $\infty$-operad $\zO^\o = \mathrm{Comm}^\o$ and $\sigma\colon \langle m \rangle \to \langle 1\rangle$ be an active map in $\mathrm{Comm}^\o$. As described in \cite[Example 3.3.1.12]{ha}, the space of extensions of $\sigma$ is supposed to be the singleton set $\langle 1\rangle^\circ$, viewed as a discrete space. However, the space $\Ext^\mathrm{HA}(\sigma)$ is not discrete.
	Indeed, consider the object $\alpha \in \Ext^\mathrm{HA}(\sigma)$ given by the following diagram
		\[
\begin{tikzcd}
\langle m\rangle \arrow[r, "\sigma"] \arrow[d, "i", hook, harpoon]                   & \langle 1\rangle \arrow[d, "\id"]               \\
\langle m+1\rangle \arrow[r, "!"] & \langle 1\rangle
\end{tikzcd}
	\]
	where $m$ is a positive integer, $!\colon \langle m+1\rangle \to \langle 1\rangle$ is the unique active map and $i$ the canonical inclusion. We claim that $\pi_1(\Ext^\mathrm{HA}(\sigma),\alpha)$ is not trivial. Let $\mu\colon \langle m+1 \rangle \to \langle m+1\rangle$ be the morphism in $\mathrm{Comm}_\mathrm{act}^\o$ that restricts to the atomic morphism $i$ on $\langle m\rangle$ and sends the remaining color $m+1$ to $1$. Then the diagram
	\begin{equation}
		\label{eqn:diag_counter_example_df_Ext_HA}
		\begin{tikzcd}
		\langle m\rangle \arrow[r, "\sigma "] \arrow[d, hook,harpoon,"i"] \arrow[dd, hook,harpoon,"i"left, bend right=50] & \langle 1\rangle \arrow[d, "\id"] \arrow[dd, "\id", bend left=40] \\
		\langle m+1\rangle \arrow[d, "\mu"] \arrow[r, "!"]                                 & \langle 1\rangle \arrow[d, "\id"]                              \\
		\langle m+1\rangle \arrow[r, "!"]                                                  & \langle 1\rangle                                              
		\end{tikzcd}
	\end{equation}
	defines a morphism $\gamma \colon \alpha \to \alpha$ in $\Ext^\mathrm{HA}(\sigma)$ with the property that $[\gamma]\neq [\id_\alpha]$ in $\pi_1(\Ext^\mathrm{HA}(\sigma),\alpha)$. Indeed, a homotopy between $\gamma$ and $\id_\alpha$ would give a retraction $\rho$ of $\mu$, which can't be.

	Note that diagram \eqref{eqn:diag_counter_example_df_Ext_HA} does not define a morphism in $\Ext(\sigma)$ since $\mu$ is not compatible with extensions. We will see that definition \ref{df:Ext(sigma)} yields the expected homotopy type for the spaces of extensions: this is the content of theorem \ref{thm:quotient_par_O1_intro}.
\end{rk}

\begin{df}[{\cite[Definition 3.3.1.9]{ha}}]
	\label{df:coherence}
	An $\infty$-operad $\zO^\o$ is \emph{coherent} if it satisfies the following conditions:
	\begin{enumerate}[label=(\alph*),ref=(\alph*)]
	 	\item \label{item:df_coherence_a}
			it is unital, 
		\item \label{item:df_coherence_b}
			its underlying $\infty$-category $\zO$ is an $\infty$-groupoid,
		\item \label{item:df_coherence_c}
			for every degenerate $3$-simplex $\sigma$
			\[
			\begin{tikzcd}
			                                  & Y \arrow[rd, "\id_Y"] \arrow[rr, "g"] &                   & Z \\
			X \arrow[ru, "f"] \arrow[rr, "f"] &                                       & Y \arrow[ru, "g"] &  
			\end{tikzcd}
			\]
			in $\oact$, the commutative diagram
			\begin{equation}
				\label{eqn:diag_df_coherence}
				\begin{tikzcd}  
					\Ext(\sigma, \{0,1\}) \ar[r,""] \ar[d,"" left]	&	\Ext(\sigma|_{\Delta^{\{0,1,3\}}}, \{0,1\}) \ar[d,""] \\
					\Ext(\sigma|_{\Delta^{\{0,2,3\}}}, \{0\}) \ar[r,""]					& 	\Ext(\sigma|_{\Delta^{\{0,3\}}}, \{0\})
				\end{tikzcd}
			\end{equation}	
			is a homotopy cocartesian square of Kan complexes.
	\end{enumerate} 
\end{df}

\begin{rk}
	Let $\sigma$ and $S = \{0,\dots, r\}$ be as in definition \ref{df:Ext(sigma)} and suppose that $\zO$ is an $\infty$-groupoid. As mentionned before, the simplicial set $\Ext(\sigma)$ is a Kan complex. By remark \cite[Remark 3.3.1.6.]{ha}, if $r < [n]$, there is a canonical map $\Ext(\sigma, S)\to \Ext(f_{r+1})$ which is trivial Kan fibration. Using these equivalences, we may rewrite the commutative square \eqref{eqn:diag_df_coherence} as
	\begin{equation}
        \label{eqn:diag_df_coherence_rewritten}
		\begin{tikzcd}  
			\Ext(\id_Y) \ar[r,""] \ar[d,"" left]	&	\Ext(g) \ar[d,""] \\
			\Ext(f) \ar[r,""]					& 	\Ext(g\circ f).
		\end{tikzcd}
	\end{equation}		
\end{rk}
Note that the previous square is only well-defined in the homotopy category of spaces.



\subsection{Symmetric monoidal envelope and its twisted arrows}
\label{sec:symmetric_monoidal_envelope_and_its_twisted_arrows}

Recall the construction of the symmetric monoidal envelope $\Env\col \Op_\infty \to \Cat_\infty^\o$, which is left adjoint to the forgetful functor from symmetric monoidal $\infty$-categories to $\infty$-operads \cite[Section 2.2.4]{ha}. This left adjoint sends an $\infty$-operad $\zP^\o$ to the $\infty$-category
$$ \Env(\zP)^\o := \zP^\o \times_{\Fun(\{0\}, \Fin_*)} \Fun^{\mathrm{act}}(\Delta^1,\Fin_*),$$
where the superscript $\mathrm{act}$ indicates the full subcategory of $\Fun(\Delta^1,\Fin_*)$ whose objects are {active} morphisms in $\Fin_*$. This $\infty$-category inherits a symmetric monoidal structure via the functor $p_1\colon \Env(\zO)^\o \to\Fin_*$ given by evaluation at $1\in \Delta^1$. 
Note that 
the underlying $\infty$-category $\Env(\zO)$ of $\Env(\zO)^\o$ can be identified with the wide subcategory $\oact$ of $\zO^\o$ consisting of all objects and only active maps between them. 
As explained in \ref{sec:categories_of_spans_and_of_twisted_arrows}, the $\infty$-category of twisted arrows $\Tw\left(\Env(\zO)\right)$ inherits a symmetric monoidal structure from that of the monoidal envelope, also denoted $p_1\colon \Tw(\Env(\zO))^\o \to \Fin_*$.

For later purposes, we let $p_0 \colon \Env(\zO)^\o \to \Fin_*$ denote the functor given by evaluation at $0$.

\begin{nota}
    \label{nota:E_T}
	For simplicity, we will write $\zE$ for $\Env(\zO)^\o$ and $\zT$ for $\Tw(\Env(\zO))^\o$.
\end{nota}

Let us unravel the definitions of $\zE$ and $\zT$.

\begin{itemize}
    \item 
An object in $\zE_{\langle n\rangle}$ is given by an object $X\in \zO^\o_{\langle k\rangle}$ together with an active map $\langle k \rangle \to \langle n\rangle$ in $\Fin_*$. 
In terms of the projection functors $p_0$ and $p_1$, we have that 
$p_0(X,\langle k\rangle \to \langle n\rangle) = \langle k \rangle$ and 
$p_1(X,\langle k\rangle \to \langle n\rangle) = \langle n \rangle$.
Thus, we may think of the object $(X, \langle k\rangle \to \langle n\rangle)$ in $\zE$ as a list of $n$ objects $(X_1,\dots,X_n)$ in $\zO^\o$, with total arity $\oplus_{i = 1}^n \, p_0(X_i) \cong \langle k\rangle$. 
    \item 
A morphism $f$ in $\zE$ from $(X,p_0(X) \to \langle n\rangle)$ to $(Y, p_0(Y) \to \langle m\rangle)$ is a morphism $X\to Y$ in $\zO^\o$ together with a commutative diagram
\[
\begin{tikzcd}  
	p_0(X) \ar[r,""] \ar[d,"" left]	&	p_0(Y) \ar[d,""] \\
	\langle n\rangle \ar[r,"\alpha"]					& 	\langle m\rangle.
\end{tikzcd}
\]
In the case where $\alpha$ is active, the morphism $f$ is $p_1$-cocartesian if and only if $X\to Y$ is an equivalence, by \cite[Lemma 2.2.4.15.]{ha}

    \item 
        An object of $\zT_{\langle n\rangle}$ is given by an active map $g\colon X\to Y$ in $\zO^\o$ together with a commutative triangle
    \[
    \begin{tikzcd}[sep = small]
    p_0(X) \ar[rr,"p_0(g)"] \ar[rd,"\mathrm{act}"below] & & p_0(Y) \ar[dl,"\mathrm{act}"]\\
    & \langle n\rangle & 
    \end{tikzcd}
    \]
    of active maps in $\Fin_*$. A morphism in $\zT$ from the previous object to $(g'\colon X'\to Y', p_0(Y')\to \langle m \rangle)$ is given as a pair of commutative diagrams
    \begin{equation}
        \label{eqn:description_morphism_in_T}
    \begin{tikzcd}  
        X \ar[r,""] \ar[d,"g" left]	&	X' \ar[d,"g'"] \\
        Y 					& 	Y' \ar[l]
    \end{tikzcd}
    \qquad \text{and}\qquad
    \begin{tikzcd}  
        p_0(X) \ar[r,""] \ar[d,"p_0(g)" left]	&	p_0(X') \ar[d,"p_0(g')"]  \\
        p_0(Y) 	\ar[d,"\mathrm{act}"]				& 	p_0(Y') \ar[l]  \ar[d,"\mathrm{act}"]  \\
        \langle n\rangle \ar[r,"\alpha"] & \langle m\rangle
    \end{tikzcd}
    \end{equation}
    respectively in $\zO^\o$ and $\Fin_*$.
    If one interprets the objects in $\zE$ as lists of objects of $\zO^\o$, then the equivalence $\zT_{\langle m\rangle}\simeq \Tw(\oact)^m$ allows to view the object $(g\colon X\to Y, p_0(Y)\to \langle n\rangle)$ in $\zT_{\langle n\rangle}$ as a list of $n$ active morphisms $(X_1\to Y_1,\dots,X_n\to Y_n)$ in $\zO^\o$.

    \item 
A morphism in $\zT$ between two objects $(\sigma_1\colon X_1\to Y_1,\dots,\sigma_n \colon X_n\to Y_n)$ and $(\sigma_1'\colon X_1'\to Y_1', \dots, \sigma_m'\colon X_m'\to Y_m')$ then corresponds to a morphism $\alpha \colon \langle n\rangle \to \langle m\rangle$ in $\Fin_*$ together with a commutative diagram in $\zO^\o$ of the form
\[
\begin{tikzcd}  
	\mathop{\bigoplus}\limits_{i\in \alpha\inv(j)}X_i \ar[r,""] \ar[d,"\mathrm{act}" right, "\oplus_{i} \sigma_i" left]	&	X_j' \ar[d,"\mathrm{act}"] \\
	\mathop{\bigoplus}\limits_{i\in \alpha\inv(j)}Y_i 					& 	Y_j' \ar[l]
\end{tikzcd}
\]
for each $j\in \langle m\rangle^\circ$.
\end{itemize}

\begin{rk}
	\label{rk:condition_morphism_in_T_p_1-cocartesian}
	Consider the morphism in $\zT$ given by diagrams \eqref{eqn:description_morphism_in_T} and assume that $\alpha$ is active. Then this morphism is $p_1$-cocartesian if and only if both maps $X\to X'$ and $Y'\to Y$ are equivalences.
\end{rk}


\subsection{Construction of the brane fibration}
\label{sec:construction_of_the_brane_fibration}

To prove theorem \ref{thm:main_thm}, we will follow the strategy developed by Mann--Robalo in \cite[Section 2.1]{robalo}. Let us recall their approach.

\subsubsection{Mann--Robalo's strategy} 
\label{ssub:mann_robalo_s_strategy}
	First, note that the datum of a map of $\infty$-operads $\zO^\o \to \Cospan(\zS)^\o$ is equivalent to that of a map of symmetric monoidal functors $\zE \to \Cospan(\zS)^\o$. By the universal property of spans (proposition \ref{prop:universal_property_Tw_Span}), this datum is equivalently that of a symmetric monoidal functor $\zT \to (\zS\op)^\amalg$ satisfying the pullback condition.
    By \cite[Proposition 2.4.1.7.]{ha}, since the monoidal structure $(\zS\op)^\amalg$ on $\zS\op$ is cartesian, this is the same as providing a weak cartesian structure $\zT\to \zS\op$ satisfying the pullback condition. \todo{[OUBLI POUR L'INSTANT] ajouter la définition de structure cartésienne faible, plus de détails aussi}
	Using the Grothendieck construction, it will suffice to construct a right fibration $\pi\col\bo \to \zT$ whose classifying functor $F_\pi\colon \zT\to \zS\op$ satisfies the conditions described above. 

The rest of this section is devoted to the construction of this fibration $\pi$.

\begin{df}[The brane fibration, following \cite{robalo}]
\label{df:bo}
Define $\bo$ as the subsimplicial set of $\Fun(\Delta^1, \zT)$ whose 
\begin{itemize}
	\item \textbf{objects} are twisted morphisms $\sigma \leadsto \sigma^+$ such that
	\begin{itemize}
		\item 
		the projection $p_1(\sigma\leadsto \sigma^+)$ in $\Fin_*$ is the unique active map $p_1(\sigma) \to \langle 1\rangle$;
		\item in the corresponding $3$-simplex in $\zO^\o$
		\begin{equation}
			\label{eqn:diag_definition_objects_BO}
			\begin{tikzcd}  
				S_0 \ar[r,"\sigma_0",hook,harpoon] \ar[d,"\sigma" left]	&	S_0^+ \ar[d,"\sigma^+"] \\
				S_1 					& 	S_1^+ \ar[l,"\sim" above, "\sigma_1" ]
			\end{tikzcd}
		\end{equation}
	the map $\sigma_0$ is atomic and $\sigma_1$ is an equivalence; 
	\end{itemize}
	\item \textbf{morphisms} from $\sigma\leadsto \sigma^+$ to $\tau\leadsto \tau^+$ are the morphisms $f$ in $\Fun(\Delta^1,\zT)$ such that
	\begin{itemize}
		\item the projection $p_1\left( 
		\begin{tikzcd}  
			\sigma \ar[r,squiggly] \ar[d,squiggly]	&	\sigma^+ \ar[d,squiggly] \\
			\tau \ar[r,squiggly]					& 	\tau^+
		\end{tikzcd}
		\right)$
		in $\Fin_*$	is the diagram
		$
		\begin{tikzcd}  
			p_1(\sigma)  \ar[r, "\mathrm{act}"] \ar[d]	&	\langle 1\rangle \ar[d, "\id"] \\
			p_1(\tau)  \ar[r, "\mathrm{act}"]					& 	\langle 1\rangle,
		\end{tikzcd}
		$
		\item in the induced square
			\begin{equation}
				\label{eqn:diag_definition_bo}
				\begin{tikzcd}  
					S_0 \ar[r,"\sigma_0",hook,harpoon] \ar[d,"f_0" left]	&	S_0^+ \ar[d,"f_0^+"] \\
					T_0 \ar[r,"\tau_0",hook,harpoon]					& 	T_0^+,
				\end{tikzcd}
			\end{equation}	
			the morphism $f_0^+$ is \emph{compatible with extension},
			in the sense that $p_0(f_0^+)$ is of the form $\langle s+1 \rangle \to \langle t+1\rangle$, sending the singleton $\langle s+1 \rangle\setminus \im(p_0(\sigma_0))$ to the singleton $\langle t+1 \rangle\setminus \im(p_0(\tau_0))$.
	\end{itemize}
\end{itemize}
Let $\pi\col \bo \to \zT$ be the composite of $\mathrm{ev}_0$ with the inclusion $\bo \subset \Fun(\Delta^1, \zT)$.
\end{df}

\begin{rk}
	\label{rk:remarks_on_objects_in_bo}
	The following properties will be useful throughout this paper.
	\begin{itemize}
		\item Since equivalences and atomic maps are active, the diagram \eqref{eqn:diag_definition_objects_BO} is in fact in $\oact$.
		\item The image of $\bo$ under $p_1$ is constant along the fibers of $\pi$, in the sense that there is a commutative diagram
		\[
		\begin{tikzcd}  
			\bo \ar[r,hook] \ar[d,"\pi" left]	&	\Fun(\Delta^1,\zT) \ar[dd,"p_1\circ -"] \\
			\zT \ar[d,"p_1" left]					& 	\\
			\Fin_*\ar[r,"\varepsilon"]					& \Fun(\Delta^1,\Fin_*).
		\end{tikzcd}
		\]
		Here $\varepsilon$ is the unique functor that sends $\langle n\rangle$ to the unique active morphism $\langle n\rangle \to \langle 1\rangle$ and such that $\mathrm{ev}_0\circ \varepsilon = \id_{\Fin_*}$ and $\mathrm{ev}_1\circ \varepsilon = \mathrm{const}_{\langle 1\rangle}$. Thanks to this observation, we will often leave implicit the description of the projection under $p_1$ of various constructions.
	\end{itemize}
\end{rk}

Let us mention the following general facts about $\bo$.

\begin{lm}
\label{lm:pi_inner_fibration}
The inclusion $\bo\subset \Fun(\Delta^1,\zT)$ is a conservative isofibration. In particular, $\bo$ is an $\infty$-category.
\end{lm}
\begin{proof}
	We have to show that $\bo$ is a \emph{replete subcategory} (in the sense of \cite[\href{https://kerodon.net/tag/01CF}{Definition 01CF}]{kerodon} and \cite[\href{https://kerodon.net/tag/01EX}{Example 01EX}]{kerodon}) of $\Fun(\Delta^1,\zT)$. First, we verify the conditions of \cite[\href{https://kerodon.net/tag/01CR}{Corollary 01CR}]{kerodon} to prove that $\bo$ is a subcategory of $\Fun(\Delta^1,\zT)$. As the condition of compatibility with extension of \ref{df:bo} only depends of the image of the morphisms in the $1$-category $\Fin_*$, one easily verifies that the set of morphisms in $\bo$ contain all identities of objects in $\bo$ is closed under homotopy and composition, as desired.\\
    Next, we turn to the proof that $\bo$ is replete. Let $f^+\colon \sigma^+\to \tau^+$ be an equivalence in $\Fun(\Delta^1,\zT)$ with $\sigma^+\in \bo$. We have to show that both $\tau^+$ and $f^+$ belong to $\bo$. Since the canonical functor $\Tw(\zO^\o)\to \zO^\o \times(\zO^\o)\op$ is conservative (being a right fibration), we deduce that in the diagram induced\todo{[PAS GRAVE] rendre ce diagramme explicite} by 
	$f^+$ in $\zO^\o$, all four morphisms $f_0\colon S_0\to T_0$, $f_0^+\colon S_0^+\to T_0^+$, $f_1\colon T_1\to S_1$, $f_1^+\colon T_1^+\to S_1^+$ are equivalences. From this and the commutativity of the square \eqref{eqn:diag_definition_bo}, one obtains that $\tau_0\colon T_0\to T_0^+$ is semi-inert, lies over an injection $\langle t\rangle \inc \langle t+1 \rangle$ and that $f_0^+$ is compatible with extension. Similarly, the morphisms $\sigma_1$, $f_1$ and $f_1^+$ are equivalences, therefore so must be $\tau_1$. This concludes the proof.
\end{proof}

\begin{lm}
	\label{lm:bo_sigma_is_Kan}
	Assume that the underlying $\infty$-category $\zO$ of $\zO^\o$ is an $\infty$-groupoid and let $\sigma\in \zT$. Then the fiber $\bo_\sigma$ of $\pi$ at $\sigma$ is a Kan complex. 
\end{lm}

\begin{proof}
	By the previous lemma, the inclusion $\bo \subset \Fun(\Delta^1,\zT)$ is a conservative isofibration. So is the map $\mathrm{ev}_0 \colon \Fun(\Delta^1,\zT)\to \zT$, hence $\pi$ is an isofibration. To prove the result, it now suffices to show that $\pi$ is conservative. Consider a morphism $f\colon \sigma^+\to \tau^+$ in $\bo$ whose image $\pi(f)$ in $\zT$ is an equivalence. The data of $f$ is that of a diagram of the following form:
\begin{equation}
\begin{tikzcd}[row sep = small]
                                             & S_0^+ \arrow[dd, "\sigma_0^+" near end] \arrow[rr, "f_0^+"] &                           & T_0^+ \arrow[dd,"\tau^+"]                    \\
S_0 \arrow[dd,"\sigma"] \arrow[rr,"f_0" near end] \arrow[ru, "\sigma_0",hook,harpoon] &                                             & T_0 \arrow[dd,"\tau", near start] \arrow[ru,"\tau_0",hook,harpoon] &                                     \\
                                             & S_1^+ \arrow[ld, "\sigma_1"]                    &                           & T_1^+. \arrow[ld,"\tau_1"] \arrow[ll, "f_1^+", near end] \\
S_1                                          &                                             & T_1 \arrow[ll,"f_1"]            &                                    
\end{tikzcd}
\end{equation}
Using that the fibration $\zT\to \zE\times\zE\op$ is conservative, we deduce that $f_0$ and $f_1$ are equivalences. 
By definition of the objects in $\bo$, the maps $\sigma_1$ and $\tau_1$ are also equivalences, therefore so is $f_1^+$. Finally, we claim that $f_0^+$ is an equivalence. To see this, write $f_0^+$ as the sum $f_0 \oplus {f_0^+}|_{+}$, where ${f_0^+}|_{+} \colon S_0^+\setminus S_0\to T_0^+\setminus T_0$ is the restriction of $f_0^+$ to the new color. Since ${f_0^+}|_{+}$ is a map in $\zO$, which by assumption is an $\infty$-groupoid, it is an equivalence; therefore so is $f$.
\end{proof}

\subsection{Proof of theorem \ref{thm:main_thm}} 
\label{sec:proof_of_the_main_theorem}


One of the key steps in the proof of theorem \ref{thm:main_thm} is the following result, whose proof is given in section \ref{sec:cartesianity_of_brane_fibration}.

\begin{thm}
	\label{thm:pi_is_cartesian_fibration}
	Let $\zO^\o$ be a unital $\infty$-operad.
	Then the functor $\pi\colon \bo\to \zT$ is a cartesian fibration.
\end{thm}

Assuming theorems \ref{thm:pi_is_cartesian_fibration} and \ref{thm:comparaison_bo_ext_intro}, whose proofs will be given in subsection \ref{sec:cartesianity_of_brane_fibration} and \ref{sec:comparison_with_lurie_s_spaces_of_extensions}, we can prove theorem \ref{thm:main_thm}.

\begin{proof}
	[Proof of theorem \ref{thm:main_thm}]
	By Mann--Robalo's argument (as described in section \ref{ssub:mann_robalo_s_strategy}), in order to prove the theorem it suffices to construct a right fibration over $\zT$, with fibers equivalent to spaces of extensions and whose associated functor $\zT\to \zS\op$ is a weak cartesian structure and satisfies the pullback condition of proposition \ref{prop:universal_property_Tw_Span}.

	Theorem \ref{thm:pi_is_cartesian_fibration} ensures that $\pi$ is a cartesian fibration. By lemma \ref{lm:bo_sigma_is_Kan}, its fibers are Kan complexes, hence $\pi$ is a right fibration. Moreover, theorem \ref{thm:comparaison_bo_ext_intro} identifies the fiber $\bo_\sigma$ over an object $\sigma\in \zT$ as its space of extensions $\Ext(\sigma)$.
    Therefore, it remains to show that the functor $F_\pi \colon \zT\to \zS\op$ classifying the right fibration $\pi$ is a weak cartesian structure and satisfies the pullback condition. The latter is exactly the condition that the $\infty$-operad $\zO^\o$ is coherent, using the equivalence $\bo_\sigma\simeq \Ext(\sigma)$\todo{[OUBLI POUR L'INSTANT] prouver la naturalité ?}. 
	For the weak cartesian condition, let $\sigma$ in $\zT$ be decomposed as a sum $\sigma \simeq \oplus_{i=1}^n \sigma_i$ of objects in $\zT_{\langle 1\rangle} \simeq \Tw(\oact)$. Since $p_1$ is constant along fibers of $\pi$ (in the sense of remark \ref{rk:remarks_on_objects_in_bo}), the fiber $\bo_\sigma$ decomposes as a disjoint union of the spaces $\bo_{\sigma_i}$, so that the natural map $F_\pi (\sigma)\stackrel{\sim}{\leftarrow} \coprod_{i=1}^n F_\pi(\sigma_i)$ in $\zS$ is an equivalence. This shows that $F_\pi$ is a lax cartesian structure. To verify that it is in fact a weak cartesian structure, let $f\colon \sigma\to \sigma'$ in $\zT$ be a $p_1$-cocartesian lift of the unique active morphism $\langle n\rangle \to \langle 1\rangle$ in $\Fin_*$. By remark \ref{rk:condition_morphism_in_T_p_1-cocartesian}, this implies that the two maps $\mathrm{source}(\sigma)\to \mathrm{source}(\sigma')$ and $\mathrm{target}(\sigma)\leftarrow \mathrm{target}(\sigma')$ are equivalences, which in turn ensures that $F_\pi(f)$ is an equivalence, as desired.
\end{proof}

\subsection{Generalized version of theorem \ref{thm:main_thm}} 
\label{sec:generalized_version_of_theorem_thm:main_thm}

The brane action given by theorem \ref{thm:main_thm} can be generalized to the setting where $\zO^\o$ is a unital $\infty$-operad, without assuming that its underlying $\infty$-category $\zO$ is an $\infty$-groupoid (condition \ref{item:df_coherence_b} in the definition of coherence given in \ref{df:coherence}).

To make this claim precise, let us say that an $\infty$-operad $\zO^\o$ is \emph{categorically coherent} if it is unital and satisfies the variant $(c')$ of condition \ref{item:df_coherence_c} in definition \ref{df:coherence} in which one requires diagram \eqref{eqn:diag_df_coherence} to be a \emph{categorical} pushout square of $\infty$-categories (instead of a homotopy pushout square). Note that if $\zO^\o$ is unital with $\zO^\o$ an $\infty$-groupoid, i.e. $\zO^\o$ satisfies conditions \ref{item:df_coherence_a} and \ref{item:df_coherence_b} from definition \ref{df:coherence}, then conditions \ref{item:df_coherence_c} and $(c')$ actually coincide, since for Kan complexes, homotopy pushout squares are automatically categorical pushout squares. As a consequence, coherent $\infty$-operads are categorically coherent. The generalized version of theorem \ref{thm:main_thm} writes as follows.

\begin{thmintrovariant}
	\label{thm:main_thm_generalized}
	Let $\zO^\o$ be a categorically coherent $\infty$-operad. Then the collection of $\infty$-categories {$\{\Ext(\id_c)\}_{c\in \zO}$} carries a canonical $\zO$-algebra structure in $\Cospan(\Cat_\infty)$, which recovers that of theorem \ref{thm:main_thm} when $\zO$ is an $\infty$-groupoid.
\end{thmintrovariant}

The proof of theorem \ref{thm:main_thm_generalized} is almost the same as the one given above for theorem \ref{thm:main_thm}, only slightly simpler. Indeed, most of the arguments, including the use of theorems \ref{thm:pi_is_cartesian_fibration} and \ref{thm:comparaison_bo_ext_intro}, do not use the assumption that $\zO$ is an $\infty$-groupoid. The only difference is that in the situation of theorem \ref{thm:main_thm_generalized}, $\pi$ is merely a cartesian fibration (as opposed to a right fibration) and therefore its classifying functor is of the form $\zT\to \Cat^\mathrm{op}_\infty$.
\bigskip

Following \cite{toen_branes}, one may go one step further in generality by dropping the assumption that $\zO^\o$ is coherent, that is assuming only that $\zO^\o$ is a unital $\infty$-operad. In this case the brane action merely gives a \emph{lax} algebra structure on the $\infty$-category $\Ext(\sigma)$ in cospans of $\infty$-categories, which is an genuine algebra structure precisely when $\zO^\o$ is coherent (in the previous generalized sense).
We refer to Kern's thesis \cite{kern_these} for more details on this lax structure.

	\section{Cartesianity of the brane fibration} 
\label{sec:cartesianity_of_brane_fibration}

This section is devoted to the proof of theorem \ref{thm:pi_is_cartesian_fibration}, asserting that the brane fibration $\pi \colon \bo\to \Tw(\Env(\zO))^\o$ of definition \ref{df:bo} is indeed a cartesian fibration. We will define particular lifts of edges along $\pi$ and then show that these are cartesian arrows in $\bo$ in the rest of the section.
Note that cartesianity of this fibration is the property ensures the existence of all the homotopical coherences involved in the definition of the $\zO$-algebra in $\Cospan(\zS)$ given by the brane action.



\subsection{Construction of cartesian lifts}
\label{sec:construction_f^+}

\begin{figure}[h]
\centering
	\label{fig:diags_steps}
\begin{tikzcd}[row sep = small, column sep = scriptsize]
                                       &  &                           & T_0^+ \arrow[dd]              \\
S_0 \arrow[dd] \arrow[rr] \arrow[rrru,"\tilde{f}_0",dashed] &  & T_0 \arrow[dd] \arrow[ru,hook,harpoon] &                               \\
                                       &  &                           & T_1^+ \arrow[ld] \arrow[llld,dashed] \\
S_1                                    &  & T_1 \arrow[ll]            &                              
\end{tikzcd}
\begin{tikzcd}[row sep = small, column sep = scriptsize]
                                               &                                       &  & T_0^+ \arrow[dd] \arrow[lldd, dashed] \\
S_0 \arrow[dd] \arrow[rrru] \arrow[rd, dashed] &                                       &  &                                       \\
                                               & S_1^+ \arrow[ld, Rightarrow, no head] &  & T_1^+ \arrow[llld] \arrow[ll, dashed] \\
S_1                                            &                                       &  &                                      
\end{tikzcd}
\begin{tikzcd}[row sep = small, column sep = scriptsize]
                                                          & S_0^+ \arrow[rr, "f_0^+", dashed]     &  & T_0^+ \arrow[dd] \arrow[lldd] \\
S_0 \arrow[dd] \arrow[rrru] \arrow[rd] \arrow[ru, dashed] &                                       &  &                               \\
                                                          & S_1^+ \arrow[ld, Rightarrow, no head] &  & T_1^+ \arrow[llld] \arrow[ll] \\
S_1                                                       &                                       &  &                              
\end{tikzcd}
\caption{Diagrams of steps 1, 2 and 3 of the construction of $f^+$. An arrow is dashed if it is added at the current step.}
\end{figure}
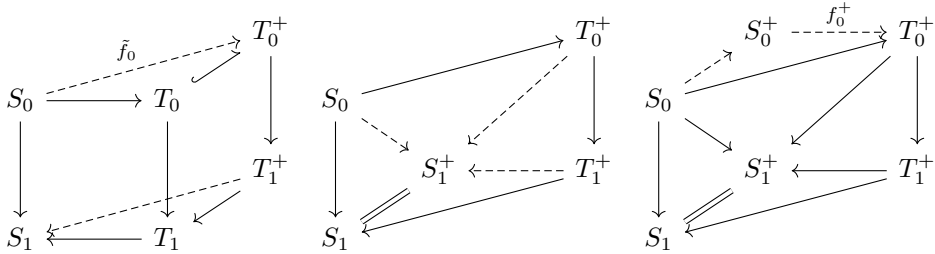

Let $f \col \sigma \leadsto \tau$ be a morphism in $\zT$ and let $e_\tau\col\tau\leadsto \tau^+$ be in the fiber $\bo_\tau$. We will construct a cartesian edge $f^+ \col \sigma^+\leadsto \tau^+$ lying $\pi$-above $f$.

\[
\begin{tikzcd}  
	\sigma^+ \ar[r,"f^+", dashed] 	&	\tau^+   & & \bo\ar[d,"\pi"] \\
	\sigma \ar[r,"f"] \ar[u,dashed]	& 	\tau\ar[u,"e_\tau"]&& \zT
\end{tikzcd}
\]
Unraveling the definition of $\bo$, we are given a diagram of the form
\begin{equation}
	\label{equ:construction_cartesian_lift_initial_data}
\begin{tikzcd}[row sep = small]
                          &  &                           & T_0^+ \arrow[dd,"\tau^+"] \\
S_0 \arrow[dd,"\sigma"] \arrow[rr,"f_0"] &  & T_0 \arrow[dd,"\tau"] \arrow[ru,"\tau_0",hook,harpoon] &                  \\
                          &  &                           & T_1^+ \arrow[ld,"\tau_1"] \\
S_1                       &  & T_1 \arrow[ll,"f_1"]            &                 
\end{tikzcd}
\end{equation}
in $\zE$ and want to extend it to one of the shape
\begin{equation}
\label{eqn:diagramme_cube}
\begin{tikzcd}[row sep = small]
                                             & S_0^+ \arrow[dd, dashed,""] \arrow[rr, dashed,"f_0^+"] &                           & T_0^+ \arrow[dd,"\tau^+"]                    \\
S_0 \arrow[dd,"\sigma"] \arrow[rr,"f_0" near end] \arrow[ru, dashed,"\sigma_0",hook,harpoon] &                                             & T_0 \arrow[dd,"\tau", near start] \arrow[ru,"\tau_0",hook,harpoon] &                                     \\
                                             & S_1^+ \arrow[ld, dashed,"\sigma_1"]                    &                           & T_1^+ \arrow[ld,"\tau_1"] \arrow[ll, dashed,"f_1^+", near end] \\
S_1                                          &                                             & T_1 \arrow[ll,"f_1"]            &                                    
\end{tikzcd}
\end{equation}
so that the resulting morphism $f^+\colon \sigma^+\leadsto \tau^+$ is in $\bo$. We proceed in several steps, depicted in figure \ref{fig:diags_steps}.

\begin{enumerate}[label=\bfseries{Step \arabic{enumi}.}]
    \item Pick a representative $\tilde{f}$ for the composite $e_\tau \circ f$ in $\zT$. In particular, this yields a $3$-simplex $S_0T_0^+T_1^+S_1$ and a $5$-simplex $S_0T_0T_0^+T_1^+T_1S_1$ extending  diagram  \eqref{equ:construction_cartesian_lift_initial_data}.

	\item Define the object $S_1^+$ as $S_1$ and the morphism $\sigma_1 \col S_1^+\to S_1$ as the identity. Since $\sigma_1$ is an equivalence, by using Joyal's lifting theorem \cite[\href{https://kerodon.net/tag/019F}{Theorem 019F}]{kerodon} and several horn fillers, we can extend the $3$-simplex $S_0T_0^+T_1^+S_1$ to a $4$-simplex $S_0T_0^+T_1^+S_1^+S_1$.

	\item 
	We now turn to the key step, namely the construction of the triangle $S_0S_0^+T_0^+$. Decompose $T_0^+$ as a sum of colors  
	$$T_0^+ = \oplus_{i\in p_0(T_0)}C_i \oplus C^+$$
	so that $C^+$ is the color lying above the element $p_0(T_0^+)\setminus \im(p_0(\tau_0))$. 
    Since $\zO^\o$ is unital, there exists an essentially unique morphism $\iota_{C^+}$ from the zero object of $\zO^\o$ to $C^+$. Define $S_0^+$ as the sum $S_0\oplus C^+$ and $\sigma_0$ as  $\id_{S_0}\oplus \iota_{C^+}$, which is clearly an atomic morphism.

	It remains to construct $f_0^+$. Note that $p_0(f_0^+)$ is required to coincide with the unique morphism $h\col p_0(S_0^+)\to p_0(T_0^+)$ that restricts to $p_0(\tilde{f}_0)$ on $p_0(S_0)$ and preserves $p_0(C^+)$. Consider the $\infty$-category	
	\[
	\zM = (\zO^\o)^{\Delta^2} \times_{(\zO^\o)^{\Delta^{\{2\}}}} \{T_0^+\} \times_{(\zO^\o)^{\Delta^{01}}} \{\sigma_0\} \times_{\Fin_*^{\Lambda_2^2}} \{ (p_0(\tilde{f}_0), h) \}
	\]
	consisting of all diagrams of the form
	\begin{equation}
		\label{eqn:diag_triangle_f0+}
	\begin{tikzcd}
	& S_0^+ \ar[rd,dashed,"a"] & \\
	S_0 \ar[ru,"\sigma_0",hook,harpoon] \ar[rr,dashed,"b"] & & T_0^+
	\end{tikzcd}
	\end{equation}
	satisfying $p_0(a) = h$ and $p_0(b) = p_0(\tilde{f}_0)$. 
	The inclusion $\Delta^{02}\inc \Delta^2$ yields a morphism
	\begin{equation}
		\label{eqn:map_composition_D}
	\zM \too \Map_{\zO^\o}^{p_0(\tilde{f}_0)} (S_0,T_0^+) \wk  \prod_{i = 1}^{s}\Map_{\zO^\o}^{\rho^i\circ h} (S_0,C_i).
	\end{equation}
	
	On the other hand, from the inner anodyne inclusion $\Lambda_1^2 \inc \Delta^2$ and the definition of $\infty$-operads, we get the following sequence of equivalences 
	\begin{eqnarray*}
	\zM & \wk & \Map_{\zO^\o}^h (S_0^+,T_0^+) \\
	& \wk & \prod_{i = 1}^{s}\Map_{\zO^\o}^{\rho^i\circ h} (S_0^+,C_i) \times 
	\Map_{\zO^\o}^{\rho^{n+1}\circ h} (S_0^+,C^+) \\
	& \wk & \prod_{i = 1}^{s}\Map_{\zO^\o}^{\rho^i\circ h} (S_0,C_i) \times \Map_\zO(C^+,C^+).
	\end{eqnarray*}
	Composing those equivalences with the projection $\Map_\zO(C^+,C^+)\to *$ recovers exactly the morphism (\ref{eqn:map_composition_D}). Therefore we see that the $\infty$-category of diagrams of the form (\ref{eqn:diag_triangle_f0+}) satisfying that $b = \tilde{f}_0$, which we identify with the fiber of the morphism (\ref{eqn:map_composition_D}) at $\tilde{f}_0$, is equivalent to $\Map_\zO(C^+,C^+)$.

	To define $f_0^+$ and a corresponding $2$-simplex of diagram  \eqref{eqn:diag_triangle_f0+}, it then suffices to specify any object in this $\infty$-groupoid $\Map_\zO(C^+,C^+)$. 

	\item At that point, we have extended the $3$-simplex $S_0T_0^+T_1^+S_1$ to a diagram of shape
	$$\Delta^{S_0T_0^+T_1^+S_1^+S_1} \cup_{\Delta^{S_0T_0^+}} \Delta^{S_0S_0^+T_0^+}.$$
	A simple computation shows that the inclusion of the latter simplicial set into $\Delta^{S_0S_0^+T_0^+T_1^+S_1^+S_1}$ is inner anodyne; this allows us to choose an extension of this diagram to a $5$-simplex $S_0S_0^+T_0^+T_1^+S_1^+S_1$. 

\end{enumerate}

This completes the construction of an edge $f^+ \col \sigma^+\leadsto \tau^+$ lifting $f$. The bulk of the proof of theorem \ref{thm:pi_is_cartesian_fibration} consists in proving that $f^+$ is cartesian.



\subsection{Outline of the proof of cartesianity}

Given a morphism $f\colon \sigma \to \tau $ in $\zT$ and an object $\tau^+$ in $\bo$, we have constructed a particular edge $f^+\colon \sigma^+\to \tau^+$ lying over $f$, which can be interpreted as an object in the $\infty$-category $\bo_{/\tau^+} \times_{\zT_{/\tau}} \zT_{/f}$.
The purpose of this section is to give an overview of the proof that $f^+$ is cartesian. The details will be dealt with in the rest of the section.


\begin{nota}
	Throughout the proof, we will make use of the notation introduced in \ref{not:lift}. In other words,
	from now on, we fix an object $\nu^+ \in \bo$, write $\zD_{\nu^+}$ for the $\infty$-category
	$\bo_{/\tau^+} \times_{\zT_{/\tau}} \zT_{/f} \times_{\bo} \{\nu^+\}$ and fix an object $u$ in it. We then consider the associated space of lifts
	\[
	\lift = \bo_{/f^+} \times_{\bo} \{\nu^+\} \times_{\zD_{\nu^+}}
	\{u\}.
	\]
\end{nota}

More explicitly, the datum of the object $u\in \zD_{\nu^+}	$ is that of a triangle $u_0$ in $\zT$ of the form
\begin{equation}
	\label{eqn:triangle_E}
 \begin{tikzcd}
 & \sigma \ar[dr,"f"] & \\
 \nu\ar[ru] \ar[rr,"g"] && \tau
 \end{tikzcd}
\end{equation}
together with a morphism $g^+ \col \nu^+\to \tau^+$ in $\bo$ lying $\pi$-above $g$.
An object in $\lift$ is a lift of $u$, that is the datum of a triangle
\begin{equation}
	\label{eqn:triangle_bo}
 \begin{tikzcd}
 & \sigma^+ \ar[dr,"f^+"] & \\
 \nu^+\ar[ru] \ar[rr,"g^+"] && \tau^+
 \end{tikzcd}
\end{equation}
in $\bo$ that lies $\pi$-above the triangle $u_0$ depicted in (\ref{eqn:triangle_E}). In the diagrams parametrized by this $\infty$-groupoid $\lift$, only the morphism $\nu^+\to \sigma^+$ and the $2$-simplex filling triangle \eqref{eqn:triangle_bo} are allowed to vary.\\

By lemma \ref{lm:equivalent_description_cartesian_edge}, proving that $f^+$ is $\pi$-cartesian amounts to showing the following result.

\begin{prop}
\label{prop:lift_is_contractible}
The space of lifts $\lift$ is contractible.	
\end{prop}




The rest of this section is devoted to the proof of proposition \ref{prop:lift_is_contractible}.
To help the reader, let us first explain the strategy of the argument: we will study the terminal morphism $q \colon \lift \to *$, decompose it as a composition
\begin{equation}
\label{eqn:q_as_composite}
	q 
	\col
	\zL
	\stackrel{q^{(0)}}{\too} \zL^{(0)}
	\stackrel{q^{(1)}}{\too} \zL^{(1)}
	\stackrel{q^{(2)}}{\too} \zL^{(2)}
	\stackrel{q^{(3)}}{\too} \zL^{(3)} \cong *
\end{equation}
and prove that each of the maps $q^{(i)}$ is an equivalence of Kan complexes. 
The idea is that each $\infty$-category $\lift^{(i)}$ parametrizes diagrams in $\zT$ of a certain shape $S^{(i)}$  and with certain data fixed. For $i>0$, the functors $q^{(i)}\colon \lift^{(i-1)}\to \lift^{(i)}$ can be interpreted as forgetful maps. 
The simplicial set $S^{(0)}$ is  $\Delta^2\times \Delta^1$ and corresponds to the shape of  diagrams in $\zE$ corresponding to triangle in $\bo$ (such as diagram \eqref{eqn:triangle_bo}).
The decreasing sequence of simplicial sets $S^{(0)}\supset S^{(1)}\supset S^{(2)}\supset S^{(3)}$ encode diagrams with fewer and fewer non-fixed data (see definition \ref{df:S0_S3}).

The following picture illustrates the decomposition of the composite functor $\lift^{(0)}\to \lift^{(3)}$.

\begin{equation}
\label{eqn:composition_q}
	\left\{
\begin{tikzcd}[cramped, column sep = 0.6em,row sep = small,execute at end picture={
\foreach \Nombre in  {V,Sp,Tp,Vp}
  {\coordinate (\Nombre) at (\Nombre.center);}
	\fill[black,opacity=0.2] 
  (V) -- (Sp) -- (Tp) --  cycle;
	\fill[black,opacity=0.1] 
  (V) -- (Vp) -- (Sp) --  cycle;
}]
                                                                                 &  & |[alias=Sp]|\sigma^+ \arrow[rrd]                      &  &                   \\
|[alias=Vp]|\nu^+ \ar[rru,dashed]\arrow[rrrr]                                                  &  &                                           &  & 
|[alias=Tp]|\tau^+\\
\\
\\
                                                                                 &  & |[alias=S]|\sigma \arrow[rrd, crossing over] \arrow[uuuu,crossing over]  &  &                   \\
|[alias=V]|\nu \ar[rruuuuu,dashed]\arrow[rru] \arrow[rrrr] \arrow[uuuu] &  &                                     &  & |[alias=T]|\tau   \ar[uuuu]             
\end{tikzcd}
\right\}
 \stackrel{q^{(1)}}{\too}
 \left\{
\begin{tikzcd}[cramped, column sep = 0.6em,row sep = small,execute at end picture={
\foreach \Nombre in  {V,S,Sp,Tp,Vp}
  {\coordinate (\Nombre) at (\Nombre.center);}
	\fill[black,opacity=0.2] 
  (V) -- (S) -- (Sp) --  cycle;
	\fill[black,opacity=0.08] 
  (V) -- (S) -- (Tp) --  cycle;
	\fill[black,opacity=0.12] 
  (S) -- (Sp) -- (Tp) --  cycle;
}]
                                                                                 &  & |[alias=Sp]|\sigma^+ \arrow[rrd]                      &  &                   \\
|[alias=Vp]|\nu^+ \arrow[rrrr]                                                  &  &                                           &  & 
|[alias=Tp]|\tau^+\\
\\
\\
                                                                                 &  & |[alias=S]|\sigma \arrow[rrd, crossing over] \arrow[uuuu,crossing over]  &  &                   \\
|[alias=V]|\nu \ar[rruuuuu,dashed]\arrow[rru] \arrow[rrrr] \arrow[uuuu] &  &                                     &  & |[alias=T]|\tau   \ar[uuuu]             
\end{tikzcd}
\right\}
\stackrel{q^{(2)}}{\too}
\left\{
\begin{tikzcd}[cramped, column sep = 0.6em,row sep = small,execute at end picture={
\foreach \Nombre in  {V,S,Tp,T}
  {\coordinate (\Nombre) at (\Nombre.center);}
	\fill[black,opacity=0.3] 
  (V) -- (S) -- (Tp) --  cycle;
	\fill[black,opacity=0.2] 
  (V) -- (S) -- (T) --  cycle;
	\fill[black,opacity=0.1] 
  (S) -- (T) -- (Tp) --  cycle;
}]
                                                                                 &  & \sigma^+ \arrow[rrd]                      &  &                   \\
\nu^+ \arrow[rrrr]                                                  &  &                                           &  & 
|[alias=Tp]|\tau^+\\
\\
\\
                                                                                 &  & |[alias=S]|\sigma \arrow[rrd, crossing over] \arrow[uuuu,crossing over]  &  &                   \\
|[alias=V]|\nu \arrow[rru] \arrow[rrrr] \arrow[uuuu] &  &                                     &  & |[alias=T]|\tau   \ar[uuuu]             
\end{tikzcd}
\right\}
 \stackrel{q^{(3)}}{\too}
 *
\end{equation}
In this description, solid arrows stand for morphisms in $\zT$ that are fixed within $\zL^{(i)}$, whereas dashed arrows indicate morphisms that are allowed to vary in that space. 
At each step of the composition, the new diagram is obtained from the previous one by removing one $3$-simplex in $\zT$ (and some simplices of smaller dimension), namely the $3$-simplices $\nu\nu^+\sigma^+\tau^+$, $\nu\sigma\sigma^+\tau^+$ and $\nu\sigma\tau\tau^+$, respectively for $q^{(1)}$, $q^{(2)}$ and $q^{(3)}$ (as indicated in grey in the picture).	
The last $\infty$-category $\lift^{(3)}\cong *$ should be thought of as the fixed data in the $\lift^{(i)}$.

The functors $q^{(2)}$ and $q^{(3)}$ are both induced by inner anodyne morphisms and will therefore be trivial Kan fibrations. The case of the functor $q^{(1)}$ is more delicate, and proving that it is also a trivial fibration will constitute the heart of the proof of proposition \ref{prop:lift_is_contractible}.

We divide the argument outlined above in three steps: each one amounts to proving that some of the functors $q^{(i)}$ are equivalences. We postpone the most technical parts of the proof to the end of the section (subsection \ref{sec:proof_of_technical_lemmas}).

\subsection{From slices to functor categories: the functor \texorpdfstring{$q^{(0)}$}{q0}} 
\label{sec:from_slices_to_functors_the_functor_}

The first step is to define the functor $q^{(0)} \colon \lift\to \lift^{(0)}$ and prove that it is a categorical equivalence. The $\infty$-category $\zL^{(0)}$ will be a slight variation of $\lift$, in that these two $\infty$-categories both parametrize triangles of the form \eqref{eqn:triangle_bo} with the following data fixed: the morphisms $f^+$ and  $g^+$ in $\bo$ and the triangle $u_0$ underlying $u$ of shape \eqref{eqn:triangle_E}. The two $\infty$-categories thus share the same objects, the difference being that $\lift$ is constructed from the slice $\infty$-category $\bo_{/f^+}$ whereas $\zL^{(0)}$ is obtained from the functor $\infty$-category $\Fun(\Delta^2,\bo)$. More precisely, we define $\zL^{(0)}$ as
\begin{equation}
	\label{eqn:df_lift^0}
\zL^{(0)} = \bo^{\Delta^2} \times_{\bo^{\Lambda_2^2}} \{(f^+,g^+)\}
\times_{\zT^{\Delta^2}} \{u_0\}.
\end{equation}

\begin{lm}
\label{lm:q0_equivalence}
There exists an equivalence of $\infty$-categories $q^{(0)} \col \lift\to \zL^{(0)}$. In particular, $\zL^{(0)}$ is a Kan complex.
\end{lm}

To construct this equivalence $q^{(0)}$, we first need a comparison between slice $\infty$-categories and corresponding $\infty$-categories of diagrams, given by the following lemma.

\begin{lm}
\label{lm:equivalence_slice_diagrams}
	Let $\zC$ be an $\infty$-category and $p \colon K \to \zC$ a diagram. Then there is a canonical equivalence of $\infty$-categories
	\[
	\zC_{/p} \stackrel{\sim}{\too} 
	\zC^{K^\triangleleft} \times_{\zC^{K}} \{p\}.
	\]
\end{lm}

For the sake of completeness, we provide a proof of this folklore result at the end of this section, see \ref{sub:proof_of_lemma_lm:equivalence_slice_diagrams}.
We can now proceed to the proof of lemma \ref{lm:q0_equivalence}.

\begin{proof}
	[Proof of lemma \ref{lm:q0_equivalence}]
	First, note that we can write
	\[
	\zL^{(0)} =
	\left(\bo^{\Delta^2} \times_{\bo^{\Delta^{12}}}  \{f^+\} \right)
	\times_\zP \{u\},
	\]
	where $\zP$ denotes the pullback
	\[
	\zP = 
	\left( \bo^{\Delta^{02}} \times_{\bo^{\Delta^{\{2\}}}} \{\tau^+\} \right)\times_{
	\left( \zT^{\Delta^{02}} \times_{\zT^{\Delta^{\{2\}}}} \{\tau\} \right)
	}
	\left( \zT^{\Delta^{2}}  \times_{\zT^{\Delta^{12}}} \{f\} \right).
	\]
	We define the functor $q^{(0)} \colon \lift\to \zL^{(0)}$ as the one induced from the commutative square
	\begin{equation}
		\label{eqn:comparison_slice_diagrams}
	\begin{tikzcd}  
		\bo_{/f^+} \ar[r,"\psi'"] \ar[d,"\xi'" left]
		&	\bo^{\Delta^2} \times_{\bo^{\Delta^{12}}}  \{f^+\} \ar[d,"\xi"] \\
		\bo_{/\tau^+}\times_{\zT_{/\tau}} \zT_{/f} \ar[r,"\psi"]					& 	\zP
	\end{tikzcd}
	\end{equation}
	by taking the fiber at $u \in \bo_{/\tau^+}\times_{\zT_{/\tau}} \zT_{/f}$. 

	To prove that $q^{(0)}$ is an equivalence of $\infty$-categories, we will use that Joyal's model structure is \emph{locally right proper} (although it is not right proper). This property holds for any model structure, and means that for any diagram $d = (X_a \to X_b \lo X_c)$ of fibrant objects in which one of the maps is a fibration, the canonical morphism 
	$X_a \times_{X_b} X_c \to 
	X_a \times^\mathrm{h}_{X_b} X_c$
	from the pullback to the homotopy pullback is an equivalence. In particular, if a morphism of such diagrams $d\to d'$ is a pointwise weak equivalence, then the induced morphism $\lim d \to \lim d'$ is a weak equivalence.

	Since each of the simplicial sets in diagram \eqref{eqn:comparison_slice_diagrams} is an $\infty$-category (hence a fibrant object in Joyal's model structure), to show that $q^{(0)}$ is an equivalence of $\infty$-categories, it suffices to establish that the following two claims~:
	\begin{enumerate}
		\item $\psi$ and $\psi'$ are categorical equivalences,
		\item $\xi$ and $\xi'$ are isofibrations.
	\end{enumerate}

	To prove the first claim, we make again use of the argument described in the previous paragraphs. Indeed, the morphism $\psi$ is itself induced from the natural transformation of diagrams
	\[
	\begin{tikzcd}
	\bo_{/\tau^+} \ar[r]	\ar[d,"\psi_a"] &
	{\zT_{/\tau}} 	\ar[d,"\psi_b"] &
	\zT_{/f} \ar[l]	\ar[d,"\psi_c"]
\\
	\bo^{\Delta^{02}} \times_{\bo^{\Delta^{\{2\}}}} \{\tau^+\} \ar[r] &
	\zT^{\Delta^{02}} \times_{\zT^{\Delta^{\{2\}}}} \{\tau\} &
	\zT^{\Delta^{2}}  \times_{\zT^{\Delta^{12}}} \{f\}. \ar[l,"\chi"]
	\end{tikzcd}
	\]
	Lemma \ref{lm:equivalence_slice_diagrams} guarantees that each of the vertical morphism are equivalences. We know that $\zT_{/f} \to {\zT_{/\tau}}$ is a right fibration (by the dual of Proposition 2.1.2.1 in \cite{htt}), hence an isofibration. 

	We now prove that the functor $\chi \colon \zT^{\Delta^{2}}  \times_{\zT^{\Delta^{12}}} \{f\} \to
	\zT^{\Delta^{02}} \times_{\zT^{\Delta^{\{2\}}}} \{\tau\}
	$
	is an isofibration. Let $v$ be an object in $ \zT^{\Delta^{2}}  \times_{\zT^{\Delta^{12}}} \{f\}$ and $\chi(v)\wk \bar{w}$ an equivalence in $\zT^{\Delta^{02}} \times_{\zT^{\Delta^{\{2\}}}} \{\tau\}$. We want to lift this equivalence to one in $\zT^{\Delta^{2}}  \times_{\zT^{\Delta^{12}}} \{f\}$. The datum of $v$ is that of a triangle in $\zT$ of the form
	\begin{equation}
	\label{eqn:diag_triangle_v}
	\begin{tikzcd}
 & \sigma \ar[dr,"f"] & \\
 \alpha \ar[ru] \ar[rr,"\ell"] && \tau.
	\end{tikzcd}
\end{equation}
The datum of the morphism $\chi(v)\to \bar{w}$ is that of a commutative square of the form
\begin{equation}
\label{diag:square_proof_isofib}
\begin{tikzcd}  
	\alpha \ar[r,"\ell"] \ar[d]\ar[dr]	&	\tau \ar[d,equal] \\
	\alpha' \ar[r]					& 	\tau'.
\end{tikzcd}
\end{equation}
As a natural transformation is an equivalence if it so pointwise, the fact that $\chi(v)\to \bar{w}$ is an equivalence translates into the statement that $\alpha \to \alpha'$ is an equivalence.
We want to extend the previous diagram into one of the form
\begin{equation}
	\label{diag:xi_isofib_1}
\begin{tikzcd}[row sep = small]
                                                & \sigma \arrow[rd, "f"]  &                             \\
\alpha \arrow[rr, "\ell" near start] \arrow[ru] \arrow[dd] &                                               & \tau \arrow[dd, equal] \\
                                                & \sigma' \arrow[rd, "f"]   \arrow[uu, crossing over, equal]                     &                             \\
\alpha' \arrow[rr] \arrow[ru]                   &                                               & \tau'                       
\end{tikzcd}
\end{equation}
which represents an equivalence in $\zT^{\Delta^{2}}  \times_{\zT^{\Delta^{12}}} \{f\}$. We do this construction in several steps~: first, by gluing diagrams \eqref{eqn:diag_triangle_v} and \eqref{diag:square_proof_isofib} and adding degenerate $2$-simplices $\sigma\tau\tau'$ and $\sigma\sigma'\tau'$, we obtain the diagram
\[
\begin{tikzcd}[row sep = small]
                                                & \sigma \arrow[rd, "f"]  &                             \\
\alpha \arrow[rr, "\ell" near start] \arrow[ru] \arrow[dd] &                                               & \tau \arrow[dd, equal] \\
                                                & \sigma' \arrow[rd, "f"] \arrow[uu, crossing over, equal]                        &                             \\
\alpha' \arrow[rr]                 &                                               & \tau'.
\end{tikzcd}
\]
From this point, the construction of diagram \eqref{diag:xi_isofib_1} is obtained using successive horn fillers in $\zT$, that is to say a sequence of choices of solutions to lifting problems, each of the form
\[
\begin{tikzcd}
\Lambda_k^n \ar[d]\ar[r] & \zT .\\
\Delta^n \ar[ru, dashed] & 
\end{tikzcd}
\]
First, choosing a filler of the horn of shape $\Lambda_2^3$ in $\alpha\sigma\tau\tau'$, we construct the $2$-simplex $\alpha\sigma\tau'$. Similarly, by filling the horn $\Lambda_1^2$ in $\alpha\sigma\sigma'$, we obtain a morphism $\alpha\sigma'$. Using a filler of the horn $\Lambda_1^2$ in $\alpha\sigma\sigma'\tau'$, we get a $2$-simplex $\alpha\sigma'\tau'$. Finally, using that the morphism $\alpha\alpha'$ is an equivalence, we can fill the horn $\Lambda_0^2$ in $\alpha\alpha'\sigma'$, as well as the horn $\Lambda_0^3$ in $\alpha\alpha'\sigma'\tau'$. This yields a diagram of the form \eqref{diag:xi_isofib_1} in which $\alpha\to\alpha'$ is an equivalence; hence an equivalence $v\wk w$ lifting the given morphism $\chi(v)\to \bar{w}$ along $\chi$. This concludes the proof of the first claim.\\

	We now come to the second claim. As $\xi'$ is a right fibration, it is in particular an isofibration. It remains to prove that $\xi$ is also an isofibration. 
	Consider an object $x\in \bo_{/\tau^+}\times_{\zT_{/\tau}} \zT_{/f}$ and an equivalence $\xi(x)\wk \bar{y}$ in $\zP$. We want to construct an equivalence $x\to y$ lifting $\xi(x)\wk \bar{y}$.
	The data of $x$ is that of a triangle in $\bo$ of the form
	\begin{equation}
	\label{diag:xi_isofib_2}
		\begin{tikzcd}
			& \sigma^+ \ar[rd, "f^+"] & \\
			\alpha^+ \ar[ru] \ar[rr,"\ell^+"] & & \tau^+.
		\end{tikzcd}
	\end{equation}
	The data of the morphism $\xi(x)\to \bar{y}$ is that of a diagram of the form \eqref{diag:xi_isofib_1} and a lift
\begin{equation}
\label{diag:xi_isofib_3}
\begin{tikzcd}  
	\alpha^+ \ar[r,"\ell^+"] \ar[d,"" left]\ar[dr]	&	\tau^+ \ar[d,equal] \\
	\alpha'^+ \ar[r,""]					& 	\tau'^+
\end{tikzcd}
\end{equation}
of its subdiagram \eqref{diag:square_proof_isofib} along $\pi$. By \cite[\href{https://kerodon.net/tag/01H4}{Corollary 01H4}]{kerodon}, since the morphism $\chi$ is an isofibration,
the maximal $\infty$-subgroupoid $\zP^{\simeq}$ of $\zP$ is given by the limit of the following diagram of $\infty$-categories
\[
 \begin{tikzcd}
 	\left(\bo^{\Delta^{02}} \times_{\bo^{\Delta^{\{2\}}}} \{\tau^+\}\right)^{\simeq} \ar[r] &
	\left(\zT^{\Delta^{02}} \times_{\zT^{\Delta^{\{2\}}}} \{\tau\}\right)^{\simeq} &
	\left(\zT^{\Delta^{2}}  \times_{\zT^{\Delta^{12}}} \{f\}\right)^{\simeq}. \ar[l,"\chi" above]
 \end{tikzcd}
\] 
The morphism $\xi(x)\to \bar{y}$ being an equivalence therefore translates into the fact that the morphism $\alpha^+\to \alpha'^+$ from is an equivalence. Our aim is to extend diagrams \eqref{diag:xi_isofib_2} and \eqref{diag:xi_isofib_3} to obtain a lift
\begin{equation}
	\label{diag:xi_isofib_goal}
\begin{tikzcd}[row sep = small]
                                                & \sigma^+ \arrow[rd, "f^+"]  &                             \\
\alpha^+ \arrow[rr,"\ell^+" near start] \arrow[ru] \arrow[dd] &                                               & \tau^+ \arrow[dd, equal] \\
                                                & \sigma'^+ \arrow[rd, "f^+"]   \arrow[uu, crossing over, equal]                     &                             \\
\alpha'^+ \arrow[rr] \arrow[ru]                   &                                               & \tau'^+              
\end{tikzcd}
\end{equation}
of the given diagram \eqref{diag:xi_isofib_1} along the functor $\pi$. The construction of diagram \eqref{diag:xi_isofib_goal} is given by solutions to the same sequence of horn filling problems as that of the proof that $\chi$ is an isofibration (in claim (1)); the only difference being that in the present case, the horn filling problems have to be considered relative to the inner fibration $\pi$, that is as problems of lifting of the form
\[
\begin{tikzcd}
\Lambda_k^n \ar[d]\ar[r] & \bo \ar[d,"\pi"]\\
\Delta^n \ar[r]\ar[ru, dashed] & \zT.
\end{tikzcd}
\]
This construction provides an equivalence $x \wk y$ as desired, ensuring that $\xi$ is an isofibration. This concludes the proof of lemma \ref{lm:q0_equivalence}.
\end{proof}

\subsection{Anodyne extensions: the functors \texorpdfstring{$q^{(2)}$ and $q^{(3)}$}{q2 and q3}} 
\label{sec:anodyne_extensions_the_functors_}

We have defined the functor $q^{(0)}\colon \lift\to \lift^{(0)}$ and proved it is an equivalence of $\infty$-groupoids.
We now spell out the rest of the decomposition \eqref{eqn:q_as_composite} of the terminal morphism $q\colon \lift \to *$, by defining the $\infty$-categories $\lift^{(i)}$, for $i\in\{1,2,3\}$, as well as the functors $q^{(i)}$.
\begin{df}
    \label{df:S0_S3}
    \begin{itemize}
        \item 
First, let $S^{(0)}$ denote the simplicial set $\Delta^{2}\times \Delta^1$ with vertices labelled with $\nu,\sigma,\tau,\nu^+,\sigma^+,\nu^+$ as in the diagrams represented in \eqref{eqn:composition_q}. 
        \item 
    For $i\in \{1,2,3\}$, we define decreasing subsimplicial sets $S^{(i)}$ of $S^{(0)}$
using the following formulas:
\begin{eqnarray}
S^{(1)} & = & 
\Delta^{\nu\nu^+\tau^+}  \cup_{\Delta^{\nu\tau^+}} \Delta^{\nu\sigma\tau\tau^+}
\cup_{\Delta^{\nu\sigma\tau^+}} \Delta^{\nu\sigma\sigma^+\tau^+},\nonumber
\\
\label{eqn:df_S^i}
S^{(2)} & = &  
\Delta^{\nu\nu^+\tau^+}  \cup_{\Delta^{\nu\tau^+}} \Delta^{\nu\sigma\tau\tau^+}
\cup_{\Delta^{\sigma\tau^+}} \Delta^{\sigma\sigma^+\tau^+},
 \\
S^{(3)} & = & 
\Delta^{\nu\nu^+\tau^+}  \cup_{\Delta^{\nu\tau^+}} \Lambda_\tau^{\nu\sigma\tau\tau^+}
\cup_{\Delta^{\sigma\tau^+}} \Delta^{\sigma\sigma^+\tau^+}.\nonumber
\end{eqnarray}
    \end{itemize}
\end{df}

{\sloppy
For $i\in\{0,1,2\}$, the simplicial set $S^{(i)}$ will encode the shape of the diagrams parametrized by $\lift^{(i)}$, whereas $S^{(3)}$ will describe the shape of diagrams that are fixed within $\lift^{(i)}$. The simplicial sets $S^{(0)},\dots,S^{(3)}$ can be pictured as }
\begin{equation}
\label{eqn:description_S^i}
\begin{tikzcd}[cramped, column sep = 0.6em,row sep = small,execute at end picture={
\foreach \Nombre in  {V,Sp,Tp,Vp}
  {\coordinate (\Nombre) at (\Nombre.center);}
	\fill[black,opacity=0.3] 
  (V) -- (Sp) -- (Tp) --  cycle;
	\fill[black,opacity=0.15] 
  (V) -- (Vp) -- (Sp) --  cycle;
}]
                                                                                 &  & |[alias=Sp]|\sigma^+ \arrow[rrd]                      &  &                   \\
|[alias=Vp]|\nu^+ \ar[rru]\arrow[rrrr]                                                  &  &                                           &  & 
|[alias=Tp]|\tau^+\\
\\
\\
\\
                                                                                 &  & |[alias=S]|\sigma \arrow[rrd, crossing over] \arrow[uuuuu,crossing over]  &  &                   \\
|[alias=V]|\nu \ar[rruuuuuu]\arrow[rru] \arrow[rrrr] \arrow[uuuuu] &  &                                     &  & |[alias=T]|\tau   \ar[uuuuu]             
\end{tikzcd}, 
\quad
\begin{tikzcd}[cramped, column sep = 0.6em,row sep = small,execute at end picture={
\foreach \Nombre in  {V,S,Sp,Tp,Vp}
  {\coordinate (\Nombre) at (\Nombre.center);}
	\fill[black,opacity=0.3] 
  (V) -- (S) -- (Sp) --  cycle;
	\fill[black,opacity=0.1] 
  (V) -- (S) -- (Tp) --  cycle;
	\fill[black,opacity=0.2] 
  (S) -- (Sp) -- (Tp) --  cycle;
}]
                                                                                 &  & |[alias=Sp]|\sigma^+ \arrow[rrd]                      &  &                   \\
|[alias=Vp]|\nu^+ \arrow[rrrr]                                                  &  &                                           &  & 
|[alias=Tp]|\tau^+\\
\\
\\
\\
                                                                                 &  & |[alias=S]|\sigma \arrow[rrd, crossing over] \arrow[uuuuu,crossing over]  &  &                   \\
|[alias=V]|\nu \ar[rruuuuuu]\arrow[rru] \arrow[rrrr] \arrow[uuuuu] &  &                                     &  & |[alias=T]|\tau   \ar[uuuuu]             
\end{tikzcd}, 
\quad
\begin{tikzcd}[cramped, column sep = 0.6em,row sep = small,execute at end picture={
\foreach \Nombre in  {V,S,Tp,T}
  {\coordinate (\Nombre) at (\Nombre.center);}
	\fill[black,opacity=0.4] 
  (V) -- (S) -- (Tp) --  cycle;
	\fill[black,opacity=0.25] 
  (V) -- (S) -- (T) --  cycle;
	\fill[black,opacity=0.15] 
  (S) -- (T) -- (Tp) --  cycle;
}]
                                                                                 &  & \sigma^+ \arrow[rrd]                      &  &                   \\
\nu^+ \arrow[rrrr]                                                  &  &                                           &  & 
|[alias=Tp]|\tau^+\\
\\
\\
\\
                                                                                 &  & |[alias=S]|\sigma \arrow[rrd, crossing over] \arrow[uuuuu,crossing over]  &  &                   \\
|[alias=V]|\nu \arrow[rru] \arrow[rrrr] \arrow[uuuuu] &  &                                     &  & |[alias=T]|\tau   \ar[uuuuu]             
\end{tikzcd}, \quad
\begin{tikzcd}[cramped, column sep = 0.6em,row sep = small,execute at end picture={
\foreach \Nombre in  {V,S,Tp,T}
  {\coordinate (\Nombre) at (\Nombre.center);}
	\fill[black,opacity=0.25] 
  (V) -- (S) -- (T) --  cycle;
	\fill[black,opacity=0.15] 
  (S) -- (T) -- (Tp) --  cycle;
}]
                                                                                 &  & \sigma^+ \arrow[rrd]                      &  &                   \\
\nu^+ \arrow[rrrr]                                                  &  &                                           &  & 
|[alias=Tp]|\tau^+\\
\\
\\
\\
                                                                                 &  & |[alias=S]|\sigma \ar[rruuuu]\arrow[rrd, crossing over] \arrow[uuuuu,crossing over]  &  &                   \\
|[alias=V]|\nu \arrow[rru] \arrow[rrrr] \arrow[uuuuu] &  &                                     &  & |[alias=T]|\tau   \ar[uuuuu]
\arrow[from=V, to=Tp] 
\end{tikzcd}.
\end{equation}

We define $\lift^{(3)}$ as the terminal $\infty$-groupoid~; we think of its unique object as the diagram $S^{(3)}\to \zT$ given by the data $(f^+,g^+,u_0)$ that we fixed earlier on. The inclusion $j^{(i)}\colon S^{(i)}\subset S^{(i-1)}$ induces a forgetful functor $p^{(i)}\colon \zT^{S^{(i-1)}} \to \zT^{S^{(i)}}$ that we use to define the $\infty$-categories
\begin{eqnarray*}
\lift^{(1)} & = & 
\zT^{S^{(1)}} \times_{\zT^{S^{(3)}}} \{(f^+,g^+,u_0)\},\\
\lift^{(2)} & = & 
\zT^{S^{(2)}} \times_{\zT^{S^{(3)}}} \{(f^+,g^+,u_0)\}.
\end{eqnarray*}
Recall that $\lift^{(0)}$ was defined by formula \eqref{eqn:df_lift^0} using terms of diagrams with values in $\bo$.
Nevertheless, the following lemma ensures that $\lift^{(0)}$ actually admits a simple equivalent description in terms of diagrams in $\zT$, following the above pattern for $\lift^{(1)}$ and $\lift^{(2)}$.
\begin{lm}
\label{lm:rewrite_L0}
We have a canonical equivalence of $\infty$-groupoids
\[
\lift^{(0)} \simeq 
\zT^{S^{(0)}} \times_{\zT^{S^{(3)}}} \{(f^+,g^+,u_0)\}.
\]
\end{lm}
\begin{proof}
	It follows from its definition that the $\infty$-groupoid $\lift^{(0)}$ fits into the commutative diagram
	\begin{equation*}
		\begin{tikzcd}
			\lift^{(0)} \ar[r]\ar[d] \arrow[dr, phantom, "{\lrcorner}" , very near start, color=black] & *\ar[d] & \\
			\bo^{\Delta^2} \ar[r]\ar[d] & \bo^{\Lambda_2^2}\times_{\zT} \zT^{\Delta^2} \ar[r]\ar[d] \arrow[dr, phantom, "{\lrcorner}" , very near start, color=black]
			& \bo^{\Lambda_2^2} \ar[d] \\
			\zT^{\Delta^2 \times \Delta^1} \ar[r] & \zT^{S^{(3)}} \ar[r] & \zT^{\Lambda_2^2\times \Delta^1}
		\end{tikzcd}	
	\end{equation*}
	in which the upper left and the bottom right squares are cartesian. To prove the lemma, it therefore suffices to show that the bottom outer square is cartesian. Note that both vertical maps $\bo^{\Lambda_2^2}\to  \zT^{\Lambda_2^2\times \Delta^1}$ and $\bo^{\Delta^2}\to \zT^{\Delta^2\times \Delta^1}$ are subcategories, meaning they are monomorphisms that are inner fibrations. This implies that the morphism $\bo^{\Lambda_2^2} \times_{\zT^{\Lambda_2^2\times \Delta^1}} \zT^{\Delta^2}\to \zT^{\Delta^2\times \Delta^1}$ is also the inclusion of a subcategory, hence it is enough to verify that the two subcategories  $\bo^{\Delta^2}$ and $\bo^{\Lambda_2^2} \times_{\zT^{\Lambda_2^2\times \Delta^1}} \zT^{\Delta^2}$ have the same objects and morphisms.
	An object in $\bo^{\Delta^2}$ (respectively in $\bo^{\Lambda_2^2} \times_{\zT^{\Lambda_2^2\times \Delta^1}} \zT^{\Delta^2}$) is a diagram
\begin{equation*}
\begin{tikzcd}[cramped, row sep = tiny, row sep = small]
& \alpha_1^+\ar[dr,"t"] &  \\
\alpha_0^+\ar[ru,"s"]\ar[rr,"w" near start] & & \alpha_2^+
\\
\\
& \alpha_1\ar[uuu,crossing over]\ar[dr] &  \\
\alpha_0\ar[uuu]\ar[ru]\ar[rr] & & \alpha_2\ar[uuu]
\end{tikzcd}
\end{equation*}
in $\zT$ in which the maps $\alpha_0\alpha_0^+$, $\alpha_1\alpha_1^+$ and $\alpha_2\alpha_2^+$ are objects in $\bo$ and such that the morphisms $s$, $t$ and $w$ (respectively only $t$ and $w$) lie $p_1$-above $\id_{\langle 1\rangle}$ and are compatible with extension (conditions of definition \ref{df:bo}). The key observation is that whenever $t$ and $w$ both satisfy these properties, then so does $s$; this fact implies that the two subcategories have the same objects. One can use a similar argument to show that the same is true for morphisms of these two subcategories, as desired.

\end{proof}

Therefore, the functors $q^{(i)}$ and the $\infty$-categories $\lift^{(i)}$ fit into the commutative diagram
\begin{equation}
	\label{eqn:diag_with_all_q^i_as_pullback}
\begin{tikzcd}[column sep = large]
\lift^{(0)}\ar[r,"q^{(1)}"]\ar[d]\arrow[dr, phantom, "{\lrcorner}" , very near start] &
\lift^{(1)}\ar[r,"q^{(2)}"]\ar[d]\arrow[dr, phantom, "{\lrcorner}" , very near start] &
\lift^{(2)}\ar[r,"q^{(3)}"]\ar[d]\arrow[dr, phantom, "{\lrcorner}" , very near start] &
\lift^{(3)}\cong * \ar[d,"{(f^+,g^+,u_0)}"] \\
\zT^{S^{(0)}} \ar[r,"p^{(1)}"] &
\zT^{S^{(1)}} \ar[r,"p^{(2)}"] &
\zT^{S^{(2)}} \ar[r,"p^{(3)}"] &
\zT^{S^{(3)}},
\end{tikzcd}
\end{equation}
where all squares are cartesian.\\

We now claim that the inclusions $j^{(2)}\colon S^{(1)}\subset S^{(2)}$ and $j^{(3)}\colon S^{(2)}\subset S^{(3)}$ are inner anodyne. For the latter morphism, this is obvious from the formulas \eqref{eqn:df_S^i}, as $j^{(3)}$ is obtained as a pushout of the inner anodyne map $\Lambda_\tau^{\nu\sigma\tau\tau^+} \subset \Delta^{\nu\sigma\tau\tau^+}$. For the former map, note that we can write $j^{(2)}$ as a pushout of the composition
\begin{equation}
\Delta^{\nu\sigma\tau^+}
\cup_{\Delta^{\sigma\tau^+}} \Delta^{\sigma\sigma^+\tau^+} \ \subset \ 
\Lambda_\sigma^{\nu\sigma\sigma^+\tau^+} \ \subset \
\Delta^{\nu\sigma\sigma^+\tau^+},\nonumber
\end{equation}
in which both maps are inner anodyne. This shows the claim. We therefore obtain that the induced functors 
$p^{(2)}$ and
$p^{(3)}$ are trivial Kan fibrations. 
As every square in diagram \eqref{eqn:diag_with_all_q^i_as_pullback} is cartesian, we deduce that the functors
$q^{(2)}$ and
$q^{(3)}$ are also trivial Kan fibrations. In particular, we see that $\lift^{(1)}$ and $\lift^{(2)}$ are contractible Kan complexes.

\subsection{Existence and uniqueness of factorizations: the functor \texorpdfstring{$q^{(1)}$}{q1}} 
\label{sec:constructing_lifts_the_functor_}

In order to prove theorem \ref{prop:lift_is_contractible}, it only remains to prove that $q^{(1)}$ is a trivial Kan fibration, which is the main step of the proof. Since the inclusion $j^{(1)}$ restricts to a bijection
$S^{(1)}_0\subset S^{(2)}_0$ on the sets of $0$-simplices, by \cite[Proposition 40.6 and footnote 30]{rezk_notes},
we obtain that $p^{(1)}\colon \lift^{(1)}\to \lift^{(2)}$ is an isofibration. Since we already know that both $\lift^{(1)}$ and $\lift^{(2)}$ are Kan complexes, it follows that $p^{(1)}$ is a Kan fibration. Since the squares in diagram \eqref{eqn:diag_with_all_q^i_as_pullback} are cartesian, we deduce that $q^{(1)}$ is also a Kan fibration between Kan complexes. To see that it is a weak equivalence, we have to show that all of its fibers are contractible, which is asserted in the following proposition.


\begin{prop}
\label{prop:main_step_theorem}
Let $d\colon S^{(1)}\to \zT$ be a diagram in $\lift^{(1)}$. Then the fiber 
$\lift^{(0)}_{d}$
 of $q^{(1)}$ at $d$ is contractible.
\end{prop}

\begin{proof}
First, we claim that the question can be restricted to the full subdiagram of $d$ on the objects $\nu,\nu^+,\sigma^+$ and $\tau^+$. Indeed, if we let $S^{(1)'}$ denote the subsimplicial set $\Delta^{\nu\nu^+\tau^+}\cup_{\Delta^{\nu\tau^+}} \Delta^{\nu\sigma^+\tau^+}$ of $S^{(0)'} := \Delta^{\nu\nu^+\sigma^+\tau^+}$ encoding commutative squares of the form 
\begin{equation}
	\label{eqn:diag_step_I}
	\begin{tikzcd}  
		\nu^+ \ar[r] & \tau^+\\
		\nu\ar[u]\ar[r]\ar[ru] & \sigma^+\ar[u],
	\end{tikzcd}
\end{equation}
we observe that $j^{(1)}\colon S^{(1)}\subset S^{(0)}$ can be written as the pushout 
\begin{equation}
	\begin{tikzcd}  
		S^{(1)'} \ar[r,""] \ar[d,"" left]\arrow[dr, phantom, "{\ulcorner}" , very near end]
											&	S^{(0)'} \ar[d,""] \\
		S^{(1)} \ar[r,"j^{(1)}" below]			& 	S^{(0)},
	\end{tikzcd}
\end{equation}
of the inclusion. Therefore we can identify the fiber of $q^{(1)}$ at $d$ as
\[
\lift_d^{(0)} \simeq \zT^{S^{(0)'}} \times_{\zT^{S^{(1)'}}} \{d'\},
\]
where $d' = d|_{S^{(1)'}}$, thus proving the claim.

Note that the diagram $d'$, which is of shape \eqref{eqn:diag_step_I}, is essentially determined by the fixed data $(f^+,g^+,u_0)$. In particular, in this diagram, the morphism $\sigma^+\tau^+$ is the edge $f^+$ constructed in section \ref{sec:construction_f^+} and the $2$-simplices $\nu\sigma^+\tau^+$ and $\nu\nu^+\tau^+$ are given.
The proof of the proposition thus consists in showing that, from this data, the remaining simplices $\nu^+\sigma^+$, $\nu\nu^+\sigma^+$, $\nu^+\sigma^+\tau^+$ and $\nu\nu^+\sigma^+\tau^+$ can be constructed in an essentially unique way. 

Until now, the simplices were written in the $\infty$-category $\zT$; we need to reformulate the problem in terms of diagrams with values in $\zO^\o$. The diagram $d'$ of shape \eqref{eqn:diag_step_I} corresponds to a certain diagram $d'_\zO \colon K\to \zO^\o$, which can be pictured as
\begin{equation}
\label{eqn:diag_7_simplex}
\begin{tikzcd}[cramped, row sep = small, column sep = 1.8em]
                                                     & V_0^+ \arrow[ddd] \arrow[rrr] &  &                                       & T_0^+ \arrow[ddd]                          \\
V_0 \arrow[ru,hook,harpoon] \arrow[ddd] \arrow[rrr,crossing over] \arrow[rrrru,crossing over] &                                                   &  & S_0^+ \arrow[ru]           &                                            \\
                                                     &                                                   &  &                                       &                                            \\
                                                     & V_1^+ \arrow[ld]                                  &  &                                       & T_1^+ \arrow[ld] \arrow[lll] \arrow[lllld] \\
V_1                                                  &                                                   &  & S_1^+ \arrow[lll] &                                           
\ar[from=2-4,to=l,crossing over]
\end{tikzcd}
\quad 
\stackrel{\mathrm{def}}{=}
\quad
\begin{tikzcd}[row sep = 1.2em, column sep = 1.8em]
                                                     & 1 \arrow[ddd] \arrow[rrr] &  &                                       & 3 \arrow[ddd]                          \\
0 \arrow[ru,hook,harpoon] \arrow[ddd] \arrow[rrr,crossing over] \arrow[rrrru,crossing over] &                                                   &  & 2 \arrow[ru]           &                                            \\
                                                     &                                                   &  &                                       &                                            \\
                                                     & 6 \arrow[ld]                                  &  &                                       & 4. \arrow[ld] \arrow[lll] \arrow[lllld] \\
7                                                  &                                                   &  & 5 \arrow[lll] &                                           
\ar[from=2-4,to=l,crossing over]
\end{tikzcd}
\end{equation}\todo{[FAIT POUR CHAP 3] mettre les  ,hook,harpoon}
Here, in order to sometimes simplify notations, we denote the objects 
$V_0$, $V_0^+$, $S_0^+$, $T_0^+$, $T_1^+$, $S_1^+$, $V_1^+$, $V_1$ as the integers $0,\dots,7$, in the same order. This way, we can write the simplicial set indexing $d'_\zO$ as the subsimplicial set
$$K = \Delta^{013467}\cup_{\Delta^{0347}}\Delta^{023457}$$
of $\Delta^{7}$. The fiber $\lift^{(0)}_d$ is therefore canonically equivalent to the space $\lift_d^{(0)} \simeq (\zO^\o)^{\Delta^7} \times_{(\zO^\o)^{K}} \{d'_\zO\}$ parametrizing extensions of the diagram $d'_\zO$ to a $7$-simplex.

The key step of the proof concerns the space of lifts of the upper part of diagram (\ref{eqn:diag_7_simplex}), namely the subdiagram indexed by the full subsimplicial set $K_0 \subseteq K$ on the objects $V_0$, $V_0^+$, $S_0^+$ and $T_0^+$. Note that this simplicial set is isomorphic to $\Delta^1\times \Delta^1$.
Define $\zZ$ as the fiber $(\zO^\o)^{\Delta^3}\times_{(\zO^\o)^{K_0}} \{d'_\zO|_{K_0}\}$; this $\infty$-category parametrizes extensions of the diagram $d'_\zO|_{K_0} \colon K_0\to \zO$ to a $3$-simplex $V_0V_0^+S_0^+T_0^+$. Since the inclusion $K_0 \subset \Delta^3$ is a bijection on objects, the induced functor $(\zO^\o)^{\Delta^3}\to (\zO^\o)^{K_0}$ is an isofibration and therefore the $\infty$-category $\zZ$ is a space. We will show the following intermediate result.
\begin{claim}
	\label{claim:Z_is_contractible}
	The space $\zZ = (\zO^\o)^{\Delta^3}\times_{(\zO^\o)^{K_0}}\{d'_\zO |_{K_0}\}$ of lifts of the upper part of diagram \eqref{eqn:diag_7_simplex} is contractible.	
\end{claim}

The argument relies on the observation that the diagram $K_0\to \zO^\o \to \Fin_*$ has
a decomposition $p\circ d'_\zO|_{K_0} = d_{\Fin_*}^- \oplus d_{\Fin_*}^+$ given by
\begin{equation}
 	\label{eqn:diag_square_fin_*}
 	\begin{tikzcd}  
 		p_0(V_0)^+ \ar[r]					&	p_0(T_0^+) \\
 		p_0(V_0) \ar[r,""]\ar[u,hook,harpoon]		& 	p_0(S_0^+)\ar[u]
 	\end{tikzcd} 
 	=
 	\begin{tikzcd}  
 		\langle n \rangle \ar[r]					&	\langle k \rangle \\
 		\langle n \rangle \ar[r,""]\ar[u,equal]		& 	\langle m \rangle\ar[u]
 	\end{tikzcd} 
 	\oplus
 	\begin{tikzcd}  
 		\langle 1 \rangle \ar[r,equal]				&	\langle 1 \rangle \\
 		\langle 0 \rangle \ar[r,""]\ar[u,hook,harpoon]		& 	\langle 1 \rangle\ar[u,equal]
 	\end{tikzcd},
\end{equation}
where $\oplus$ stands for the operation of pointwise disjoint union of diagrams in $\Fin_*$. Using the identity maps in \eqref{eqn:diag_square_fin_*}, one readily sees that both diagrams $ d_{\Fin_*}^-$ and $d_{\Fin_*}^+$ extend uniquely to $3$-simplices $\widetilde{d}_{\Fin_*}^-$  and $\widetilde{d}_{\Fin_*}^+$ in $\Fin_*$, which implies that the diagram $p\circ d'_\zO|_{K_0}$ also extends uniquely to a $3$-simplex, namely $\widetilde{d'}_{\Fin_*} := \widetilde{d}_{\Fin_*}^-\oplus\widetilde{d}_{\Fin_*}^+$. In particular, any diagram $\Delta^3\to \zO^\o$ in $\zZ$ will be a lift of $\widetilde{d'}_{\Fin_*}$. This shows that we can rewrite the space $\zZ$ as
\[
	\zZ \simeq
	(\zO^\o)^{\Delta^3} \times_{(\zO^\o)^{K_0}} \left\{d'_\zO|_{K_0}\right\} \times_{\Fin_*^{\Delta^3}} \{\widetilde{d'}_{\Fin_*}\}.
\]
We will make use of decomposition \eqref{eqn:diag_square_fin_*} to obtain a splitting of the space $\zZ$, using the next lemma. Recall that a simplicial set is said to be \emph{braced} if every face of a nondegenerate simplex remains nondegenerate \cite[\href{https://kerodon.net/tag/00XU}{Tag 00XU}]{kerodon}.

\begin{lm}
\label{lm:operad_decomposition}
	Let $\zO^\o$ be any $\infty$-operad. 
	Let $J$ be a braced simplicial set and $I\subseteq J$ be a subsimplicial set. Consider two diagrams $q$ and $r$ making the following square commute:
	\begin{equation}
		\label{eqn:lm_key_statement}
	\begin{tikzcd}  
		I \ar[r,"r"] \ar[d,hook, left]	&	\zO^\o \ar[d,"p"] \\
		J \ar[r,"q"]					& 	\Fin_*
	\end{tikzcd} 
	\end{equation}
	and assume that $q$ decomposes as a disjoint union $q = \oplus_{i=1}^n q_i$ of diagrams $J\to \Fin_*$.
	Then there exists a decomposition $r \simeq \oplus_{i=1}^n r_i$ such that the $\infty$-category of lifts in the square \eqref{eqn:lm_key_statement} splits as the following direct product:
	\begin{equation}
		\label{eqn:lm_key_equivalence}
	(\zO^\o)^J \times_{(\zO^\o)^I} \{r\} \times_{\Fin_*^J} \{q\} \simeq 
	\prod_{i=1}^n
	(\zO^\o)^J \times_{(\zO^\o)^I} \{r_i\} \times_{\Fin_*^J} \{q_i\}.
	\end{equation}
\end{lm}
The proof of this lemma  is given at the end of this section, in section \ref{sec:proof_lemma_operad_decomposition}.

	We can now complete the proof of lemma \ref{lm:operad_decomposition}.
	Through the first of the equivalences in \eqref{eqn:decompositions_E^I_fiber}, the object $r\in (\zO^\o)^I$ is identified with an object $\oplus_i r_i$.
	Now observe that $\zD$ fits in a diagram of pullback squares
	\[
	\begin{tikzcd}
	\zD \ar[r]\ar[d] 
	\arrow[dr, phantom, "{\lrcorner}" , very near start]& 
	(\zO^\o)^J \times_{\Fin_*^J} \{q\} \ar[r]\ar[d,"\rho'"]
	\arrow[dr, phantom, "{\lrcorner}" , very near start]& 
	(\zO^\o)^J \ar[d,"\rho"] \\
	* \ar[r,"r"]\ar[rd, equal, bend right = 20] &
	(\zO^\o)^I\times_{\Fin_*^I} \{q|_I\}
	\arrow[dr, phantom, "{\lrcorner}" , very near start]
	 \ar[r] \ar[d] &
	(\zO^\o)^I \times_{\Fin_*^I} \Fin_*^J \ar[d]\\
	& * \ar[r,"q"] & \Fin_*^J.
	\end{tikzcd}
	\]
	Since the vertical functor $\rho$ is an isofibration, so is its pullback $\rho'$. Therefore the top left pullback square is invariant under equivalence of $\infty$-categories. Using the above two decompositions, we deduce that $\zD$ itself can be written as a product	
	$$\zD\,\simeq \,
	\prod_{i=1}^n (\zO^\o)^J \times_{\Fin_*^J} \{q_i\} \mathop{\times}\limits_{(\zO^\o)^I \times_{\Fin_*^I} \{q_i|_I\}} \{r_i\} \,\simeq \,
		\prod_{i=1}^n
	(\zO^\o)^J \times_{(\zO^\o)^I} \{r_i\} \times_{\Fin_*^J} \{q_i\},
	$$
	as desired.
\end{proof}


We now come back to the proof of proposition \ref{prop:main_step_theorem}. Using the previous lemma, since $\Delta^3$ is braced, we obtain a decomposition $d'_\zO|_{K_0} = d^-_\zO \oplus d^+_\zO$ lifting that of equation \eqref{eqn:diag_square_fin_*} and a corresponding splitting $\zZ = \zZ^- \times \zZ^+$, with 
components given by
\[
\zZ^\pm = (\zO^\o)^{\Delta^3}\times_{(\zO^\o)^{K_0}} \left\{ d^\pm_{\zO} \right\} 
\times_{\Fin_*^{\Delta^3}} \left\{ \widetilde{d}^\pm_{\Fin_*}  \right\}.
\]
As before, we note that any diagram $\Delta^3\to \zO^\o$ lifting $d^{\pm}_\zO$ will automatically be a lift of $\widetilde{d}^{\pm}_{\Fin_*}$. Therefore, the space $\zZ^\pm$ is equivalent to  $(\zO^\o)^{\Delta^3}\times_{(\zO^\o)^{K_0}} \{ d^\pm_{\zO}\} $.
We note that the arrow in $d^-_{\zO}$ that lifts the left vertical map $\id_{\langle n \rangle}$ in diagram \eqref{eqn:diag_square_fin_*} is an equivalence, since by assumption $V_0\to V_0^+$ is semi-inert. Similarly, the arrow in $d^+_{\zO}$ that lifts the right vertical map $\id_{\langle 1\rangle}$ in diagram \eqref{eqn:diag_square_fin_*} is necessarily an equivalence, by construction of the map $S_0^+\to T_0^+$.
We claim that these properties force the spaces $\zZ^-$ and $\zZ^+$ to be contractible. To see this, we need the following version of Joyal's lifting theorem.

\begin{lm}
\label{lm:generalization_joyal_thm}
Let $\zC$ be an $\infty$-category and $\alpha$ an equivalence in $\zC$. Consider an outer horn $\bar{\alpha}\colon \Lambda_0^n\to \zC$, with $n\geq 2$, whose restriction along $\Delta^1\stackrel{\langle 01\rangle}{\too}\Lambda_0^n$ is $\alpha$. Then the fiber $\zC^{\Delta^n}\times_{\zC^{\Lambda_0^n}} \{\bar{\alpha}\}$ parametrizing extensions of the form
\begin{equation}
	\label{eqn:joyal_lifting_diagram}
	\begin{tikzcd}[cramped]
	\Delta^1 \ar[r,"\langle 01\rangle" below, hook] \ar[rr,bend left, "\alpha"] & \Lambda_0^n \ar[r,"\bar{\alpha}" near start] \ar[d,hook] & \zC\\
	& \Delta^n \ar[ru,dashed] & 
	\end{tikzcd}
\end{equation}
is a contractible $\infty$-groupoid.
\end{lm}
\begin{proof}
	[Proof of lemma \ref{lm:generalization_joyal_thm}]
	We want to show that any morphism $\beta \colon \partial\Delta^m \to \zC^{\Delta^n}\times_{\zC^{\Lambda_0^n}} \{\bar{\alpha}\}$ extends to an $m$-simplex, for all $m\geq 0$.
	This problem is equivalent to finding a lift in the commutative square
	\[
	 \begin{tikzcd}  
	 	\Lambda_0^n \ar[r,"\mathrm{diag}\,\circ\, \bar{\alpha}"] \ar[d,hook,]	&	\zC^{\Delta^m} \ar[d,"i^*"] \\
	 	\Delta^n \ar[r,"\beta'" near end]\ar[ru,dashed]					& 	\zC^{\partial\Delta^m}
	 \end{tikzcd}
	 \] 
	where $\mathrm{diag}$ is the diagonal functor $\zC\to\zC^{\Delta^m}$, $i^*$ is the inner fibration induced by the inclusion $i\colon \partial\Delta^m\to\Delta^m$ and $\beta'$ is adjoint to $\partial\Delta^m\stackrel{\beta}{\to} \zC^{\Delta^n}\times_{\zC^{\Lambda_0^n}} \{\bar{\alpha}\} \to \zC^{\Delta^n}$. As $\alpha$ is an equivalence, so is $\mathrm{diag}\,\circ\,\alpha \colon  \Delta^1 \to \zC^{\Delta^m}$. Therefore, using Joyal's lifting theorem \cite[\href{https://kerodon.net/tag/019F}{Theorem 019F}]{kerodon}, we obtain the existence of lifts in the above square, as desired.
\end{proof}

From the above argument and the previous lemma, we deduce that the space $\zZ$ is contractible. This proves claim \ref{claim:Z_is_contractible}.
To finish the proof of proposition \ref{prop:main_step_theorem}, consider the pushout $\widetilde{K} = K\cup_{K_0} \Delta^{3}$ and the inclusion $\iota \colon \widetilde{K}\to \Delta^7$. The fibers $\lift^{(0)}_d$ and $\zZ$ fit in the commutative diagram
\[
\begin{tikzcd}  
	\lift^{(0)}_d \ar[r] 
	\arrow[dr, phantom, "{\lrcorner}" , very near start]
	\ar[d]			&	\zZ \ar[d] \ar[r,"\sim"] 
	\arrow[dr, phantom, "{\lrcorner}" , very near start]
	& \{d'_\zO\} \ar[d]\\
	(\zO^\o)^{\Delta^7} \ar[r,"\iota^*"]			& 	(\zO^\o)^{\widetilde{K}}\ar[r]\ar[d]
	\arrow[dr, phantom, "{\lrcorner}" , very near start]
	& (\zO^\o)^K \ar[d]\\
	& (\zO^\o)^{\Delta^3}\ar[r] & (\zO^\o)^{K_0}.
\end{tikzcd}
\]
Here, both the right outer and the bottom right squares are cartesian, therefore so is the top right square. As the top outer square is cartesian, so must be the top right square. To complete the proof that $\lift^{(0)}_d$ is contractible, we need a careful analysis of the morphism $\iota$. 
The result makes use of the notion of \emph{right anodyne morphism} recalled in the appendix in definition \ref{df:right_anodyne_morphisms} and writes as follows.

\begin{lm}
\label{lm:iota_almost_inner_anodyne}
	Consider the simplicial set $\Delta^7$ as a marked simplicial set, with $6\to 7$ as the only nondegenerate marked edge. Then the inclusion $\iota \col\widetilde{K} \too \Delta^7$ is right marked anodyne.
\end{lm}

This result is proved using a tedious explicit calculation that we defer to the end of this section, in section \ref{sub:appendix_i_is_inner_anodyne}.

Using the previous lemma and the observation that the morphism $67$ in diagram $d'|_\zO$ is an equivalence, it then follows from lemma \ref{lm:generalization_joyal_thm} that the functor $\lift^{(0)}_d\to \zZ$ induced by $\iota^*$ is an equivalence.

This concludes the proof of proposition \ref{prop:main_step_theorem}, hence that of theorem \ref{prop:lift_is_contractible}.


\subsection{Proof of technical lemmas}
\label{sec:proof_of_technical_lemmas}

In this section, we complete the proof of theorem \ref{thm:pi_is_cartesian_fibration} by providing proofs to lemmas \ref{lm:equivalence_slice_diagrams}, \ref{lm:operad_decomposition} and \ref{lm:iota_almost_inner_anodyne}.

\subsubsection{Proof of lemma \ref{lm:equivalence_slice_diagrams}} 
\label{sub:proof_of_lemma_lm:equivalence_slice_diagrams}

\begin{proof}
	[Proof of lemma \ref{lm:equivalence_slice_diagrams}]
	Let $W$ be a simplicial set. By definition of the slice $\infty$-category, we have a natural bijection
	\[
	\Hom(W,\zC_{/p})= \Hom_{{ K}_/} (W\star  K,\zC),
	\]
	where the index ${{ K}_/}$ denotes the subset of those morphisms $W\star  K\to \zC$ that restrict to $p$ on $ K$. On the other hand, we have a natural bijection
	\[
	\Hom(W,\zC^{K^\triangleleft} \times_{\zC^{ K}} \{p\})
	= \Hom_{{ K}_/} (\bar{W}, \zC),
	\]
	where we use the simplicial set
	\[
	\bar{W} = (W\times K^\triangleleft)	 \ \mathop{\amalg}\limits_{W\times  K} \  K.
	\]
	We now construct a categorical equivalence $\varphi_W \colon \bar{W}\to W\star K$, natural in $W$.
	It is obtained by the universal property of the pushout $\bar{W}$, induced by the canonical inclusion $K\to W\star K$ and a certain morphism $W\times K^{\triangleleft}\to W\star K$. To describe the latter, recall that maps from a simplicial set $X$ to $W\star K$ can be identified with triples of morphisms $(X\to \Delta^1, X_0\to W, X_1\to K)$, where $X_i = \{i\}\times_{\Delta^1}X$ is the fiber at $i$. Using this description, the morphism $W\times K^\triangleleft\to W\star K$ corresponds to the triple
	\[
	(\mathrm{can} \circ \mathrm{proj} \colon W\times K^\triangleleft \to K^\triangleleft \to \Delta^1,
	\quad
	W\times \Delta^0 \cong W,
	\quad
	\mathrm{proj}\colon W\times K\to K).
	\]
	We now prove that $\varphi_W$ is a categorical equivalence. The argument relies on the fact that the canonical morphism $c_{A,B}\colon A \diamond B\to A\star B$, comparing the two join constructions, is a categorical equivalence for all simplicial sets $A,B$ \cite[\href{https://kerodon.net/tag/01HV}{Theorem 01HV}]{kerodon}. Observe that $\bar{W}$ fits in a commutative diagram
	\[
	\begin{tikzcd}[column sep = large]
	W\times K \ar[r, hook]\ar[d,"\mathrm{proj}"] 
	\arrow[dr, phantom, "{\ulcorner}" , very near end, color=black] &
	W\times (\Delta^0\diamond K) \ar[r,"\id_W \times c_{\Delta^0,K}" below]\ar[d]
	\arrow[dr, phantom, "{\ulcorner}" , very near end, color=black] &
	W\times K^\triangleleft \ar[d] \\
	K \ar[r] &
	W\diamond K \ar[r,"\bar{c}_{W,K}"] &
	\bar{W}
	\end{tikzcd}
	\]
	in which all three squares are cocartesian. Since the top horizontal maps $W\times K \to W\times (\Delta^0\diamond K)$ and $W\times K \to W\times K^\triangleleft$ are monomorphisms, the left and outer squares are categorical pushout squares; therefore, so must be the right one. As the top right morphism $\id_W\times c_{\Delta^0,K}$ is a categorical equivalence, so is $\bar{c}_{W,K}$. Now observing that the comparison equivalence $c_{W,K} \colon W\diamond K\to W\star K$ factors as
	\[
	\begin{tikzcd}
	W\diamond K \ar[r,"{\bar{c}_{W,K}}"] &
	\bar{W} \ar[r,"{\varphi_W}"] &
	W\star K,
	\end{tikzcd}
	\]
	we conclude that $\varphi_W$ is a categorical equivalence.
	
	The morphisms $\varphi_W$ for varying $W$ induce a functor of $\infty$-categories $\varphi \colon \zC_{/p} \to \zC^{K^\triangleleft} \times_{\zC^{ K}} \{p\}$. 
	To prove that $\varphi$ is an equivalence, we will show that for each simplicial set $W$, the induced morphism
	\begin{equation}
		\label{eqn:varphi_comparison_slice}
	\varphi_* \colon 
	\pi_0\left(\Fun(W,\zC_{/p}) ^{\simeq}\right) \too
	\pi_0\left(\Fun(W,\zC^{K^\triangleleft} \times_{\zC^{ K}} \{p\}) ^{\simeq}\right)
	\end{equation}
	is a bijection.
	By \cite[\href{https://kerodon.net/tag/01KV}{Tag 01KV}]{kerodon}, we have a canonical bijection
	\[
	\pi_0\left(\Fun(W,\zC_{/p}) ^{\simeq}\right) \cong
	\pi_0\left(\Fun_{{ K}_/}(W\star K,\zC) ^{\simeq}\right).
	\]
	We claim that one has a similar bijection for the target of $\varphi_*$, namely:\vspace*{3pt}\\
	{\bfseries Claim.} There is a canonical bijection
	\begin{equation}
	\label{eqn:claim_proof_lm_slice}
	\pi_0\left(\Fun\left(W,\zC^{K^\triangleleft} \times_{\zC^{ K}} \{p\} \right)^{\simeq}\right)
	\cong 
	\pi_0\left(\Fun_{{ K}_/}(\bar{W},\zC) ^{\simeq}\right).
	\end{equation}
	Assuming this claim, we see that the morphism $\varphi_*$ from \eqref{eqn:varphi_comparison_slice} corresponds, under the previous two bijections, to the map 
	\[
	\varphi_W^* \colon
	\pi_0\left(\Fun_{{ K}_{/}}(W\star K,\zC)^{\simeq}\right) {\too}
	\pi_0\left(\Fun_{{ K}_{/}}(\bar{W},\zC)^{\simeq}\right),
	\]
	induced by $\varphi_W$. By assumption, $\varphi_W$ is a categorical equivalence compatible with restricting to $ K$; thus $\varphi_W^*$ is a bijection, for all $W$, as desired.

	To complete the proof, we now prove claim \eqref{eqn:claim_proof_lm_slice}. Let $\alpha_0, \alpha_1$ be two functors $W\to \zC^{K^\triangleleft} \times_{\zC^{ K}} \{p\}$ and let $\bar{\alpha}_0, \bar{\alpha}_1$ denote the corresponding objects in $\Fun_{{ K}_/}(\bar{W},\zC)$ under the bijection $
	\Hom\left(W,\zC^{K^\triangleleft} \times_{\zC^{ K}} \{p\} \right)
	\cong 
	\Hom_{{ K}_/}(\bar{W},\zC)$. We wish to prove that $\alpha_0$ and $\alpha_1$ are equivalent if and only if $\bar{\alpha}_0$ and $\bar{\alpha}_1$ are. 
	To this end, we will use a characterization of equivalences in functor categories.
	Consider a categorical mapping cylinder of $W$, that is a factorization of $(\id_W,\id_W)$ of the form
	\begin{equation}
	\label{eqn:categorical_mapping_cylinder}	
	\begin{tikzcd}
	W\amalg W \ar[r,"{(s_0,s_1)}", hook] &
	RW \ar[r,"\sim" below, "p" above] &
	W
	\end{tikzcd}
	\end{equation}
	where $p$ is a categorical equivalence and $(s_0,s_1)$ is a monomorphism. From \cite[\href{https://kerodon.net/tag/01KD}{Corollary 01KD}]{kerodon}, we know that the objects $\alpha_0$ and $\alpha_1$ are equivalent if and only if the following condition is satisfied:
	\begin{enumerate}
		\item there exists $\alpha \colon RW\to \zC^{K^\triangleleft}\times_{\zC^K} \{p\}$ such that $\alpha\circ s_0 = f_0$ and $\alpha\circ s_1 = f_1$.
	\end{enumerate}
	By definition of the functor $X\mapsto \bar{X}$, one observes that the latter condition is equivalent to the following:
	\begin{enumerate}
		\setcounter{enumi}{1}
		\item there exists some $\alpha'\colon \bar{RW} \to \zC$ satisfying $\alpha' \circ \bar{s_0} = \bar{f}_0$ and $\alpha' \circ \bar{s_1} = \bar{f}_1$. 
	\end{enumerate}
	Using again \cite[\href{https://kerodon.net/tag/01KD}{Corollary 01KD}]{kerodon}, one sees that this last condition is verified if and only if the objects $\bar{f}_0$ and $\bar{f}_1$ are equivalent, provided that the factorization	
	\begin{equation}
	\begin{tikzcd}
	\bar{W}\amalg_K \bar{W} \ar[r,"{(\bar{s}_0,\bar{s}_1)}"] &
	\bar{RW} \ar[r, "\bar{p}" above] &
	\bar{W}
	\end{tikzcd}
	\end{equation}
	is a categorical mapping cylinder for $\bar{W}$ relative to $K$; the proof of this last fact is an easy verification. This proves claim \eqref{eqn:claim_proof_lm_slice} and finishes the proof of lemma \ref{lm:equivalence_slice_diagrams}.

\end{proof}

\subsubsection{Proof of lemma \ref{lm:operad_decomposition}}
\label{sec:proof_lemma_operad_decomposition}

\begin{proof}[Proof of lemma \ref{lm:operad_decomposition}]
	For ease of notation, let $\zD$ denote the left hand side of \eqref{eqn:lm_key_equivalence}. 
	Consider a diagram $X\in \zD$.
	For each vertex $j\in J$, 
	we can write the object $X(j)$ in the form $\oplus_i X(j)_i$, with the resulting diagram $X_i$ lying over $q_i$. Then, for each $1$-simplex $f\col j_0\to j_1$ of $J$, using the definition of $\infty$-operads, we see that the space $\Map_{\zO}^{q(f)}(X(j_0), X(j_1))$ decomposes canonically as $\prod_i \Map_{\zO}^{q_i(f)}(X(j_0)_i, X(j_1)_i)$. Therefore, up to equivalence, we can write the morphism $X(f)$ as a disjoint union $\oplus_i X(f)_i$ with each component lying over $q_i(f)$.

	We now consider the general case. Let $\sigma$ be a simplex of $J$ of dimension $n$ and let $\alpha \col \mathrm{Sp}^n \inc \Delta^n$ denote the spine inclusion. Since $\mathrm{Sp}^n$ is $1$-skeletal, the previous part of the proof gives a decomposition of the space 
	parametrizing diagrams $\mathrm{Sp}^n\to \zO^\o$ lifting $q(\sigma\circ\alpha)$
	as a product
	\begin{equation}
	\label{eqn:proof_lm_key_equiv_sp^n}
		\prod_{i=1}^n
		(\zO^\o)^{\mathrm{Sp}^n} \times_{\Fin_*^{\mathrm{Sp}^n}} \{q_i(\sigma\circ\alpha)\}
		\stackrel{\simeq}{\too}
		(\zO^\o)^{\mathrm{Sp}^n} \times_{\Fin_*^{\mathrm{Sp}^n}} \{q(\sigma\circ\alpha)\},
	\end{equation}
	with the equivalence given by disjoint union.
	Now, as $\alpha$ is anodyne, we have a canonical equivalence between the spaces of diagrams $(\zO^\o)^{\Delta^n} \wk (\zO^\o)^{\mathrm{Sp}^n}$, from which we can extend \eqref{eqn:proof_lm_key_equiv_sp^n} to an equivalence
	\begin{equation}
		\label{eqn:proof_lm_key}
		\prod_{i=1}^n
		(\zO^\o)^{\Delta^n} \times_{\Fin_*^{\Delta^n}} \{q_i(\sigma)\}
		\stackrel{\simeq}{\too}
		(\zO^\o)^{\Delta^n} \times_{\Fin_*^{\Delta^n}} \{q(\sigma)\}.
	\end{equation}
	\begin{claim}
        \label{claim:proof_operad_decomposition}
		Using the previous equivalence for every simplex $\sigma$ of $I$ and $J$, we obtain two decompositions
		\begin{equation}
			\label{eqn:decompositions_E^I_fiber}
		\prod_{i=1}^n (\zO^\o)^I\times_{\Fin_*^I} \{q_i|_I\}\simeq (\zO^\o)^I\times_{\Fin_*^I} \{q|_I\}
		,\qquad
		\prod_{i=1}^n (\zO^\o)^J\times_{\Fin_*^J} \{q_i\}\simeq (\zO^\o)^J\times_{\Fin_*^J} \{q\}.
		\end{equation}
	\end{claim}
	\begin{proof}
		[Proof of the claim]
		We prove the result for $J$, the case of $I$ being completely similar.
		Consider the category $\Delta \downarrow J$ of simplices of $J$ and let $\iota$ denote the inclusion $(\Delta\downarrow J)_{\mathrm{nd}}\subset \Delta\downarrow J$ of the full subcategory consisting of all nondegenerate simplices. Writing $\cF$ and $\cF_i$ for the simplicial presheaves 
		$\cF(\sigma) = (\zO^\o)^{\Delta^n} \times_{\Fin_*^{\Delta^n}} \{q(\sigma)\}$ and 
		$\cF_i(\sigma)=(\zO^\o)^{\Delta^n} \times_{\Fin_*^{\Delta^n}} \{q_i(\sigma)\}$ on $\Delta \downarrow J$,
		the natural equivalences \eqref{eqn:proof_lm_key} give a natural transformation $\gamma \colon \prod_i\cF\Rightarrow \cF$ which is a levelwise categorical equivalence. 
		Our goal is to show that the equivalence \eqref{eqn:decompositions_E^I_fiber} is obtained from $ \iota^*\gamma \colon \iota^*\prod_i\cF_i \Rightarrow \iota^*\cF$ by taking the colimit over $(\Delta\downarrow J)_{\mathrm{nd}}$.

		First, we show that $\iota^*\prod_i\cF_i$ and $\iota^*\cF$ are isofibrant diagrams, in the sense of \cite[\href{https://kerodon.net/tag/0349}{Tag 0349}]{kerodon}. We begin by proving that the forgetful functor $\cU\colon (\Delta\downarrow J)_{\mathrm{nd}}\to \sSet$ sending $\Delta^n\to J$ to $\Delta^n$ is projectively cofibrant (that is cofibrant in the projective global model structure on the diagram category $\Fun((\Delta\downarrow J)_{\mathrm{nd}},\sSet)$). Observe that each simplicial level $\cU_k$ of $\cU$ decomposes as a coproduct
		$$
		\cU_k = \cU^\mathrm{nd}_k \coprod \cU^{\mathrm{deg}}_k
		$$
		of subfunctors of nondegenerate and degenerate simplices respectively. Each of these two subfunctors is a coproduct of representables, so that both are projectively cofibrant. By \cite[Corollary 9.4]{dugger_universal_homotopy_theories}, we deduce that $\cU$ is projectively cofibrant. Therefore, the simplicial presheaves $\Hom(\cU,\zO^\o)$ and $\Hom(\cU,\Fin_*)$ on $(\Delta\downarrow J)_{\mathrm{nd}}$ are isofibrant. Taking pullback, it follows that the simplicial presheaf $\cF$ is also isofibrant. The same argument proves that $\prod_i\cF_i$ is isofibrant.

		Since $\iota^*\gamma$ is a levelwise categorical equivalence between isofibrant diagrams, the induced map $\lim(\iota^*\gamma) \colon \lim \prod_i\iota^*\cF_i \to \lim \iota^*\cF$ on limits is a categorical equivalence. The only remaining step is to identify $\lim \prod_i\iota^*\cF_i$ (respectively $\lim \iota^*\cF$) with the left (resp. right) hand side of the second equivalence of \eqref{eqn:decompositions_E^I_fiber}. To this end, note that since $J$ is braced by assumption, the inclusion $\iota$ admits a left adjoint and therefore is cofinal. This implies that $J \cong \colim \cU \cong \colim \iota^*\cU$, from which the result follows easily.
	\end{proof}
    The proof of lemma \ref{lm:operad_decomposition} now follows from claim \ref{claim:proof_operad_decomposition}.
\end{proof}
\todo{[SEMBLE OK] vérifier le changement E vs O à la relecture}

\subsubsection{Proof of lemma \ref{lm:iota_almost_inner_anodyne}}
\label{sub:appendix_i_is_inner_anodyne}

\begin{proof}
	[Proof of lemma \ref{lm:iota_almost_inner_anodyne}]
	Recall that the morphism $\iota$ is the inclusion
	\[
	\iota \colon \widetilde{K} = \left(\Delta^{013467}\bigcup_{\Delta^{0347}}\Delta^{023457}\right)
\bigcup_{\Delta^{013}\mathop{\cup}\limits_{\Delta^{03}} \Delta^{023}} \Delta^{0123}
	\too \Delta^7
	\]
	using the numbering introduced in diagram \eqref{eqn:diag_7_simplex}.
	First, note that $\iota$ factors through $\hat{K} = \widetilde{K}\mathop{\cup}\limits_{\Lambda^{567}_7}\Delta^{567}$ and the right anodyne inclusion $\widetilde{K}\to\hat{K}$ satisfies the conditions of the statement. It therefore suffices to show that the induced map $\hat{K}\to \Delta^7$ is inner anodyne. Note that the spine inclusion $\mathrm{Sp}^7\to \Delta^7$, which is inner anodyne, factors through $\hat{K}$. By the right cancellation property for inner anodyne morphisms, it suffices to show that $\mathrm{Sp}^7\to \hat{K}$ is inner anodyne. We decompose the latter inclusion in several steps; first, it is easy to see that the two morphisms
	\[
	\mathrm{Sp}^7\to 
	\mathrm{Sp}^7\mathop{\cup}\limits_{\mathrm{Sp}^{0123}}\Delta^{0123} \to 
	\left(
	\mathrm{Sp}^7\mathop{\cup}\limits_{\mathrm{Sp}^{0123}}\Delta^{0123}
	\right)
	\mathop{\cup}\limits_{\mathrm{Sp}^{4567}}
	\Lambda_{7}^{4567}
	=: S
	\]
	are inner anodyne. As the inclusion $\mathrm{Sp}^7\to\Delta^7$ factors through the composite map $\mathrm{Sp}^7\to S$, the remaining steps consists in adding to $S$ the simplices $\Delta^{013467}$ and $\Delta^{023457}$. One verifies that the intersection between $S$ and $\Delta^{013467}$ is given by
	\[
	S \cap \Delta^{013467} = \Delta^{013}\mathop{\cup}\limits_{\Delta^{\{3\}}} \Delta^{34} \mathop{\cup}\limits_{\Delta^{\{4\}}} \Delta^{467}
	\]
	so that the inner anodyne inclusion $\mathrm{Sp}^{013467}\to \Delta^7$ factors through $S \cap \Delta^{013467}$ as an inner anodyne map. Therefore, by right cancellation, we deduce that $S \cap \Delta^{013467}\to \Delta^{013467}$ is inner anodyne.
	Since the inclusion of $S$ into the simplicial set $\hat{S} =  S\mathop{\cup}\limits_{S\cap\Delta^{013467}} \Delta^{013467}$ is the pushout of the inclusion $S \cap \Delta^{013467}\to \Delta^{013467}$, it is also inner anodyne. Finally, we will prove that the same holds for the inclusion $\hat{S}\to\hat{S} \cup \Delta^{023567} = \hat{K}$. We proceed as before: the intersection of $\hat{S}$ and $\Delta^{023457}$ is given by
	\[
	\hat{S} \cap \Delta^{023457} = 
	\Delta^{023}\mathop{\cup}\limits_{\Delta^{\{3\}}}
	\Delta^{34}
	\mathop{\cup}\limits_{\Delta^{\{4\}}}
	\Delta^{457}
	\]
	and the inclusion $\mathrm{Sp}^{023457}\to\Delta^{023457}$ factors through it. To prove that $\hat{S}\to\hat{K}$ is inner anodyne, it is then enough to show that $\mathrm{Sp}^{023457}\to \hat{S} \cap \Delta^{023457}$ has this property, which is obvious. 
\end{proof}


	\section{Comparison with Lurie's spaces of extensions} 
\label{sec:comparison_with_lurie_s_spaces_of_extensions}

We now turn to the problem of comparison between Mann--Robalo's and Lurie's models for spaces of extensions, which we solve through theorem \ref{thm:comparaison_bo_ext}.

\subsection{Statement of the results} 
\label{sec:statement_of_the_results_cha_comparison_Lurie}

Let $\zO^\o$ be a unital $\infty$-operad and fix an active morphism $\sigma$, viewed as an object in the twisted arrow $\infty$-category $\zT := \Tw(\Env(\zO))^\o$ of the monoidal envelope of $\zO^\o$.

The goal of the present section is to prove the equivalence between the fiber $\bo_\sigma$ of the brane fibration $\bo\to \zT$ (introduced in definition \ref{df:bo}) at the operation $\sigma$ coincide with the $\infty$-category of extensions $\Ext(\sigma)$ of $\sigma$, thereby proving the corresponding claim in \cite{robalo}
\footnote{
    More precisely, for the claim $\bo_\sigma\simeq\Ext(\sigma)$ to be correct, one needs to adopt definition \ref{df:Ext(sigma)} as one's definition of the $\infty$-category of extensions, instead of \cite[Definition 3.3.1.4.]{ha} used in \cite{robalo}.
This minor change is argued in remark \ref{rk:difference_df_Ext_with_HA}.
}.
Our method consist in providing an explicit zigzag of categorical equivalences between $\bo_\sigma$ and $\Ext(\sigma)$.

\todo{[OK] Vérifier les hypothèses sur $\zO^\o$ : $\infty$-groupoïde (pas besoin), unitalité (osef), naturalité des équivalences (?!)}

\begin{thm}[Theorem \ref{thm:comparaison_bo_ext_intro}]
	\label{thm:comparaison_bo_ext}
	Let $\sigma$ be an active morphism in a unital $\infty$-operad $\zO^\o$. Then the fiber $\bo_\sigma$ of the brane fibration and the $\infty$-category of extensions $\Ext(\sigma)$ are equivalent in $\Cat_\infty$.
\end{thm}

\begin{cor}
	\label{cor:comparaison_bo_ext_case_Kan_complex}
	Let $\sigma$ be an active morphism in a unital $\infty$-operad $\zO^\o$, whose underlying $\infty$-category is an $\infty$-groupoid. Then $\bo_\sigma$ and $\Ext(\sigma)$ are equivalent Kan complexes.
\end{cor}

The difference between the simplicial sets $\Ext(\sigma)$ and $\bo_\sigma$ can be observed from their set of $0$-simplices. Both sets parametrize diagrams in $\zO^\o$, respectively of the form
\begin{equation}
	\label{eqn:ext_vs_bo}
	\begin{tikzcd}  
		\bullet \ar[r,"\sigma"] \ar[d,hook,harpoon]\ar[dr]	&	\bullet \ar[d,"\wr"] \\
		\bullet \ar[r,""]					& 	\bullet
	\end{tikzcd}
	\qquad \text{and} \qquad
	\begin{tikzcd}  
		\bullet \ar[r,"\sigma"] \ar[d,hook,harpoon]\ar[dr]	&	\bullet \\
		\bullet \ar[r,""]		\ar[ru, crossing over]			& 	\bullet, \ar[u,"\wr" right]
	\end{tikzcd}
\end{equation}
that is, of respective shape $\Delta^1\times \Delta^1$ and $\Delta^3$,
such that the following conditions are satisfied:
\begin{itemize}
	\item the top horizontal arrow is sent to $\sigma$;
	\item the left vertical arrow is sent to an atomic map; 
	\item the right vertical arrow is sent to an equivalence in $\zO^\o$.
\end{itemize}

Informally, since the right vertical arrow is marked as an equivalence, both diagrams in \eqref{eqn:ext_vs_bo} encode the same data, namely that of a triangle of the form
\begin{equation}
	\begin{tikzcd}  
		\bullet \ar[r,"\sigma"] \ar[d,hook,harpoon]	&	\bullet \\
		\bullet \ar[ru]					& 
	\end{tikzcd}
\end{equation}
Our proof follows this idea and relies on finding suitable generalizations of the previous diagram for higher dimensional simplices in $\Ext(\sigma)$ and $\bo_\sigma$.




\subsection{Definition of spaces interpolating between \texorpdfstring{$\bo_\sigma$}{BO sigma} and \texorpdfstring{$\Ext(\sigma)$}{Ext(sigma)}} 
\label{sec:definition_of_spaces_interpolating}

We first discuss the shape of diagrams represented by general simplices of $\Ext(\sigma)$ and $\bo_\sigma$. Let $K$ be a simplicial set. 

\begin{figure}[!ht]
    \label{fig:zigzag_2}
    \centering
	\begin{eqnarray*}
	\begin{tikzcd}[scale cd = .8, column sep = tiny, row sep = small]
	{\triangleleft} \arrow[rddd, hook,harpoon, crossing over] \arrow[dd, hook,harpoon] \arrow[rrr, "\sigma"] &                                     &  & {\bar{\triangleleft}} \arrow[rddd, "\wr"] \arrow[dd, "\wr"] &           \\  
	& &  & &   \\
	{k_0} \arrow[rrr] \arrow[rd]                                &                                     &  & {\bar{k_0}} \arrow[rd]                                      &           \\
	& {k_1} \arrow[rrr] &  &                                                           & {\bar{k_1}}
	\end{tikzcd}
	&
    \stackrel{i_0}{\too}
	&
	\begin{tikzcd}[scale cd = .8, column sep = tiny, row sep = small]
	{\triangleleft} \arrow[rddd, hook,harpoon, crossing over] \arrow[dd, hook,harpoon] \arrow[rrr, "\sigma"] &                                     &  & {\bar{\triangleleft}} \arrow[rddd, "\wr"] \arrow[dd, "\wr"] &           \\
	& &  &&           \\
	{k_0} \arrow[rrr] \arrow[rrruu] \arrow[rd]                                &                                     &  & {\bar{k_0}} \arrow[rd,"\sim" left]                                      &           \\
	& {k_1} \arrow[rrr] \arrow[rruuu, crossing over] &  & & {\bar{k_1}}
	\end{tikzcd}
    \stackrel{i_1}{\lla}
	\begin{tikzcd}[scale cd = .8, column sep = tiny, row sep = small]
	\triangleleft \arrow[rddd, hook,harpoon] \arrow[dd, hook,harpoon] \arrow[rrr, "\sigma"] &                   &  & \triangleright \\
	&                   &  &                \\
	k_0 \arrow[rrruu] \arrow[rd]                                            &                   &  &                \\
	 & k_1 \arrow[rruuu] &  &               
	\end{tikzcd}\\
	&
    \stackrel{i_2}{\too}
	&
    \begin{tikzcd}[scale cd = .8, column sep = tiny, row sep = small]
    \triangleleft \arrow[rddd, hook,harpoon] \arrow[dd, hook,harpoon] \arrow[rrr, "\sigma"] &                 &  & \overline{\triangleleft}         &                                               \\
     &                 &  &                                  &                                               \\
    k_0 \arrow[rrr] \arrow[rd]                                              &                 &  & \overline{k_0} \arrow[uu, "\wr"] &                                               \\
     & k_1 \arrow[rrr] &  &                                  & \overline{k_1} \arrow[lu,"\sim"] \arrow[luuu, "\wr"]
\end{tikzcd}
	\stackrel{p}{\lla}
	\begin{tikzcd}[scale cd = .8, column sep = tiny, row sep = tiny]
	{(k_0,0)} \arrow[rddd, hook,harpoon, crossing over] \arrow[dd, hook,harpoon] \arrow[rrr, "\sigma", crossing over] \arrow[rd, Rightarrow, no head] &                                            &  & {\overline{(k_0,0)}}  \arrow[from=dd, "\wr" near start] \arrow[from=dddr]                 &                                                                       \\
		& {(k_1,0)} \arrow[dd, hook,harpoon, crossing over] \arrow[rrr, "\sigma", near start, crossing over] &  & & {\overline{(k_1,0)}} \arrow[lu, Rightarrow, no head]                  \\
	{(k_0,1)} \arrow[rrr] \arrow[rd] & &  & {\overline{(k_0,1)}}  &                                                                       \\
		& {(k_1,1)} \arrow[rrr]                      &  & & {\overline{(k_1,1)}} \arrow[uu, "\wr"] \arrow[lu] 
	\end{tikzcd}
	\end{eqnarray*}
	\caption{Zigzag diagram \eqref{eqn:zigzag_F_n^+(K)} for $K = \{k_0\to k_1\} \cong \Delta^1$}
\end{figure}
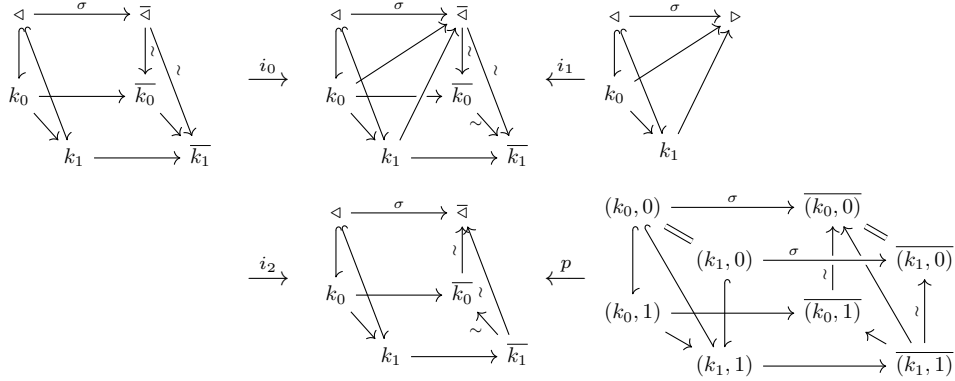

To help the reader with the cumbersome definitions and notations, we recommend to look at figures \ref{fig:zigzag_1} and \ref{fig:zigzag_2} where the different diagrams introduced in this section are depicted, for the cases $K = \Delta^0$ and $K = \Delta^1$.

\subsubsection{Diagrams indexed by \texorpdfstring{$\Ext(\sigma)$}{Ext(sigma)}}
Let $F_0(K)$ denote 
the simplicial set $K^\triangleleft \times \Delta^1$ and consider the canonical projection $F_0(K) \to \Delta^1$. The vertices in the fiber over $0$ will be denoted as $x$, whereas those in the fiber over $1$ will be written $\bar{x}$, where $x\in K^\triangleleft$. Unravelling definition \ref{df:Ext(sigma)}, we see that
morphisms $K\to \Ext(\sigma)$ are identified with those diagrams $\alpha \colon F_0(K) \to \zO^\o$  that send all morphisms in $K\times\{1\}\subset K^\triangleleft \times \Delta^1$ to equivalences and moreover satisfy the following condition.

\vspace{5pt}
\hypertarget{text:condition_F_0}{
\textbf{{Condition ${(\star)_{0,\sigma}}$}}:}
for every vertex $k\in K$, 
\begin{itemize}
	\item the morphism $\triangleleft \to k$ in $F_0(K)$ is sent to an atomic map;
	\item the morphism $k \to \bar{k}$ in $F_0(K)$ is sent to an active map;
    \item the restriction of $\alpha$ to $\{\triangleleft\}\times \Delta^1$ is $\sigma$; and
    \item for every morphism $k_0\to k_1$ in $K$, the corresponding morphism in $F_0(K)$ is compatible with extensions.
\end{itemize}
Letting $F_0^+(K)$ denote the marked simplicial set obtained from $F_0(K)$ by further marking the $1$-simplices of $K^\triangleleft \times \{1\}$, one obtains a bijection between $\Hom(K,\Ext(\sigma))$ and a subset $\Hom^{\sigma}(F_0^+(K),\zO^\o)$ of $\Hom_{\sSet^+}(F_0^+(K), \zO^{\o,\natural})$ given by those diagrams that satisfy condition \hyperlink{text:condition_F_0}{$(\star)_{0,\sigma}$}.

\subsubsection{Diagrams indexed by \texorpdfstring{$\bo_\sigma$}{BO sigma}}
\label{sub:diag_bo_sigma}

Observe that the set of diagrams from $K$ to $\bo_\sigma$ is a certain subset of
\[
	\Hom(K,\zT^{\Delta^1} \times_{\zT} \{\sigma\})
	= \Hom(K\times \Delta^1, \zT) \mathop{\times}\limits_{\Hom(K,\zT)} \{\sigma \circ \mathrm{proj}_K\}.
\]
Using remark \ref{rk:remarks_on_objects_in_bo}, we may then identify $\Hom(K,\bo_\sigma)$ with a particular subset of $\Hom(G(K), \zO^\o)$, where $G(K)$ denotes the pushout 
$$ G(K) := s_*(K\times \Delta^1) \mathop{\amalg}\limits_{s_*(K\times\{0\})} s_*\{0\}.$$
\sloppy
Let $G^+(K)$ denote the marked simplicial set obtained from $G(K)^\flat$ by further marking the arrows of the form $\bar{(k,1)}\to \bar{(k,0)}$, for all vertices $k$ in $K$. Unraveling the definition of $\bo$, we will identify $\Hom(K,\bo_\sigma)$ with the subset $\Hom^\sigma(G^+(K), \zO^{\o,\natural})$ of $\Hom_{\sSet^+}(G^+(K), \zO^{\o,\natural})$ consisting of those diagrams verifying the following condition:

\hypertarget{text:condition_bo_sigma}{
\textbf{{Condition ${(\star)_{G,\sigma}}$}:}}
\begin{itemize}
	\item the arrow $s_*\{0\}\cong \Delta^1$ in $G^+(K)$ is sent to $\sigma$;
	\item for every $k\in K$, the $3$-simplex $s_*(\{k\}\times\Delta^1)$ given by the vertices $(k,0)(k,1)\bar{(k,1)}\bar{(k,0)}$ 
	is sent to an object of $\bo$;
	\item for every $1$-simplex $f\col k\to k'$ in $K$, the corresponding diagram $s_*(f\times \Delta^1)$ is sent to a morphism of $\bo$; more explicitly, this conditions means that the arrow $f\times \{1\} \col (k,1)\to(k',1)$ in $G^+(K)$ is sent to a map in $\zO^\o$ that is compatible with extensions (in the sense of definition \ref{df:bo}).
\end{itemize}

\subsubsection{Definition of intermediate steps}
\label{sub:diag_intermediate_steps}
We will relate $\Ext(\sigma)$ and $\bo_\sigma$ using intermediate $\infty$-categories $\zC_i$, for $i\in \{1,2,3\}$, whose construction is given by the following procedure. 
For $K$ a simplicial set, the set $\Hom(K,\zC_i)$ is identified with the subset $\Hom^\sigma(F_i^+(K),\zO^{\o,\natural})$ of $\Hom_{\sSet^+}(F_i^+(K),\zO^{\o,\natural})$ consisting of all diagrams $F_i^+(K)\to \zO^{\o,\natural}$ that satisfy a certain condition denoted $(\star)_{i,\sigma}$. Here, the marked simplicial sets $F_i^+(K)$ are to be thought of as shapes which interpolate between $F_0^+(K)$ and $G^+(K)$. These marked simplicial sets will fit into a zigzag of the form
\begin{equation}
	\label{eqn:zigzag_F_n^+(K)}
	\begin{tikzcd}
		F_0^+(K) 
		\ar[r,"i_0(K)"] &
		F_1^+(K) &
		F_2^+(K) \ar[l,"i_1(K)" above]\ar[r,"i_2(K)"]&
		F_3^+(K) &
		G^+(K), \ar[l,"p(K)" above]
	\end{tikzcd}
\end{equation}
natural in $K$, which yields a zigzag of functors of $\infty$-categories
\begin{equation}
	\label{eqn:zigzag_comparaison_categories}
	\begin{tikzcd}
		\Ext(\sigma) &
		\zC_1 \ar[l,"i_0^*" above]\ar[r,"i_1^*"] &
		\zC_2 &
		\zC_3 \ar[l,"i_2^*" above]\ar[r,"p^*" above]&
		\bo_\sigma .
	\end{tikzcd}
\end{equation}

\begin{nota}
	Unless ambiguous, we will write $i_n$ instead of $i_n(K)$.
\end{nota}

We now describe the different marked simplicial sets $F_i^+(K)$, their associated condition $(\star)_{i,\sigma}$ and the morphisms $i_0$, $i_1$, $i_2$  and $p$.

\begin{enumerate}
    \item[$(F_1^+)$]
        The marked simplicial set $F_1^+(K)$ is obtained from the simplicial set $$F_1(K) = \mathop{\colim}\limits_{\Delta^m\to K^\triangleleft} (\Delta^m \ast {\Delta^m}),$$ by marking all edges of the second copy of $\Delta^m$. 
    \item[$(i_0)$] Using that $F_0(K)$ can be rewritten as $\mathop{\colim}\limits_{\Delta^m \to K^\triangleleft} (\Delta^m \times \Delta^1)$, the canonical inclusions $\Delta^m \times \Delta^1 \to \Delta^m\ast \Delta^m$ induce a map of marked simplicial sets $i_0 \col F_0^+(K) \to F_0^+(K)$. Similarly to the case of $F_0(K)$, we label vertices of $F_1(K)$ as $x$ or $\bar{x}$ using the obvious projection $F_1(K) \to \Delta^1$, where $x\in K^\triangleleft$.

    \item[$(F_2^+)$] 
        The marked simplicial set $F_2^+(K)$ is defined as $F_2(K)^\flat$, where 
        $$F_2(K) = K^{\triangleleft\triangleright}.$$
    \item[$(i_1)$] Using the isomorphism $F_2(K) \cong \mathop{\colim}\limits_{\Delta^m\to K^\triangleleft} (\Delta^m)^{\triangleright}$, we define the inclusion $i_1 \col F_2^+(K) \to F_1^+(K)$ as induced by the morphisms $i_{1,m}\colon (\Delta^m)^{\triangleright} \to F_1^+(K)$ sending $\Delta^m$ to the first copy of itself in $\Delta^m\ast\Delta^m$ and $\triangleright$ to $\bar{\triangleleft}$.

    \item[$(F_3^+)$] 
        The marked simplicial set $F_3^+(K)$ is obtained from the simplicial set
        \[
            F_3(K) = s_*(K^{\triangleleft}) \cong \mathop{\colim}\limits_{\Delta^m\to K^\triangleleft} (\Delta^m\ast \Delta^{m,\mathrm{op}})
        \]
        by marking all edges of the second copy $\Delta^{m,\mathrm{op}}$.

    \item[$(i_2)$]
        The morphism $i_2\colon F_2^+(K)\to F_3^+(K)$ is induced by the inclusions $(\Delta^m)^{\triangleright}\to F_3^+(K)$, that send $\triangleright$ to $\bar{\triangleleft}\in F_3^+(K)$ and $\Delta^m$ to the first copy of itself in $\Delta^m\ast\Delta^{m,\mathrm{op}}$.

    \item[$(p)$]
We now define the morphism $p\colon G^+(K)\to F_3^+(K)$. 
We will make use of the canonical isomorphism
$$s_*(K\times\Delta^1)\cong \mathop{\colim}\limits_{\Delta^n\to K} s_*(\Delta^n\times\Delta^1).$$
Given an $n$-simplex $\tau_n \colon \Delta^n\to K$ of $K$, consider 
the projection $\tilde{p}_n\colon \Delta^n\times \Delta^1\to (\Delta^n)^{\triangleleft}$ that sends the vertex $(k,1)$ to $k$ and all vertices of the form $(i,0)$ to $\triangleleft$. The composition $\tau_n\circ \tilde{p}_n\colon \Delta^n\times \Delta^1\to K^\triangleleft$ yields a map
$$
s_* \tilde{p}_n \colon s_*(\Delta^n\times \Delta^1)\to s_*\left((\Delta^n)^\triangleleft\right) \subseteq
\mathop{\colim}\limits_{\Delta^m\to K^\triangleleft} s_*(\Delta^m) = F_3(K).
$$
For varying $n$-simplices $\tau_n$, the maps $s_*\tilde{p}_n$  organize into a map $\tilde{p}\colon s_*(K\times \Delta^1)\to F_3^+(K)$, whose restriction to $s_*(K\times \{0\})$ factors through $s_*\{0\}$ as the inclusion $s_*\{0\}\cong \Delta^{\triangleleft \bar{\triangleleft}} \subseteq F_3^+(K)$ and therefore induces our map $p\colon G(K)\to F_3(K)$.

\item[$(\star)_{n,\sigma}$] \hypertarget{text:condition_F_n}{The conditions $(\star)_{n,\sigma}$ for $n\in \{1,2,3\}$} are given \emph{mutatis mutandis} by condition \hyperlink{text:condition_F_0}{$(\star)_{0,\sigma}$}.
\end{enumerate}



\begin{figure}[!ht]
    \label{fig:zigzag_1}
    \centering
	\begin{equation*}
        \begin{tikzcd}[scale cd = 0.9]
	{\triangleleft} \arrow[r, "\sigma"] \arrow[d, hook,harpoon] \arrow[rd] & {\bar{\triangleleft}} \arrow[d, "\wr"] \\
	{k_0} \arrow[r]                                                & {\bar{k_0}}                           
	\end{tikzcd}
	\!\xrightarrow{i_0}\!
	\begin{tikzcd}[scale cd = 0.9]
	{\triangleleft} \arrow[r, "\sigma"] \arrow[d, hook,harpoon] \arrow[rd] & {\bar{\triangleleft}} \arrow[d, "\wr"] \\
	{k_0} \arrow[r] \arrow[ru, crossing over]                                     & {\bar{k_0}}                           
	\end{tikzcd}
	\!\xleftarrow{i_1}\!
    \begin{tikzcd}[scale cd = 0.9]
	\triangleleft \arrow[r, "\sigma"] \arrow[d, hook,harpoon] & \triangleright \\
	k_0 \arrow[ru]                                    &               
	\end{tikzcd}
	\!\xrightarrow{i_2}\!
	\begin{tikzcd}[scale cd = 0.9]
	{\triangleleft} \arrow[r, "\sigma"] \arrow[d, hook,harpoon] \arrow[rd] & {\overline{\triangleleft}}                  \\
	{k_0} \arrow[ru, crossing over] \arrow[r]                           & {\overline{k_0}} \arrow[u, "\wr"]
	\end{tikzcd}
\!\xleftarrow{p}\!
	\begin{tikzcd}[scale cd = 0.9]
	{(k_0,0)} \arrow[r, "\sigma"] \arrow[d, hook,harpoon] \arrow[rd] & {\overline{(k_0,0)}}                  \\
	{(k_0,1)} \arrow[ru, crossing over] \arrow[r]                           & {\overline{(k_0,1)}} \arrow[u, "\wr"]
	\end{tikzcd}
	\end{equation*}
	\caption{Zigzag diagram \eqref{eqn:zigzag_F_n^+(K)} for $K = \{k_0\} \cong \Delta^0$}
\end{figure}
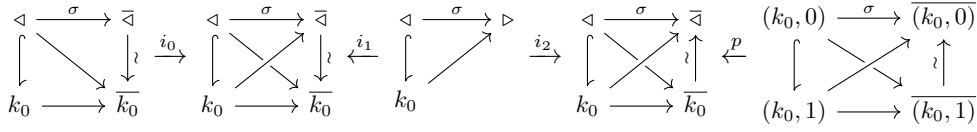


\begin{lm}
    \label{lm:condition_F_n_preserved_by_i_n}
    Let $n\in \{1,2,3\}$ and let $K$ be  a simplicial set.  
    Consider the corresponding morphism $i_n\colon F_a^+(K)\to F_b^+(K)$, where the indices $a$ and $b$ are determined by $n$. 
    Then a functor $\alpha \colon F_a^+(K)\to \zO^{\o,\natural}$ satisfies condition $(\star)_{a,\sigma}$ if and only if $i_n(K)^*(\alpha)$ satisfies condition \hyperlink{text:condition_F_0}{$(\star)_{b,\sigma}$}.
\end{lm}
\begin{proof}
    The lemma follows from inspection of the different conditions.
\end{proof}

We now establish that the definition of $(\zC_n)_\bullet$ as $\Hom^{\sigma}(F_n^+(\Delta^\bullet),\zO^{\o,\natural})$ yields an $\infty$-category.

\begin{lm}
	\label{lm:C_n_is_quasicategory}
	For $n \in \{1,2,3\}$, the simplicial set $\zC_n$ is an $\infty$-category.
\end{lm}

\begin{proof}
	We first prove the result for the simplicial set $\zC_1$. Let $k,m\in \bN$ with $0 < k < m$ and $f \in \Hom^\sigma(F_1^+(\Lambda_k^m),\zO^{\o,\natural})$. We will show that the existence of a lift in the following diagram
	\[
	\begin{tikzcd}  
		F_1^+(\Lambda_k^m) \ar[r,"f"] \ar[d,"" left]	&	\zO^{\o,\natural} \\
		F_1^+(\Delta^m) \ar[ru,dashed]					& 	
	\end{tikzcd}
	\]
	by decomposing the vertical map as a sequence of inner anodyne morphisms. First, let $X_0$ denote the simplicial set $\Delta^{\triangleleft 0 \dots n} \cup F_1^+(\Lambda_k^m)\cup \Delta^{\bar{\triangleleft}\bar{0}\dots\bar{m}}$. It is clear that the inclusion $F_1^+(\Lambda_k^m)\to X_0$ is inner anodyne. Since $\emptyset\to \Delta^m$ is both right and left anodyne, we may choose an increasing filtration 
	\begin{equation}
		\label{eqn:filtration_Lambda_j}
	\emptyset = \Lambda_0 \subset \Lambda_1 \subset \dots \subset \Lambda_r = \Delta^m
	\end{equation}
	of subsimplicial sets of $\Delta^m$ such that each inclusion $\Lambda_j\to \Lambda_{j+1}$ is the pushout of a horn inclusion $\Lambda_{k_j}^{m_j}\to \Delta^{m_j}$, which either is inner anodyne or satisfies $m_j \leq 1$. Note that in the latter case, the horn inclusion is left or right anodyne. 
	Introduce the simplicial sets
	\begin{eqnarray*}
		X_j & = & X_0 \cup (\Lambda_j^\triangleleft\ast \bar{\Lambda_j^\triangleleft})\\
		Y_j & = & X_0\cup (\Lambda_{j+1}^\triangleleft \ast \bar{\Lambda_j^\triangleleft})
		\cup 
		(\Lambda_j^\triangleleft \ast \bar{\Lambda_{j+1}^\triangleleft}),
	\end{eqnarray*}
	so that the inclusion $X_0 \to F_1^+(\Delta^m)$ can be written as the sequence of inclusions
	\[
	X_0 \,\subseteq Y_0 \,\subseteq 
	X_1 \,\subseteq Y_1 \,\subseteq 
	\dots \,\subseteq
	X_{r-1} \,\subseteq Y_{r-1} \,\subseteq 
	X_r = F_1^+(\Delta^m).
	\]
	Each inclusion $Y_j\to X_{j+1}$ is the pushout of the morphism 
	$(\Lambda_{j}^\tri\subseteq \Lambda_{j+1}^\tri)\boxast 
	(\bar{\Lambda_{j}^\tri} \subseteq \bar{\Lambda_{j+1}^\tri})$; lemma \ref{lm:calcul_box_join} implies that it is inner anodyne. For $p \in [m]$, let $\Lambda_j(p)$ denote the intersection of $\Lambda_j$ with the face of $\Delta^m$ opposed to vertex $p$. For every $j$ with $0\leq j < r$, let $p_j$ denote the index of the unique face of $\Delta^m$ that contains $\Lambda_{j+1}\setminus \Lambda_{j}$. Now consider the inclusion $X_j\to Y_j$. By construction of the filtration \eqref{eqn:filtration_Lambda_j}, this inclusion is obtained as a pushout of the map
	\[
	\begin{tikzcd}
	\left(\Lambda^\tri_{j}\ast \bar{\Lambda^\tri_{j}(p_j)}\cup \Lambda^\tri_{j}(p_j)\ast\bar{\Lambda^\tri_{j+1}(p_j)}\right)
	\cup 
	\left(\Lambda^\tri_{j}(p_j)\ast \bar{\Lambda^\tri_{j}}\cup \Lambda^\tri_{j+1}(p_j)\ast\bar{\Lambda^\tri_{j}(p_j)}\right)
	\ar[d]\\
	\left(\Lambda^\tri_{j}\ast\bar{\Lambda^\tri_{j+1}(p_j)}\right)\cup 
	\left(\Lambda^\tri_{j+1}(p_j)\ast\bar{\Lambda^\tri_{j}}\right).
	\end{tikzcd}
	\]
	This last map is the pushout of the map $(\Lambda^\tri_{j}(p_j)\subseteq \Lambda^\tri_j)\boxast (\Lambda^\tri_j(p_j)\subseteq \Lambda^\tri_{j+1}(p_j))$ with its symetric; it is therefore an inner anodyne map by lemma \ref{lm:calcul_box_join}, as desired. Therefore $f$ extends to a map $\tilde{f}\colon F_1^+(\Delta^m)\to \zO^{\o,\natural}$. The fact that $\tilde{f}$ belongs to the subset $\Hom^\sigma(F_1^+(\Delta^m),\zO^{\o,\natural})$ follows directly from the similar hypothesis for $f$; this concludes the proof that $\zC_1$ is an $\infty$-category.

	One can give a very similar proof for $\zC_3$, simply reversing the direction of the edges in $\bar{\Delta^m}$; we shall therefore omit the details.

	We now turn to the case of $\zC_2$. As before, let $k,m\in \bN$ with $0 < k < m$. It is enough to prove that the inclusion $F_2^+(\Lambda_k^m)\to F_2^+(\Delta^m)$ is inner anodyne. One simply observes that this inclusion can be described as successively adding to $F_2^+(\Lambda_k^m) = (\Lambda_k^m)^{\tri\triangleright}$ fillers of the inner horns $\Lambda_k^{0\dots m}$, $\Lambda_{k}^{\tri 0\dots m}$, $\Lambda_{k}^{0\dots m\triangleright}$ and $\Lambda_k^{\tri 0\dots m\triangleright}$, which proves the claim.
\end{proof}


\subsection{Proof of theorem \ref{thm:comparaison_bo_ext}} 
\label{sec:proof_of_theorem_}

\subsubsection{Strategy of proof}

We explain our approach to proving that 
the functors $i_0^*$, $i_1^*$ and $i_2^*$
\begin{equation*}
	\begin{tikzcd}
		\Ext(\sigma) &
		\zC_1 \ar[l,"i_0^*" above]\ar[r,"i_1^*"] &
		\zC_2 &
        \zC_3 \ar[l, "i_2^*" above]
	\end{tikzcd}
\end{equation*}
from zigzag \eqref{eqn:zigzag_comparaison_categories} are all equivalences of $\infty$-categories. The remaining case of the functor $p^*\colon \zC_3\to \bo_\sigma$ is treated separately, with different arguments, in section \ref{sub:study_map_p}.

Fix $n \in \{1,2,3\}$ and consider the associated natural transformation $i_n\colon F_a^+\to F_{b}^+$, where the indices $a,b\in \{0,1,2, 3\}$ are determined by $n$.

\todo{justifier que $i_n(K)$ pour différents $K$ induit un morphisme sur les $\pi_0\Map$.}

\begin{nota}
    \label{nota:cylinder_J}
    Let $\cJ$ denote the nerve of free-living isomorphism (i.e. the contractible groupoid on $2$ objects), which is an interval object  for Joyal's model structure (see the book by Cisinski \cite{cisinski-book} and also \cite[Appendix A]{campbell_joyals_cylinder_conjecture_2021}). 
    Given a simplicial set $K$, 
    we introduce the pushout
    \[
    F_{a,b,\cJ}^+(K) = F_{a}^+(K\times \cJ) \mathop{\amalg}\limits_{F_{a}^+(K)^{\amalg 2}} F_{b}^+(K)^{\amalg 2}.
    \]
    Note that $i_n(\bar{K})$ factors as a composite
    \begin{equation}
        \label{eqn:factorization_i_n_as_j_n_k_n}
        F_{a}^+(K\times \cJ) \stackrel{j_n}{\too} F_{a,b,\cJ}^+(K) \stackrel{k_n}{\too} F_{b}^+(K\times \cJ).
    \end{equation} 
\end{nota}

The following result is central in our approach.


\begin{lm}
    \label{lm:i_n_induce_bijections_on_pi_0}
    Let $n\in \{1,2,3\}$. Suppose that for all simplicial sets $K$, the morphism $i_n(K) \col F_a^+(K)\to F_{b}^+(K)$ is  marked anodyne 
    and the induced morphism $k_n\colon F_{a,b,\cJ}^+({K})\to F_b^+({K}\times \cJ)$ from factorization \eqref{eqn:factorization_i_n_as_j_n_k_n} is a monomorphism which is bijective on $0$-simplices.
    Then the induced map 
    \begin{equation*}
        i_n(K)^* \col \pi_0\Map(K,\zC_{b}) \to \pi_0\Map(K,\zC_a)
    \end{equation*}
    is a bijection for all $K$.
\end{lm}
\begin{proof}
    Using the left lifting property of marked anodyne morphisms against morphisms of the form $X^{\natural}\to *$ for $X$ an $\infty$-category (lemma \ref{lm:left_lifting_property_marked_anodyne}),
    we obtain that the map 
    $$\Hom_{\sSet^+}(F_{b}^+(K),\zO^{\o,\natural}) \to \Hom_{\sSet^+}(F_{a}^+(K),\zO^{\o,\natural})$$ is surjective. Hence so is its quotient $$i_n(K)^*\colon \pi_0\Map(K,\zC_{b}) \to \pi_0\Map(K,\zC_{a}).$$
    \todo{[OUBLI ?] why does it pass to the quotient?}
	We now prove that $i_n(K)^*$ is also injective. Let $\alpha, \alpha' \col F_{b}^+(K) \to \zO^{\o,\natural}$ be such that $i_n^*(\alpha) \simeq i_n^*(\alpha')$ in $\Map(K,\zC_n)$. Using the standard categorical cylinder $K\amalg K\to {K}\times\cJ \to K$ of $K$, the latter condition means that we can fill the following diagram of solid lines
	\begin{equation}
		\label{eqn:diag_extension_cylinder}
		\begin{tikzcd}[row sep = small]
			F_a^+(K) \ar[rd]\ar[rrd, bend left =25, "{i_n^*(\alpha)}"] & & \\
			& F_a^+({K}\times \cJ) \ar[r,dotted, "{\bar{\alpha}}"] & \zO^{\o,\natural}.\\
			F_a^+(K) \ar[ru]\ar[rru, bend right=25, "{i_n^*(\alpha')}" below] & & \\
		\end{tikzcd}
	\end{equation}
	so that $\bar{\alpha}$ respects conditions $(*)_{a,\sigma}$. To show that $\alpha\simeq \alpha'$ in $\Map(K,\zC_{b})$, we have to prove that the corresponding diagram for $F_{b}^+$, namely
	\begin{equation}
		\label{eqn:diag_extension_cylinder_F_n+1}
		\begin{tikzcd}[row sep = small]
			F_{b}^+(K) \ar[rd]\ar[rrd, bend left =25, "{\alpha}"] & & \\
			& F_{b}^+({K}\times\cJ) \ar[r,dotted, "{\tilde{\alpha}}"] & \zO^{\o,\natural},\\
			F_{b}^+(K) \ar[ru]\ar[rru, bend right=25, "{\alpha'}" below] & & \\
		\end{tikzcd}
	\end{equation}
	can be filled by a map $\tilde{\alpha}\col F_{b}^+(K)\to \zO^{\o,\natural}$ that satisfies conditions $(*)_{b,\sigma}$. Now observe that $i_n({K}\times\cJ)$ factors through the pushout
	$F_{a,b,\cJ}^+({K})$.
    By hypothesis, $i_n(K)$ is marked anodyne, whence so is $j_n$. Using the right simplification property of marked anodyne morphisms (Proposition \ref{prop:right_simplification_marked_anodyne}), we deduce that there exists a lift $\tilde{\alpha}$ in the following diagram
	\begin{equation}
		\label{eqn:lm_extending_along_marked_anodyne}
		\begin{tikzcd}
		F_{a}^+({K}\times\cJ) \arrow[dd, "i_n"', bend right=60] \arrow[d, "j_n"] \arrow[rrd, "\bar{\alpha}", bend left] & & \\
		{F_{a,b,\cJ}^+(K)} \arrow[d, "k_n"] \arrow[rr, "{(\alpha\amalg\alpha',\bar{\alpha})}"]&  & \zO^{\o,\natural}. \\
		F_{b}^+({K}\times\cJ) \arrow[rru, "\widetilde{\alpha}", dotted, bend right]&  & 
		\end{tikzcd}
	\end{equation}
    By lemma \ref{lm:condition_F_n_preserved_by_i_n}, 
    any lift $\widetilde{\alpha}$ as above will automatically satisfies conditions $(*)_{b,\sigma}$. This shows that $\alpha$ and $\alpha'$ are equivalent as points in $\Map(K,\zC_{b})$, as desired.

\end{proof}

We will prove that the maps $i_0$, $i_1$ and $i_2$ are marked anodyne and use the above lemma to deduce that the $\infty$-categories $\Ext(\sigma)$, $\zC_1$, $\zC_2$ and $\zC_3$ are equivalent. The remaining step will be to deal with the map $p$, which is a marked equivalence but not a marked anodyne morphism, so that we cannot use lemma \ref{lm:i_n_induce_bijections_on_pi_0}; we will therefore rely on a different argument, involving a careful analysis of categorical cylinders in $\zC_3$ and $\bo_\sigma$.



\subsubsection{Study of the map \texorpdfstring{$i_0$}{i0}} 

\begin{lm}
    \label{lm:i_0_mono_bij_on_objects}
    For every simplicial set $K$, the morphism $k_0(K)\colon F_{0,1,\cJ}^+(K)\to F_1^+(K)$ is a monomorphism which is bijective on $0$-simplices.
\end{lm}
\begin{proof}
    Let $K\in \sSet$. Observe that $i_0(K)\colon F_0^+(K)\to F_1^+(K)$ is a bijection on $0$-simplices. Since $K$ is arbitrary, the map $i_0(K\times \cJ)$ also has this property. 
    Moreover, so does the  map $j_0(K)\colon F_0^+(K\times \cJ)\to F_{0,1,\cJ}^+(K)$, as it is obtained as a cobase change of $i_0(K)^{\amalg 2}$.
    By 2-out-of-3 property, we deduce that $k_0(K)$ also has this property.

    Let us show that $k_0(K)$ is a monomorphism. Suppose that $t$ and $t'$ are two $m$-simplices of $F_{0,1,\cJ}^+(K)$ with the same image under $k_0(K)$. 
    We separate three cases:
    \begin{enumerate}
        \item either both $t$ and $t'$ lift to simplices of $F_0^+({K}\times\cJ)$, or 
        \item both lift to $F_1^+(K)^{\amalg 2}$, or 
        \item $t$ lifts to a simplex in $F_0^+({K}\times\cJ)$ and $t'$ to one in $F_1^+(K)^{\amalg 2}$.
    \end{enumerate}
    For the first case, note that $j_0({K}\times\cJ)$ is defined as a pushout of a monomorphism, hence is a monormophism, so that $t=t'$ as desired.
    For the second case, writing $\widetilde{t}$ and $\widetilde{t'}$ for choices of lifts of $t$ and $t'$ in $F_1^+(K)^{\amalg 2}$, we easily see that either $\widetilde{t}$ and $\widetilde{t'}$ belong to the same of the two copies of $F_1^+(K)$, or they factor through $(\Delta^{\triangleleft \bar{\triangleleft}})^{\amalg 2}$. 
    Finally, in the third case, the simplices $t$ and $t'$ have to factor through one of the two copies of $F_0^+(K)$; since $F_0^+(K)\to F_1^+({K}\times\cJ)$ is a monomorphism, the simplices $t$ and $t'$ must coincide.
\end{proof}

\begin{prop}
	\label{prop:i_0_marked_anodyne}
	For every simplicial set $K$, the map $i_0 \col F_0^+(K) \to F_1^+(K)$ is marked anodyne.
\end{prop}

\begin{proof}
	Recall the marking on these simplicial sets: on $F_0^+(K) = K^\triangleleft\times \Delta^1$, the marked edges are those of the form $\bar{\triangleleft} \to \bar{y}$, with $y\in K$, whereas for $F_1^+(K) = \mathop{\colim}\limits_{\Delta^m\to K^\triangleleft}(\Delta^m \ast \Delta^m)$, all edges $\bar{x} \to \bar{y}$ with $x\to y$ in $K^\triangleleft$ are marked. The morphism $i_0$ factors through the marked simplicial set $\tilde{F}_0^+(K)$ obtained from $F_0^+(K)$ by further marking the edges $(x,1) \to (y,1)$, for $x\to y$ in $K^\triangleleft$. The resulting inclusion $F_0^+(K) \to \tilde{F}_0^+(K)$ is easily seen to be marked anodyne. 

	Let $\tilde{i}_0 \col \tilde{F}_0^+(K) \to F_1^+(K)$ denote the map induced by $i_0$. We will prove that it is marked anodyne as well. By construction of $i_0$ as a colimit, it is enough to prove the claim for $K = \Delta^{m-1}$, in which case the morphism $\tilde{i}_0$ is the canonical inclusion $\Delta^m\times\Delta^1 \to \Delta^m\ast \Delta^m$, where the non-trivial marking occurs in the second copy of $\Delta^m$. Therefore, we are left with showing the following combinatorial result.
\end{proof}

\begin{lm}
	\label{lm:delta^m_combinatorial}
	For every $m \in \bN$, consider the canonical inclusion $\tilde{i}_0 \col \Delta^m\times\Delta^1 \to \Delta^m\ast \Delta^m$, where both simplicial sets are endowed with the minimal marking that makes $(\Delta^m)^\sharp \times \{1\}$ a marked simplicial subset of them. Then $\tilde{i}_0$ is marked anodyne.
\end{lm}

\begin{proof}
	Note that the inclusion $\tilde{i}_0$ factors through the simplicial set 
	\begin{equation}
		\label{eqn:def_A_m}
	A^m = (\Delta^m\times \Delta^1) \mathop{\cup}\limits_{\Delta^{\bar{0}\bar{1}\dots\bar{m}} \mathop{\cup}\limits_{\Delta^{\bar{m}}} \Delta^{m\bar{m}}} \Delta^{m\bar{0}\bar{1}\dots\bar{m}}.
	\end{equation}
	First, we want to show that the inclusion $A^m\to \Delta^m\ast \Delta^m$ is marked anodyne. To this purpose, we consider the inclusion $s_m\col S_m \to \Delta^m\ast\Delta^m$, where $S_m$ denotes the spine $\mathrm{Sp}^{0\dots m \bar{0}\dots \bar{m}}$, endowed with the maximal marking that turns $s_m$ into a morphism of marked simplicial sets. As $s_m$ is clearly marked anodyne and factors through $A^m$, it will suffice to show that the inclusion $S_m \to A^m$ is marked anodyne.
	The underlying simplicial set of $A^m$ is a union of $m+2$ simplices of dimension $(m+1)$, denoted $\tau_0$, $\dots$, $\tau_m$, $\bar{\tau}$ and defined as follows:
	\begin{itemize}
		\item for $0 \leq k \leq m$, the simplex $\tau_k$ is defined as $\Delta^{01\dots k\, \bar{k}\, \bar{(k+1)} \dots \bar{m}}$,
		\item the simplex $\bar{\tau}$ is $\Delta^{m\bar{0}\,\bar{1}\dots\bar{m}}$. 
	\end{itemize}
	In each case, the marking is induced by that of $A^m$.
	Writing $T_k$ for $S_m \cup \bar{\tau}\cup \tau_m \cup \dots \cup \tau_k$, we get a filtration of $A^m$ of the form
	\begin{equation}
		\label{eqn:filtration_A^m}
		S_m \ \subset \ 
		S_m \cup \bar{\tau} \ \subset \ 
		T_m \ \subset \ 
		T_{m-1} \ \subset \ 
		\dots \ \subset\ 
		T_0 = A^m.
	\end{equation}
	We will prove that at each step, the inclusion is a marked anodyne morphism. 
	This is clear for $S_m \subset S_m\cup\bar{\tau}$. The second inclusion $S_m\cup\bar{\tau} \subset T_m$ is obtained as the pushout
	\begin{equation}
		\label{eqn:pushout_S_m_T_m}
		\begin{tikzcd}
		(\mathrm{Sp}^{0\dots m\bar{m}})^{\flat} \arrow[r] \arrow[d] \arrow[rd, "\ulcorner" near end, phantom] & S_m\cup \bar{\tau} \arrow[d] \\
		\tau_m \arrow[r]                                                                 & T_m
		\end{tikzcd}
	\end{equation}
	and is therefore marked anodyne, since so is the left vertical map. We now prove that the inclusion $T_k \subset T_{k-1}$ is marked anodyne for every $0 < k \leq m$. 
	We introduce the marked simplicial set $\tau_k\langle k \rangle$ defined as the face opposed to vertex ${k}$ in $\tau_k$, or more explicitly as $\Delta^{0\dots (k-1) \bar{k} \dots \bar{m}}$, endowed with the induced marking. The intersection $\tau'_k = T_k \cap \tau_{k-1}$ can then be expressed as
	\begin{equation}
		\label{eqn:definition_tau'_k}
		\tau'_k = \tau_k\langle k \rangle \mathop{\cup}\limits_{(\Delta^{\bar{k}\dots\bar{m}})^\sharp} (\Delta^{\bar{(k-1)}\, \bar{k} \dots \bar{m}})^\sharp
	\end{equation}
	and the inclusion $T_k \subset T_{k-1}$ as the pushout
	\begin{equation}
		\label{eqn:pushout_T_k_T_k-1}
		\begin{tikzcd}
		\tau'_k	\arrow[r] \arrow[d] \arrow[rd, "\ulcorner" near end, phantom] & T_k \arrow[d] \\
		\tau_{k-1} \arrow[r]                                                                 & T_{k-1}.
		\end{tikzcd}
	\end{equation}
	To see that $T_k\to T_{k-1}$ is marked anodyne, it therefore suffices to show that $\tau_k'\to \tau_{k-1}$ has this property. Looking at equation \eqref{eqn:definition_tau'_k}, we observe that the latter inclusion is of the form $(I_0\subset I) \boxast (J_0\subset J)$, with $I_0 = \emptyset$, $I =[k-1]$, $J_0 = \{\bar{k},\dots,\bar{m}\}$ and $J = J_0\cup\{\bar{(k-1)}\}$. By lemma \ref{lm:combinatorial_inclusion}, we deduce that this map is marked anodyne, as desired.

	It remains to prove that the inclusion $\Delta^m\times\Delta^1 \to A^m$ is marked anodyne. By equation \eqref{eqn:def_A_m}, it suffices to prove that the inclusion $\Delta^{\bar{0}\dots\bar{m}}\cup \Delta^{mF\bar{m}}\inc \Delta^{m\bar{0}\dots\bar{m}}$ is marked anodyne. Denoting $B_i = \Delta^{\bar{0}\dots\bar{m}}\cup \Delta^{m\bar{i}\dots\bar{m}}$, the latter map can be written as the composite
	\begin{equation*}
		B_m 
		\:\subset\: B_{m-1}\:\subset\: \dots \:\subset\: 
		B_0.
	\end{equation*}
	Each inclusion $B_{i+1}\subset B_i$ is induced by its restriction $B_{i+1}\cap \Delta^{m\bar{i}\dots\bar{m}} \subset B_i\cap \Delta^{m\bar{i}\dots\bar{m}}$, which is marked anodyne by lemma \ref{lm:calcul_box_join} applied with $I_0 = \emptyset$, $I = \{m\}$, $J = \{\bar{i},\dots,\bar{m}\}$ and $J_0 = J\setminus \{\bar{i}\}$; therefore so is the map  $\Delta^m\times\Delta^1 \to A^m$. This completes the proof of lemma \ref{lm:delta^m_combinatorial}.
\end{proof}

\subsubsection{Study of the maps \texorpdfstring{$i_1$}{i1} and \texorpdfstring{$i_2$}{i2}}

\begin{lm}
    \label{lm:i_1_i_2_mono_bij_on_objects}
    For $K$ a simplicial set, the morphisms $k_1(K)\colon F_2^+(K)\to F_{2,1,\cJ}^+(K)$ and $k_2(K)\colon F_2^+(K)\to F_{2,3,\cJ}^+(K)$ are monomorphisms which are bijective on $0$-simplices.
\end{lm}
\begin{proof}
    The argument is completely analogous to that of the proof of lemma \ref{lm:i_0_mono_bij_on_objects}.
\end{proof}

\begin{lm}
\label{lm:i_1_inner_anodyne}
For all simplicial set $K$, the maps $i_1\colon F_2^+(K)\to F_1^+(K)$ and $i_2\colon F_3^+(K)\to F_2^+(K)$ are marked anodyne.
\end{lm}

\begin{proof}
	We only give the proof for $i_1$, the case of the map $i_2$ being very similar.
	Recall that $i_1$ is defined by taking the colimit over all simplices $\Delta^m\to K^{\triangleleft}$ of the morphisms $i_{1,m}\colon (\Delta^m)^\triangleright \to F_1^+(K)$, sending $\Delta^m$ to the first copy of itself in $(\Delta^m)^{\ast 2} \subseteq F_1^+(K)$ and $\triangleright$ to $\bar{\triangleleft}$. 
	Note that we may restrict the colimit to those simplices $\Delta^m$ that contains the cone point $\triangleleft$, in which case the map $i_{1,m}$ factors through $(\Delta^m)^{\ast 2}$. Let $\Delta^m\to K^{\triangleleft}$ be such a simplex; it now suffices to show that the factored map $\bar{i}_{1,m} \colon (\Delta^m)^\triangleright\to (\Delta^m)^{\ast 2}$ is marked anodyne.

	For the rest of this proof, we relabel $\triangleleft$ as $0$, so that we can identify $\bar{i}_{1,m}$ with the obvious inclusion of $(\Delta^m)^\triangleright\cong \Delta^{0 \dots m\bar{0}}$ into $(\Delta^m)^{\ast 2} \cong \Delta^{0 \dots m\bar{0}\dots \bar{m}}$. Now this map is the composite of the sequence
	\[
	\Delta^{0\dots m\bar{0}} \subset 
	\Delta^{0\dots m\bar{0}\bar{1}} \subset 
	\dots 
	\Delta^{0\dots m\bar{0}\bar{1}\dots \bar{j}} \subset 
	\dots 
	\Delta^{0\dots m\bar{0}\bar{1}\dots \bar{m}}.
	\]
	Since all edges of the form $\bar{\ell}\to\bar{p}$ are marked, by lemma \ref{lm:combinatorial_inclusion} we deduce that each of these inclusions is marked anodyne, whence the result.
\end{proof}




\subsubsection{Study of the map \texorpdfstring{$p$}{p}}
\label{sub:study_map_p}

This subsection is devoted to the last step of the comparison described in zigzag \eqref{eqn:zigzag_comparaison_categories}, namely we prove the following result.

\begin{prop}
	\label{prop:p_is_equivalence}
    The functor $p^*\colon \zC_3\to\bo_\sigma$  
	induced by the natural transformation $p\colon G^+\to F_3^+$  is an equivalence.
\end{prop}

We start with some preliminary results.

\begin{lm}
\label{lm:G_preserves_monos}
The functors $F_3, G\colon \sSet \to \sSet$ both preserve monomorphisms.
\end{lm}

\begin{proof}
	Since $s_*$ and $(-)^\triangleleft$ preserves monomorphisms, so does their composite $F_3$. We now turn to the case of $G$.
	Let $\varphi \colon A\to B$ be a monomorphism of simplicial sets and $n\in \bN$. Let $a_0$ and $a_1$ be two $n$-simplices of $G^+(A)$ whose image under $G^+(\varphi)$ coincide. We wish to prove that $a_0 = a_1$. Since $s_*(A\times \Delta^1)\to G^+(A)$ is an epimorphism, we may choose lifts $\widetilde{a_0}$ and $\widetilde{a_1}$ in $(s_*(A\times\Delta^1))_n$ of the two simplices. We distinguish several cases in the argument.
	\begin{enumerate}
		\item If both $\widetilde{a_0}$ and $\widetilde{a_1}$ belong to the subset $(s_*A)_n$, then their common image $\varphi(a_0) = \varphi(a_1)$ actually belongs to $(s_*\{0\})_n$, which is a subset of $G^+(A)_n$, so that $a_0 = a_1$.
		\item If none of $\widetilde{a_0}$ and $\widetilde{a_1}$ belongs to the subset $(s_*A)_n$, since $s_*(A\times\Delta^1)_n\setminus s_*(A)_n \to G^+(B)_n$ is an inclusion, we deduce that $\widetilde{a_0} = \widetilde{a_1}$, so that again $a_0 = a_1$.
		\item Finally, the case where exactly one of  $\widetilde{a_0}$ and $\widetilde{a_1}$ belong to $s_*(A)_n$ is contradictory: indeed, by construction of $G^+(A)$ the subset $s_*(A)_n$ and its complement in $s_*(A\times \Delta^1)_n$ remain disjoint in the quotient $G^+(A)_n$, hence also in $G^+(B)_n$.
	\end{enumerate}
\end{proof}

\begin{lm}
\label{lm:p(K)_equivalence}
For every simplicial set $K$, the morphism $p(K)\colon G^+(K)\to F_3^+(K)$ is an equivalence of marked simplicial sets.
\end{lm}

\begin{proof}
    Let $\mathcal{U}$ be the class of simplicial sets $X$ for which $p(X)$ is an equivalence. We will show that $\mathcal{U} = \sSet$ by proving that $\mathcal{U}$ contains every representable $\Delta^m$ and is stable under isomorphisms, small coproducts, pushouts along monomorphisms and sequential colimits along monomorphisms. In other words, we will show that $\mathcal{U}$ is \emph{saturated by monomorphisms} in the sense of \cite[Definition 1.3.9]{cisinski-book} that contains all representables.

    We first show that all standard simplices are in $\mathcal{U}$.
	Let $m \in \bN$ and consider $p(\Delta^m) \colon G^+(\Delta^m)\to F_3^+(\Delta^m)$. This morphism admits a section $e$, \emph{non-naturally} in the variable $[m]\in \Delta$, whose underlying morphism of simplicial sets is defined as the composition
	\[
	e\colon 
	F_3(\Delta^m) = s_*\Delta^{\triangleleft 0 \dots m} \cong s_*\Delta^{00,01,\dots,0m} \subset s_*(\Delta^m\times\Delta^1) \to G(\Delta^m).
	\]
	where the middle identification of $(m+1)$-dimensional simplices sends $\triangleleft$ to $00$ and the vertex $i$ to $0i$ (following the notations of sections \ref{sub:diag_bo_sigma} and \ref{sub:diag_intermediate_steps}). Let $D_i$ be the non-degenerate $(m+1)$-simplex of $\Delta^m\times \Delta^1$ containing the edge $i0\to i1$, where we declare all edges contained in $\Delta^m\times \{0\}$ to be marked. Let $D_{i,i+1}$ denote the face of $D_i$ opposed to vertex $i1$. 
	The morphism $e$ then factors through a filtration
	\[
	F_3^+(\Delta^m) \stackrel{e_0}{\too} G_0^+\stackrel{e_1}{\too} G_1^+ \stackrel{}{\too} \dots \stackrel{e_m}{\too} G_m^+ = G^+(\Delta^m),
	\]
	where $G_i^+$ is the image of $s_*(\bigcup_{j = 0}^i D_j)$ by the quotient map $s_*(\Delta^m\times \Delta^1) \to G^+(\Delta^m)$, with the induced marking. Since the edge $i0\to (i+1)0$ is marked, the inclusion $s_*(D_{i,i+1})\to s_*(D_{i+1})$ is marked anodyne, hence so is its pushout
	\[
	s_*\bigcup_{j = 0}^i D_j \stackrel{\sim}{\too}
	s_*\bigcup_{j = 0}^{i+1} D_j.
	\]
	We thus obtain a diagram
	\[
	\begin{tikzcd}
	s_*\{0\} \arrow[d, equal] & {s_*\Delta^{00,01,\dots,i0}} \arrow[d, "\wr"] \arrow[l] \arrow[r, tail] & s_*\bigcup_{j = 0}^i D_j \arrow[d, "\wr"] \\
	s_*\{0\}                                & {s_*\Delta^{00,01,\dots,(i+1)0}} \arrow[r, tail] \arrow[l]              & s_*\bigcup_{j = 0}^{i+1} D_j             
	\end{tikzcd}
	\]
	where each row defines a cofibrant diagram in the projective model structure on the category of diagrams $\Fun(\bullet\leftarrow \bullet \rightarrow \bullet, \sSet^+)$. Taking pushouts yields a weak equivalence
	$e_i \colon G_i \stackrel{\sim}{\too} G_{i+1}$. This proves that $e$ is a marked equivalence, hence $\Delta^m\in \mathcal{U}$.

	It is clear that $\mathcal{U}$ is stable by isomorphisms. Consider a simplicial set of the form $K = \coprod_{i\in J} K_i$ with all the $K_i$ in $\mathcal{U}$. As $(-)^{\triangleleft}\colon \sSet \to \sSet_*$ and $s_*$ both preserve colimits, we can compute
	\begin{eqnarray*}
	F_3^+(K) & = 	& s_*\left( \left(\coprod_{i\in J}K_i\right)^\triangleleft\right)\\
		& \cong & s_*\left( \left(\coprod_{i\in J} K_i^\triangleleft\right)/{\coprod_{i\in J} \triangleleft_i}\right)\\
		& \cong & \mathop{\colim} \left(
		\coprod_{i\in J} s_*(K_i^\triangleleft) \lla
		\coprod_{i\in J} s_*\{\triangleleft_i\} \too
		s_*\{*\}
		\right) 
	\end{eqnarray*}
	and
	\begin{eqnarray*}
	G^+(K) & \cong & \mathop{\colim} \left(
		\coprod_{i\in J} s_*(K_i\times \Delta^1) \lla
		\coprod_{i\in J} s_*(K_i\times\{0\}) \too
		s_*\{0\}
		\right)\\
		& \cong & \mathop{\colim} \left(
		\coprod_{i\in J} G^+(K_i) \lla
		\coprod_{i\in J} s_*\{0\} \too
		s_*\{0\}
		\right).
	\end{eqnarray*}
		Since each map $p(K_i) \colon G^+(K_i)\to F_3^+(K_i) = s_*(K_i^\triangleleft)$ is a marked equivalence and the diagrams are cofibrant, we deduce that $p(K)$ is also an equivalence. This proves that $\mathcal{U}$ is stable under small colimits.

		Suppose now that $K$ is a pushout $B\coprod_A C$, with $A$, $B$ and $C$ in $\mathcal{U}$ and $A\to B$ a monomorphism. We will prove that $K \in \mathcal{U}$. As before, rewriting the colimits gives isomorphisms
		\begin{equation*}
			F_3^+(K) \cong F_3^+(B) \coprod_{F_3^+(A)} F_3^+(C) 
			\qquad \text{and} \qquad
			G^+(K) \cong G^+(B) \coprod_{G^+(A)} G^+(C).
		\end{equation*}
		To prove that $p(K)$ is a weak equivalence, it suffices to show that the pushout diagrams 
		$F_3^+(B) \leftarrow {F_3^+(A)} \to F_3^+(C) $ and
		$G^+(B) \leftarrow {G^+(A)} \to G^+(C) $ are cofibrant. This is a consequence of the fact that $F_3^+$ and $G^+$ both preserve monomorphisms; in the first case, this is clear whereas in the latter, it is given by lemma \ref{lm:G_preserves_monos}. Finally, the proof that $\mathcal{U}$ is stable by sequential colimits along monomorphisms comes from a similar argument, therefore we omit it.
\end{proof}

Recall from \ref{nota:cylinder_J} the notation $\cJ$ for the standard interval object.
We will use the fact that Joyal's model structure on $\sSet$ can be obtained à la Cisinski using $\mathcal{J}\times (-)$ as an exact cylinder \cite{cisinski-book}.


\begin{lm}
\label{lm:F_and_G_preserve_cylinder}
For every simplicial set $K$, the image of the categorical equivalence $q\colon K\times \mathcal{J}\to K$ under the functors $F_3^+$ and $G^+$ is a marked equivalence.
\end{lm}

\begin{proof}
	Let $\mathcal{U}_F$ (respectively $\mathcal{U}_G$) be the class of simplicial sets $X$ such that $F_3^+(q)\colon F_3^+(X\times \mathcal{J})\to F_3^+(X)$ (resp. $G^+(q)\colon G^+(X\times \mathcal{J})\to G^+(X)$) is a marked equivalence. We aim at proving that $\mathcal{U}_F = \mathcal{U}_G = \sSet$. Using arguments similar to those of the proof of lemma \ref{lm:p(K)_equivalence}, one easily shows that both $\mathcal{U}_F$ and $\mathcal{U}_G$ are stable under isomorphisms, small coproducts, pushouts along a monomorphism and sequential colimits along monomorphisms. It is thus enough to prove that $\mathcal{U}_F$ and $\mathcal{U}_G$ contains every representable $\Delta^m$. The key point is the observation that both $F_3$ and $G$ restricts to an endofunctor on the full subcategory of $\sSet$ given by the essential image of the nerve functor from $1$-categories. More precisely, we have
	\[
	F_3(\Delta^m\times\mathcal{J}) \cong N \left(s_*([m]\times \sJ)^\triangleleft\right)
	\]
	and
	\[
	G(\Delta^m\times\mathcal{J}) \cong N \left(s_*([m]\times \sJ\times [1])\mathop{\coprod}_{s_*([m]\times\sJ\times\{0\})}s_*\{0\}\right),
	\]
	from which one readily verifies that $F_3^+(q)$ and $G^+(q)$ are the image under the nerve functor of an equivalence of $1$-categories, hence are marked equivalence.
\end{proof}

We now turn to the proof of the last step in our comparison of $\bo_\sigma$ and $\Ext(\sigma)$, namely the proof that $p^* \colon \zC_3\to \bo$ is an equivalence of $\infty$-categories.

\begin{proof}[Proof of proposition \ref{prop:p_is_equivalence}]
	Let $K$ be a simplicial set. 
    Recall that $\bo_\sigma$ is defined by identifying, naturally in $K$, the set $\Hom(K,\bo_\sigma)$ with the subset $\Hom^\sigma(G^+(K), \zO^{\o,\natural})$ of $\Hom_{\sSet^+}(G^+(K), \zO^{\o,\natural})$ of morphisms satisfying conditions \hyperlink{text:condition_bo_sigma}{$(\star)_{G,\sigma}$}. Similarly, $\Hom(K,\zC_3)$ is identified with the set $\Hom^\sigma(F_3^+(K),\zO^{\o,\natural})$. Since $p(K)$ preserves the conditions \hyperlink{text:condition_bo_sigma}{$(\star)_{G,\sigma}$}, the map 
	$$p(K)^*\colon \Hom(F_3^+(K),\zO^{\o,\natural})\to \Hom(G^+(K),\zO^{\o,\natural})$$
	restricts to a map $\Hom^\sigma(F_3^+(K),\zO^{\o,\natural})\to \Hom^\sigma(G^+(K),\zO^{\o,\natural})$. \relax

    We will need to consider two quotients of these $\Hom$ sets, whose associated equivalence relations we will denote $\sim$ and $\approx$. They correspond respectively to the homotopy relations in $\Fun(K,\zC_3)$ and $\Map^\flat(F_3^+(K),\zO^{\o,\natural})$.
    We only describe those two relations in the case of $\Hom^\sigma(F_3^+(K),\zO^{\o,\natural})$, the case of $\Hom^\sigma(G^+(K),\zO^{\o,\natural})$ being similar.
	\begin{itemize}
		\item \textit{(Definition of $\sim$).} First, we consider the set of connected components $\pi_0\Map(K,\zC_3)$, which is defined as the quotient of $\Hom(F_3^+(K),\zO^{\o,\natural})$ by the homotopy equivalence relation $\sim$ in the functor $\infty$-category $\Fun(K, \zC_3)$.
		\item \textit{(Definition of $\approx$).}  \sloppy Second, we consider the set of connected components $\pi_0\Map^\sharp(F_3^+(K),\zO^{\o,\natural})$, that is the quotient of $\Hom(F_3^+(K),\zO^{\o,\natural})$ by the homotopy equivalence relation $\approx$ in the functor $\infty$-category $\Map^\flat(F_3^+(K),\zO^{\o,\natural})$.
	\end{itemize}
	Using the characterization of equivalences in functor $\infty$-categories of \cite[\href{https://kerodon.net/tag/01KA}{Theorem 01KA}]{kerodon} (or more precisely, a slight generalization of this result to marked simplicial sets), we can rephrase the definitions of $\sim$ and $\approx$ more explicitly.

    Let $f_0$ and $f_1$ be two maps $F_3^+(K)\to \zO^{\o,\natural}$. 
	\begin{itemize}
		\item[\framebox{$\sim$}] Both of the following conditions are equivalent to asserting that $f_0\sim f_1$.
		\begin{enumerate}
			\item[$(i_\sim)$] There exists a factorization of the fold map $(\id_K,\id_K)$ as $K\coprod K \to \bar{K}\stackrel{\rho}{\to} K$, with $\rho$ a categorical equivalence, and a lift in the diagram
			\[
			\begin{tikzcd}  
				F_3^+(K)^{\amalg 2} \ar[r,"{(f_0,f_1)}"] \ar[d,"{(F_3^+(i_0),F_3^+(i_1))}" left]	&	\zO^{\o,\natural}\\
				F_3^+(\bar{K}) \ar[ru,dashed, bend right,"\bar{f}"]					& 	
			\end{tikzcd}
			\]
			such that $\bar{f}$ satisfies conditions \hyperlink{text:condition_bo_sigma}{$(\star)_{G,\sigma}$}.
			\item[$(ii_\sim)$] For every factorization of the fold map $(\id_K,\id_K)$ as $K\coprod K \stackrel{(s_0,s_1)}{\too} \bar{K}{\to} K$, where $s_0$ and $s_1$ are disjoint monomorphisms, there exists a lift in the diagram
			\[
			\begin{tikzcd}  
				F_3^+(K)^{\amalg 2} \ar[r,"{(f_0,f_1)}"] \ar[d,"{(F_3^+(s_0),F_3^+(s_1))}" left]	&	\zO^{\o,\natural}\\
				F_3^+(\bar{K}) \ar[ru,dashed, bend right,"\bar{f}"]					& 	
			\end{tikzcd}
			\]
			such that $\bar{f}$ satisfies conditions \hyperlink{text:condition_bo_sigma}{$(\star)_{\sigma}$}.
		\end{enumerate}
		\vspace{2mm}
		\item[\framebox{$\approx$}] Both of the following conditions are equivalent to asserting that $f_0\approx f_1$.
		\begin{enumerate}
			\item[$(i_\approx)$] There exists a factorization of the fold map $\left(\id_{F_3^+(K)},\id_{F_3^+(K)}\right)$ as 
                $$F_3^+(K)\coprod_{\Delta^1} F_3^+(K) \stackrel{\iota}{\too} \bar{F_3^+(K)}\stackrel{\rho}{\too} F_3^+(K),$$ with $\rho$ a cartesian equivalence (in the sense of \cite{htt}) 
            \todo{[plus ou moins fait] harmoniser la terminologie sur les "cartesian equiv"}, and a lift in the diagram
			\[
			\begin{tikzcd}  
				F_3^+(K)\coprod_{\Delta^1} F_3^+(K) \ar[r,"{(f_0,f_1)}"] \ar[d,"\iota" left]	&	\zO^{\o,\natural}.\\
				\bar{F_3^+({K})} \ar[ru,dashed, bend right,"\bar{f}"]					& 	
			\end{tikzcd}
			\]
			\item[$(ii_\approx)$] For every factorization of the fold map $\left(\id_{F_3^+(K)},\id_{F_3^+(K)}\right)$ as 
			$$F_3^+(K)\coprod_{\Delta^1} F_3^+(K) \stackrel{\iota}{\too} \bar{F_3^+(K)}\stackrel{\rho}{\too} F_3^+(K),$$
			where $\iota$ is a monomorphism, there exists a lift in the diagram
			\[
			\begin{tikzcd}  
				F_3^+(K)\coprod_{\Delta^1} F_3^+(K) \ar[r,"{(f_0,f_1)}"] \ar[d,"\iota" left]	&	\zO^{\o,\natural}.\\
				\bar{F_3^+({K})} \ar[ru,dashed, bend right,"\bar{f}"]					& 	
			\end{tikzcd}
			\]
		\end{enumerate}
	\end{itemize}
	Using these descriptions, one easily sees that $p(K)^*$ induces maps on the quotients by the equivalence relations $\sim$ and $\approx$, that we denote respectively
	\begin{eqnarray*}
	p^*_{\sigma,\sim} 	& \colon \Hom^\sigma(F_3^+(K), \zO^{\o,\natural})_{/\sim} \too \Hom^\sigma(G^+(K), \zO^{\o,\natural})_{/\sim},\\
	p^*_{\sigma,\approx}  & \colon \Hom^\sigma(F_3^+(K), \zO^{\o,\natural})_{/\approx} \too \Hom^\sigma(G^+(K), \zO^{\o,\natural})_{/\approx}.	
	\end{eqnarray*}

	By lemma \ref{lm:p(K)_equivalence},
	we know that $p(K)$ is marked weak equivalence, so that
	\[
	p^*_\approx \colon \Hom(F_3^+(K), \zO^{\o,\natural})_{/\approx} \too \Hom(G^+(K), \zO^{\o,\natural})_{/\approx}	
	\]
    is a bijection. It is an easy observation that a morphism $f\colon F_3^+(K)\to \zO^{\o,\natural}$ verifies the conditions \hyperlink{text:condition_F_n}{$(\star)_{3,\sigma}$}  
    if and only if $f\circ p(K)\colon G^+(K)\to \zO^{\o,\natural}$ satisfies the corresponding condition \hyperlink{text:condition_bo_sigma
    }{$(\star)_{G,\sigma}$}. Therefore the induced map $p^*_{\sigma, \approx}$ is a bijection. The fact that $p^*_{\sigma, \sim}$ is a bijection is now a consequence of the following result.\\

	\textbf{Claim.} The equivalence relations $\sim$ and $\approx$ coincide, both on $\Hom(F_3^+(K),\zO^{\o,\natural})$ and on $\Hom(G^+(K),\zO^{\o,\natural})$.\\

	We prove the claim for the functor $F_3^+$, the case of $G^+$ being similar. Consider two morphisms $f_0, f_1\colon F_3^+(K)\to \zO^{\o,\natural}$.
	\begin{enumerate}
		\item[\framebox{$\sim\,\Rightarrow\, \approx$}] Suppose $f_0\sim f_1$. By assumption $(ii_\sim)$, there exists a lift $\bar{f}$ in the diagram
		\[
		\begin{tikzcd}  
			F_3^+(K)^{\amalg 2} \ar[r,"{(f_0,f_1)}"] \ar[d,"{(F_3^+(s_0),F_3^+(s_1))}" left]	&	\zO^{\o,\natural}\\
			F_3^+(K \times \mathcal{J}) \ar[ru,dashed, bend right,"\bar{f}"]					& 	
		\end{tikzcd}
		\]
		and by lemma \ref{lm:F_and_G_preserve_cylinder}, the map $F_3^+(K\times \mathcal{J})\to F_3^+(K)$ is a marked equivalence, so that condition $(i_\approx)$ is satisfied. 
		\item[\framebox{$\sim\,\Leftarrow\, \approx$}] Suppose that $f_0\approx f_1$. We will show that condition $(i_\sim)$ holds. The factorization of $(\id_K,\id_K)$ through $q\colon K\times\mathcal{J}\to K$ induces a factorization
		\[
		\left(\id_{F_3^+(K)},\id_{F_3^+(K)}\right) \colon F_3^+(K)\coprod_{\Delta^1}F_3^+(K){\too} F_3^+(K\times\mathcal{J}) \too F_3^+(K)
		\]
		where the first map is a monomorphism. We can thus apply assumption $(ii_{\approx})$ to obtain a lift $\bar{f}$ in the diagram
		\[
		\begin{tikzcd}  
			F_3^+(K)^{\amalg 2} \ar[r,"{(f_0,f_1)}"] \ar[d,"{(F_3^+(s_0),F_3^+(s_1))}" left]	&	\zO^{\o,\natural}.\\
			F_3^+(K \times \mathcal{J}) \ar[ru,dashed, bend right,"\bar{f}"]					& 	
		\end{tikzcd}
		\]
		It now suffices to prove that $\bar{f}$ satisfies conditions \hyperlink{text:condition_F_n}{$(\star)_{3,\sigma}$}, which is a consequence of the fact that $f_0$ and $f_1$ both do and $(K\times\mathcal{J})_0\cong (K^{\coprod 2})_0$.
	\end{enumerate}

	This shows the above claim and therefore  completes our proof of proposition \ref{prop:p_is_equivalence}.
\end{proof}

\begin{proof}
	[Proof of theorem \ref{thm:comparaison_bo_ext}] 
	Combining lemma \ref{lm:i_n_induce_bijections_on_pi_0}, proposition \ref{prop:i_0_marked_anodyne}, lemma \ref{lm:i_1_inner_anodyne} and proposition \ref{prop:p_is_equivalence}, we obtain that each map in the zigzag \eqref{eqn:zigzag_comparaison_categories} is an equivalence, so that $\Ext(\sigma)$ and $\bo_\sigma$ are equivalent Kan complexes.
\end{proof}

	\section{Homotopy type of spaces of extensions} 
\label{sec:homotopy_type_of_spaces_of_extensions}

\subsection{Statement of the results} 
\label{sec:statement_of_the_results_cha_homotopy_type}

\subsubsection{Motivation}
\label{sub:motivation_cha_homotopy_type_spaces_extensions}

Let $\zO^\o$ be a unital $\infty$-operad, that we now assume to be monochromatic and such that the underlying $\infty$-category is an $\infty$-groupoid. Consider an operation $\sigma \in \zO(n)$ of arity $n$.

In the previous section, we provided a zigzag of equivalences between two models for the space of extensions of $\sigma$: on the one hand, the fiber $\bo_\sigma$ of Mann--Robalo's brane fibration $\pi\colon \bo\to \Tw(\Env(\zO))^\o$ and on the other hand, Lurie's space $\Ext(\sigma)$. 
One problem remains: neither of these two models is suitable for applications, since  computing their homotopy type of spaces of extensions is difficult \emph{even for simple examples of $\infty$-operads}.

Luckily, there is a third possible model for the space of extensions of the operation $\sigma$, namely the homotopy fiber
$$\zO(n+1)\mathop{\times}\limits^{\mathrm{h}}_{\zO(n)} \{\sigma\}$$ 
of the morphism $i^*\colon \zO(n+1)\to \zO(n)$ given by precomposition with a \emph{chosen atomic map  $i$}. We call it Toën's model for spaces of extensions.

In this section, we compare this new model to Lurie's space $\Ext(\sigma)$. We will show that, contrary to what one might expect, the models are not equivalent in general, unless the space $\zO(1)$ of unary operations is contractible. More precisely, we will prove the following result.

\begin{thm}[Theorem {\ref{thm:quotient_par_O1_intro}}]
	\label{thm:quotient_par_O(1)}
	Let $\zO^\o$ be a monochromatic unital $\infty$-operad whose underlying $\infty$-category $\zO$ is an $\infty$-groupoid and $\sigma\in \zO(n)$ an operation of arity $n$. Choose an atomic morphism  $i\colon \langle n \rangle \to \langle n+1\rangle$ in $\zO^\o$. Then there is a homotopy cartesian square
	\begin{equation}
		\label{eqn:homotopy_cartesian_square_thm_BO(1)}
		\begin{tikzcd}  
			\zO(n+1)\mathop{\times}\limits^{\mathrm{h}}_{\zO(n)} \{\sigma\} \ar[r,""] \ar[d,"" left]	&	\Ext(\sigma) \ar[d,""] \\
			* \ar[r,""]					& 	B\zO(1),
		\end{tikzcd}
	\end{equation}
	well-defined in the homotopy category of spaces, which exhibits $\Ext(\sigma)$ as a homotopy quotient of $\zO(n+1)\mathop{\times}\limits^{\mathrm{h}}_{\zO(n)}\{\sigma\}$ by a $\zO(1)$-action.
\end{thm}

\begin{rk}
	\label{rk:generalization_to_coloured_case}
	The restriction to the monochromatic situation is merely there to make the comparison with Toën's model more transparent and to slightly simplify the notations. The results of this section readily generalize to the general case of (possibly coloured) unital $\infty$-operads.
	\todo{écrire l'énoncé général dans ce cas}
\end{rk}

\begin{nota}
    \label{nota:abus_notation_colors_comme_dans_Fin}
    Throughout this section, when considering a monochromatic $\infty$-operad $\zO^\o$ with unique color $c$, we shall use the slightly abusive notation of writing objects of $\zO^\o$ in the form $\langle n \rangle$, instead of say $c^{\oplus n}$.
\end{nota}


\subsubsection{Difference with the existing literature}
\label{sub:difference_existing_literature}

Theorem \ref{thm:quotient_par_O(1)} contradicts the  statement \cite[Remark 5.1.1.10]{ha},  
in the case of $\infty$-operads with non-contractible spaces of unary operations. 
This statement  is a key result in Lurie's proof of coherence of the little disks $\infty$-operad $\bE_n$. 
Note that this statement is, however, \emph{only used} for this example of $\bE_n$, which satisfies that $\bE_n(1)\simeq *$ so that our theorem \ref{thm:quotient_par_O(1)} actually implies that the  conclusion of Lurie's statement is true, \emph{in this particular case}.

We then provide through proposition \ref{prop:example_computation_rmk_HA_fails} an explicit example of an $\infty$-operad for which this statement is incorrect, without appealing to the above theorem \ref{thm:quotient_par_O(1)}.

Let us now explain Lurie's statement.
Let $\zO_\Delta$ be a unital fibrant {simplicial} operad, with underlying $\infty$-operad $\zO^\o = \zN^\o(\zO_{\Delta})$. By unitality, the canonical inclusion $i_{\Fin_*}\colon \langle m\rangle \to \langle m+1\rangle$ in $\Fin_*$ lifts uniquely to a morphism $i$ in the simplicial category $\zO_{\Delta}^\o$, which induces a map of simplicial sets
\[
	i^*\colon 
	\Map^\mathrm{act}_{\zO_{\Delta}^\o} (\langle m+1\rangle, \langle n\rangle)\to
	\Map^\mathrm{act}_{\zO_{\Delta}^\o} (\langle m\rangle, \langle n\rangle).
\]
Given an active morphism $f\colon \langle m \rangle\to \langle n \rangle$ in $\zO_{\Delta}^\o$, Lurie defines in \cite[Notation 5.1.1.8.]{ha} the \emph{space of strict extensions of $f$}, denoted $\Ext_\Delta(f)$, as the fiber of $i^*$ at $f$.
Now consider an $n$-simplex $\sigma$ of $\zO^\o$ corresponding to a sequence of $n$ composable active morphisms
\[
	\langle m_0\rangle \stackrel{f_1}{\too}
	\langle m_1\rangle \stackrel{f_2}{\too}
	\dots \stackrel{f_n}{\too}
	\langle m_n\rangle
\]
in $\zO_\Delta^\o$. In \cite[Construction 5.1.1.9.]{ha}, a comparison map $\theta\colon \Ext_\Delta(f_n)\to \Ext(\sigma)$ is defined. Then \cite[Remark 5.1.1.10]{ha} asserts that $\theta$ can be identified with the canonical map $\mathrm{fib}_{f_n}(i^*)\to \mathrm{hofib}_{f_n}(i^*)$. In particular, $\Ext(\sigma)$ is supposed to be equivalent to the fiber of $i^*$ at $f_n$.
However, as a consequence of theorem \ref{thm:quotient_par_O(1)}, this equivalence can only hold when the group of unary operations $\zO(1)$ is trivial.

For a direct counterexample to \cite[Remark 5.1.1.10]{ha} for $\zO(1)\not\simeq *$, consider the operad $\zO_\Delta = \mathrm{AssInv}$ encoding associative algebras together with an involution. It is the monochromatic operad in sets freely generated by two operations $\mu \in \zO_\Delta(2)$ and $\tau\in \zO_\Delta(1)$ satisfying the relations $\tau^2 = \id$ and $\tau\circ \mu(a,b) = \mu(\tau b,\tau a)$; it can also be described as a semi-direct product $\mathrm{Ass} \rtimes \Sigma_2$. 
Its homotopy coherent nerve $\zO^\o = \zN(\zO_\Delta^\o)$ is a unital monochromatic discrete $\infty$-operad, in which morphisms from $\langle m\rangle$ to $\langle n \rangle$ are given by a map $\alpha \colon \langle m\rangle \to \langle n \rangle$ in $\Fin_*$, a linear order on each preimage $\alpha^{-1}\{i\}$ and a choice of sign $\varepsilon \colon \langle m\rangle\to \{+,-\} $. 
The previous morphism will be denoted more compactly 
$$\left( k_{i1}^{\varepsilon(k_{i1})} \dots k_{im_i}^{\varepsilon(k_{im_i})} \right)_{i\in \langle n\rangle^\circ}, \text{ with } \alpha\inv\{i\} = \{k_{i1} < \dots < k_{im_i}\}.$$
Composition is defined so that negative signs reverse the linear order.

\begin{prop}
    \label{prop:example_computation_rmk_HA_fails}
    For $\zO^\o = \zN(\mathrm{AssInv})^\o$ the $\infty$-operad of associative algebras with involution and $\sigma$ the identity operation on the unique color ${\langle 1\rangle}$,
    the spaces of extensions $\Ext(\sigma)$ and that of strict extensions $\Ext_\Delta(\sigma)$ are \emph{not} homotopy equivalent.
\end{prop}
\begin{proof}
    On the one hand, since $\zO_\Delta$ is a discrete simplicial operad, the homotopy fiber $\zO_\Delta(2)\mathop{\times}\limits^{\mathrm{h}}_{\zO_\Delta(1)} \{\sigma\}$ coincide with the actual fiber $\Ext_\Delta(\sigma)$, which is the $4$-elements set $\{\mu(a,b), \mu(a,\tau b), \mu(b,a), \mu(\tau b,a)\}$. 

    On the other hand, we claim that the set of connected components $\pi_0\Ext(\sigma)$ consists of only two elements. To prove this, recall the description of $k$-simplices of $\Ext(\sigma)$ given in \eqref{eqn:description_smplices_Ext_discrete_case}. 
    In particular, the objects of $\Ext(\sigma)$ are given by commutative squares
	\begin{equation}
        \begin{tikzcd}
            \langle 1\rangle \arrow[r, "\sigma=\id"] \arrow[d, "i"left, hook, harpoon]  & \langle 1\rangle \arrow[d, "\sim", "\varphi"left] \\
            \langle 2\rangle \arrow[r, "\alpha"]         & \langle 1\rangle.
        \end{tikzcd}
    \end{equation}
    with $i$ atomic, $\alpha$ active and $\varphi$ an equivalence.
    The morphisms in $\Ext(\sigma)$ are given by diagrams
	\begin{equation}
        \begin{tikzcd}
            \langle 1\rangle \arrow[r, "\sigma=\id"] \arrow[dd, "i'"left, hook, harpoon, bend right=50]\arrow[d, "i"left, hook, harpoon]  & \langle 1\rangle \arrow[dd, bend left=50, "\sim", "\varphi'"left]\arrow[d, "\sim", "\varphi"left] \\
            \langle 2\rangle \arrow[r, "\alpha"]\ar[d,"f"]& \langle 1\rangle \ar[d,"\sim","\psi"left]\\
            \langle 2\rangle \arrow[r, "\alpha'"]         & \langle 1\rangle
        \end{tikzcd}
    \end{equation}
    with active morphisms and with $f$ compatible with extensions.
    Consider two objects $x = (i,\varphi,\alpha)$ and $x' = (i',\varphi',\alpha')$. There exists a unique morphism $\psi$ such that $\psi\varphi=\varphi'$, whereas there are always two distinct morphisms $f$ that are compatible with extensions and satisfy $fi=i'$. 
    For this data $(f,\psi)$ to define a morphism in $\Ext(\sigma)$, we further need equation $\alpha' f = \psi \alpha$ to be satisfied. But since there are $4$ active morphisms $\langle 2\rangle \to \langle 1\rangle$ (namely $1^+2^+$, $1^+2^-$, $2^+1^+$ and $2^+1^-$), only half of the pairs $(x,x')$ are in the same connected component; this concludes the computation.

\end{proof}




\subsection{Auxiliary models for \texorpdfstring{$\Ext(\sigma)$}{Ext(sigma)} and the homotopy fiber of \texorpdfstring{$\zO(n+1)\to\zO(n)$}{O(n+1)->O(n)}} 
\label{sub:alternative_models_for_Ext_and_O(n+1)x_O(n)_sigma}

The first step is to represent the map $i^*\colon \zO(n+1)\to \zO(n)$, which is only well-defined in the homotopy category of spaces, by a zigzag of spaces 
$$\zO(n+1)\stackrel{\sim}{\longleftarrow} \widetilde{\zO(n+1)}\too \zO(n).$$
This will give a strict model of the homotopy fiber $\zO(n+1)\mathop{\times}\limits^{\mathrm{h}}_{\zO(n)}\{\sigma\}$.

Consider the subsimplicial set 
$$\Delta^2_2 = \Delta^{01}\coprod \Delta^{\{2\}}$$
of $\Delta^2$. Using that $\Lambda_1^2$ can be written as the pushout
$
\Lambda^2_1 = \Delta^2_2 \mathop{\coprod}\limits_{\partial\Delta^{12}} \Delta^{12},
$
we obtain that the mapping space $\zO(n+1)$ is isomorphic to the fiber
\begin{eqnarray*}
\zO(n+1) & := & \Fun(\Delta^{12},\oact) \mathop{\times}\limits_{\Fun(\partial\Delta^{12},\oact)} \left\lbrace(\langle n+1\rangle, \langle 1\rangle)\right\rbrace\\
& \cong & \Fun(\Lambda^2_1,\oact) \mathop{\times}\limits_{\Fun(\Delta^2_2,\oact)} \left\lbrace(i,\langle 1\rangle)\right\rbrace.
\end{eqnarray*}
Similarly, writing $\Lambda_0^2$ as $\Delta^2_2 \mathop{\coprod}\limits_{\partial\Delta^{02}}\Delta^{02}$, we can identify $\zO(n)$ as the fiber
\[
\zO(n)\cong \Fun(\Lambda^2_0,\oact) \mathop{\times}\limits_{\Fun(\Delta_2^2,\oact)} \left\lbrace(i,\langle 1\rangle)\right\rbrace.
\]
Now let $\widetilde{\zO(n+1)}$ denote the fiber of the restriction 
$$\Fun(\Delta^2,\oact)\too \Fun(\Delta^2_2,\oact)$$
at $(i,\langle 1\rangle)$. 
Since $\Delta^2_2 \inc \Delta^2$ factors through both $\Lambda_0^2$ and $\Lambda_2^2$, we obtain a commutative diagram
\begin{equation}
	\label{eqn:diag_O(n+1)tilde_O(n)}
	\begin{tikzcd}
		\zO(n+1) \ar[d,"{(i,\mathrm{id})}"] & 
		\widetilde{\zO(n+1)} 
		\arrow[ld, "\llcorner" near start, phantom]
		\arrow[rd, "\lrcorner" near start, phantom]
		\ar[d]\ar[r,two heads]\ar[l,"\sim", two heads] &
		\zO(n)\ar[d,"{(i,\mathrm{id})}"]\\
		\Fun(\Lambda_1^2, \oact) &
		\Fun(\Delta^2, \oact)\ar[r,two heads]\ar[l,"\sim", two heads] &
		\Fun(\Lambda_0^2, \oact).
	\end{tikzcd}
\end{equation}

\begin{lm}
	\label{lm:O(n+1)x_O(n)_as_strict_pullback}
	The above diagram yields equivalences of fibers
	$$\zO(n+1)\mathop{\times}\limits^{\mathrm{h}}_{\zO(n)} \{\sigma\}
	\:\simeq\: \widetilde{\zO(n+1)}\mathop{\times}\limits_{\zO(n)} \{\sigma\}
	\:\cong\: \Fun(\Delta^2, \oact)\mathop{\times}\limits_{\Fun(\Lambda_0^2, \oact)} \{(i,\sigma)\}.$$
\end{lm}

\begin{proof}
	In diagram \eqref{eqn:diag_O(n+1)tilde_O(n)}, observe that the top row is obtained as the fiber of the bottom row at the point $(i,\langle 1\rangle)$ of $\Fun(\Delta^2_2,\oact)$. Therefore both squares in the diagram are cartesian. In particular, since the bottom right morphism is a Kan fibration, so is $\widetilde{\zO(n+1)} \to \zO(n)$. Similarly, since the bottom left morphism is a trivial Kan fibration, so is $\widetilde{\zO(n+1)} \to \zO(n+1)$. This gives the first homotopy equivalence of the lemma. The second equivalence, which is an isomorphism of simplicial sets, is obtained by taking the fiber of the right cartesian square at the object $\sigma \in \zO(n)$, whose image in $\Fun(\Lambda_0^2, \oact)$ is $(i,\sigma)$.
\end{proof}

We now replace $\Ext(\sigma)$ with a more convenient model, that we shall denote $\Ext^\square(\sigma)$, defined as a certain subcategory of the functor $\infty$-category $\Fun(\Delta^1\times\Delta^1,\oact)$ of commutative squares of active maps in $\zO^\o$. The squares will be indexed as follows:
\[
\begin{tikzcd}  
	00 \ar[r,""] \ar[d,"" left]	&	01 \ar[d,""] \\
	10 \ar[r,""]					& 	11.
\end{tikzcd}
\]
We will require that the left vertical map is atomic, the right vertical one is an equivalence and the top one is precisely the fixed morphism $\sigma$. In order to define $\Ext^\square(\sigma)$, the following notations will be convenient.

\begin{nota}
	Let $\atom$ denote the non-full subcategory of $\Fun(\Delta^1,\oact)$ whose objects are atomic morphisms (see definition \ref{df:semi-inert_atomic_maps}).
	Let $\atom(n)$ be the full subcategory of $\atom$ whose objects are maps with codomain $\langle n \rangle$.
\end{nota}

\begin{nota}
	The marked simplicial set obtained from the square $\Delta^1\times\Delta^1$ by further marking the edge $\Delta^1\times\{1\}$ will be denoted $\square$.
\end{nota}

\begin{df}
	[Definition of $\Ext^\square(\sigma)$]
	\label{df:Ext^square}
	Define the simplicial set $\Ext^\square(\sigma)$ as the iterated fiber product
	\[
	\Ext^\square(\sigma) = \lim \left(
		\begin{tikzcd}[column sep = 0.4em]
		\atom \ar[d,hook] &
		\Map^\flat(\square, (\oact)^\natural) \ar[ld]\ar[rd] &
		\{\sigma\} \ar[d] \\
		\Fun(\Delta^1\times\{0\},\oact) & & 
		\Fun(\{0\}\times\Delta^1,\oact)
		\end{tikzcd}
	\right)
	\]
	where the two diagonal morphisms are the obvious restriction maps.

In other words, $\Ext^\square(\sigma)$ is the subcategory of $\Fun(\Delta^1\times\Delta^1,\oact)$ whose
\begin{itemize}
	\item \textbf{objects} are commutative diagrams
	\[
	\begin{tikzcd}  
		\langle n \rangle \ar[r,"\sigma"] \ar[d,hook,harpoon,"\mathrm{atom}" left]	&	\langle 1 \rangle \ar[d,"\wr"] \\
		\langle n+1 \rangle \ar[r,""]					& 	\langle 1 \rangle 
	\end{tikzcd}
	\]
	in which the left vertical map is atomic and the right vertical map is an equivalence,
\item \textbf{morphisms} are compatible with extensions, i.e. preserve the new color $\langle n+1\rangle \setminus\mathrm{Im}(\langle n\rangle)$.
\end{itemize}
\end{df}

It is easy to see that $\Ext^\square(\sigma)$ is an $\infty$-category.



\begin{lm}
	\label{lm:ext_equivalent_Ext^square}
	The space $\Ext(\sigma)$ is equivalent to $\Ext^{\square}(\sigma)$.
\end{lm}
\begin{proof}
	By inspection of definition \ref{df:Ext(sigma)}, one easily verifies that all diagrams involved in the definition contains only active maps. Therefore diagrams to $\Ext(\sigma)$ factor through the subcategory $\Fun(\Delta^1,\oact)_{\sigma/}$ of $\Fun(\Delta^1,\zO^\o)_{\sigma/}$.
	Now recall the canonical equivalence of $\infty$-categories
	\[
	\gamma\colon 
	\Fun(\Delta^1,\oact)_{\sigma/} \stackrel{\sim}{\too} 
	\Fun(\Delta^1,\oact)^{\sigma/}
	\]
	from the slice to the alternative slice, the latter being defined as the fiber at $\sigma$ of the restriction map
	\[
	\Fun(\Delta^1\times\Delta^1,\oact)\too 
	\Fun(\{0\}\times\Delta^1,\oact).
	\]
	The restriction of $\gamma$ to $\Ext(\sigma)$ factors through the obvious inclusion $\Ext^\square(\sigma)\to \Fun(\Delta^1,\oact)^{\sigma/}$. Moreover, by inspection of the objects of these two $\infty$-categories, one sees that this functor $\gamma|_{\Ext(\sigma)}\colon \Ext(\sigma)\to \Ext^\square(\sigma)$ is essentially surjective. To prove the lemma, it therefore suffices to show fully faithfullness of $\gamma|_{\Ext(\sigma)}$. Given two extensions $X,X'\in \Ext(\sigma)$, consider the commutative diagram
	\[
	\begin{tikzcd}  
		\Map_{\Ext(\sigma)}(X,X') \ar[r,""] \ar[d,"\gamma|_{\Ext(\sigma)}" left]	&	\Map_{\Fun(\Delta^1,\oact)_{\sigma/}}(X,X') \ar[d,"\gamma" right, "\wr" left] \\
		\Map_{\Ext^{\square}(\sigma)}(\gamma(X),\gamma(X')) \ar[r,""]&	\Map_{\Fun(\Delta^1,\oact)^{\sigma/}}(\gamma(X),\gamma(X')).
	\end{tikzcd}
	\]
	As $\gamma$ is an equivalence of $\infty$-categories, the right vertical map is a homotopy equivalence.
	Observe that, given two equivalent morphisms $f_0\simeq f_1 \colon X\to X'$ in $\Fun(\Delta^1,\oact)_{\sigma/}$, $f_0$ is compatible with extensions if and only if $f_1$ has this property, and similarly for morphisms $\gamma(X)\to\gamma(X')$.
	Consequently, the horizontal maps in the above diagram are both inclusions of the connected components determined by the condition of preservation of the new color in the extensions. Therefore the restriction $\gamma|_{\Ext(\sigma)}$ is a homotopy equivalence, as desired.
\end{proof}

To compare $\Ext^\square(\sigma)$ with $\zO(n+1)\mathop{\times}\limits^{\mathrm{h}}_{\zO(n)} \{\sigma\}$, we first give an alternative description of the former $\infty$-groupoid.
By definition of $\atom(n)$, we have a commutative diagram 
\begin{equation}
	\label{eqn:diag_atom(n)}
	\begin{tikzcd}
	{\atom(n)} 		\arrow[rd, "\lrcorner" near start, phantom]\ar[rd,"j",bend right=20,dotted, hook]
	\arrow[rr] \arrow[dd] &                        & {\{\sigma\}} \arrow[d] \\
	                         & {\Fun(\Lambda_0^2, \oact)} 
	                     		\arrow[rd, "\lrcorner" near start, phantom]
	                         	\arrow[r] \arrow[d] & {\Fun(\{0\}\times\Delta^1, \oact)} \arrow[d] \\
	{\atom} \arrow[r]             & {\Fun(\Delta^1\times\{0\}, \oact)} \arrow[r]           & {\Fun(\{0\}\times\{0\}, \oact)}          
	\end{tikzcd}
\end{equation}
in which both squares are cartesian, using the identification 
\begin{equation}
	\label{eqn:identification_Lambda_0^2_as_square}
	\Lambda_0^2 =  \Delta^{01}\cup\Delta^{02}
	\cong 
	\Delta^1\times\{0\}\cup \{0\}\times\Delta^1,
\end{equation}
and the map $j$ is induced by the universal property of pullbacks.
\begin{lm}
	\label{lm:Ext^square_as_fiber_of_Atom(n)}
	There is a canonical isomorphism 
	\begin{equation}
		\label{eqn:ext^square_=_fiber_of_Atom(n)}
		\Ext^\square(\sigma) \cong \Map^\flat(\square, (\oact)^\natural) \mathop{\times}\limits_{\Fun(\Lambda_0^2, \oact)} \atom(n).
	\end{equation}
\end{lm}

\begin{proof}
	Consider the following commutative diagram, extending diagram \eqref{eqn:diag_atom(n)}:
	\begin{equation*}
		\begin{tikzcd}[cramped, scale cd = 0.7, column sep = 0em, row sep = scriptsize]
		\Ext^\square(\sigma) 		\arrow[rrdd, "\lrcorner" very near start, phantom]
		\arrow[dd] \arrow[rr] \arrow[rd]    &                                & \Map^\flat(\square,\oact)\mathop{\times}\limits_{{(\oact)^{\{0\}\times\Delta^1}}}\{\sigma\} \arrow[rrrddd, "\lrcorner" very near start, phantom]\arrow[rd] \arrow[rrrd] \arrow[dd] &                                                 &  &                                      \\
		                                                         & \atom(n) 		\arrow[rrdd, "\lrcorner" very near start, phantom]\arrow[dd] \arrow[rr] &
	& { {(\oact)^{\Delta^1\times\{0\}}}\mathop{\times}\limits_{\oact} \{\sigma\}} \arrow[rrdd, "\lrcorner" very near start, phantom]\arrow[rr] \arrow[dd]                        &  & \{\sigma\} \arrow[dd]                \\
		{\atom\mathop{\times}\limits_{{(\oact)^{\Delta^1\times\{0\}}}}\Map^\flat(\square,\oact)} \arrow[rrrddd, "\lrcorner" very near start, phantom] \arrow[rd] \arrow[rddd] \arrow[rr] &                                & \Map^\flat(\square,\oact) \arrow[rd] \arrow[rrrd] \arrow[rddd]               &                                                 &  &                                      \\
		                                                         & {\atom\mathop{\times}\limits_{\oact}{(\oact)^{\{0\}\times\Delta^1}}}\arrow[rrdd, "\lrcorner" very near start, phantom]
 \arrow[rr] \arrow[dd]       &                                                               & {(\oact)^{\Lambda_0^2}} 	\arrow[rrdd, "\lrcorner" very near start, phantom]\arrow[rr] \arrow[dd] &  & {(\oact)^{\{0\}\times\Delta^1}} \arrow[dd] \\
		                                                         &                                &                                                               &                                                 &  &                                      \\
	& \atom \arrow[rr] & & {(\oact)^{\Delta^1\times\{0\}}} \arrow[rr]  &  & {(\oact)^{\{0\}\times\{0\}}}.                 
		\end{tikzcd}
	\end{equation*}
	In the above, certain squares are cartesian by construction, namely:
		\begin{itemize}
			\item all the  squares whose arrows are either vertical or horizontal
			\item the two squares that contains $\Map^\flat(\square,\oact)$ and either $\{\sigma\}$ or $\atom$.
		\end{itemize}
	Using the usual transitivity rule for pullback squares, one deduce that any square in the top left cube is cartesian, from which the desired isomorphism follows.
\end{proof}

\subsection{Proof of theorem \ref{thm:quotient_par_O(1)}} 
\label{sec:proof_of_theorem_thm:quotient_par_o(1)}

There are two differences between the homotopy fiber of $\zO(n+1)\to \zO(n)$ at $\sigma$, modelled as $\mathrm{fib}_{(i,\sigma)}(\Fun(\Delta^2, \oact)\to \Fun(\Lambda_0^2, \oact)$, and $\Ext^\square(\sigma)$:
\begin{enumerate}
	\item objects of $\mathrm{hofib}_{\sigma}(\zO(n+1)\to\zO(n))$ are given by commutative \emph{triangles} in $\oact$, whereas objects of $\Ext^\square(\sigma)$ are commutative \emph{squares},
	\item in $\mathrm{hofib}_{\sigma}(\zO(n+1)\to\zO(n))$, the map $\langle n\rangle \to \langle n+1\rangle$ is the fixed morphism $i$ whereas in $\Ext^\square(\sigma)$, any atomic map is allowed.
\end{enumerate}
As we will see, the first difference does not affect the homotopy type of the spaces, but the second difference explains the origin of the quotient by the action of $\zO(1)$. To make this remark precise, we will introduce variants of $\Ext^\square(\sigma)$ that differ according to the previous two parameters.\\

Consider the morphism $r\colon \Delta^1\times\Delta^1\to \Delta^2$ induced from the map of posets $[1]\times[1]\to [2]$ given by
\[
r(0,0) = 0, \quad
r(0,1) = 1, \quad
r(1,0) = 1, \quad
r(1,1) = 2.
\]
It induces a morphism of marked simplicial sets $\square\to (\Delta^2)^\flat$, which yields a restriction map $r^*\colon {\Fun(\Delta^2,\oact)}\to \Map^\flat(\square,\oact)$. Moreover, $r$ extends the identification \eqref{eqn:identification_Lambda_0^2_as_square} to a commutative square
\[
\begin{tikzcd}  
	\{0\}\times\Delta^1\cup \Delta^1\times\{0\} \ar[r,"\cong"] \ar[d,hook]	&	\Lambda_0^2 \ar[d, hook] \\
	\Delta^1\times\Delta^1 \ar[r,"r"]					& 	\Delta^2.
\end{tikzcd}
\]

\begin{df}
	Let	$\Ext^\triangle(\sigma)$ and $\Ext^\square(\sigma,i)$ be the $\infty$-categories fitting in the following diagram of cartesian squares
	\begin{equation}
		\label{eqn:diag_def_variants_of_Ext(sigma)}
		\begin{tikzcd}
		\zO(n+1)\mathop{\times}\limits_{\zO(n)}^{\mathrm{h}}\{\sigma\}\ar[rd, "\lrcorner" very near start, phantom] \arrow[r] \arrow[d] & \Ext^\triangle(\sigma)\ar[rd, "\lrcorner" very near start, phantom] \arrow[r] \arrow[d] & {\Fun(\Delta^2,\oact)} \arrow[d,"r^*"]      \\
		{\Ext^\square(\sigma,i)}\ar[rd, "\lrcorner" very near start, phantom] \arrow[r] \arrow[d]                                       & \Ext^\square(\sigma)\ar[rd, "\lrcorner" very near start, phantom] \arrow[r] \arrow[d]   & {\Map^\flat(\square,\oact)} \arrow[d] \\
		\ast \arrow[r, "i"] \arrow[rr, "\sigma", bend right=25]                              & \atom(n) \arrow[r,"j"]                         & {\Fun(\Lambda_0^2,\oact)}.       
		\end{tikzcd}
	\end{equation}
\end{df}
The fact that $\zO(n+1)\mathop{\times}\limits_{\zO(n)}^{\mathrm{h}}\{\sigma\}$ and $\Ext^\square(\sigma)$ fit in the above diagram is a reformulation of lemmas \ref{lm:O(n+1)x_O(n)_as_strict_pullback} and \ref{lm:Ext^square_as_fiber_of_Atom(n)}.

\begin{lm}
	\label{lm:equivalences_triangle_carres}
	In the top left square of diagram \eqref{eqn:diag_def_variants_of_Ext(sigma)}
	\begin{equation}
	\begin{tikzcd}
		\zO(n+1)\mathop{\times}\limits_{\zO(n)}^{\mathrm{h}}\{\sigma\}\ar[rd, "\lrcorner" very near start, phantom] \arrow[r] \arrow[d] & 
		\Ext^\triangle(\sigma)\arrow[d]   \\
		{\Ext^\square(\sigma,i)}	\arrow[r]  & \Ext^\square(\sigma),
	\end{tikzcd}
	\end{equation}
	the vertical maps are equivalences.
\end{lm}
\begin{proof}
	A simple computation shows that $r$ is an equivalence of marked simplicial sets $\square = (\Delta^1\times\Delta^1,\Delta^1\times\{1\})\to (\Delta^2)^\flat$. Therefore $r^*$ is an equivalence of $\infty$-categories. Since $\Map^\flat(\square,\oact)$ and $\Fun(\Delta^2,\oact)$ are fibrant over $\Fun(\Lambda_0^2,\oact)$, taking pullback along the morphisms $j$ and $\sigma\colon \ast\to \Fun(\Lambda_0^2,\oact)$ gives the desired equivalences.
\end{proof}

\begin{lm}
	\label{lm:atom(n)_equivalent_O}
    Let $\zO^\o$ be a monochromatic \todo{[preuve à réécrire plus tard] monochromatic ? en tout cas, réécrire la preuve pour avoir le cas général ici} unital $\infty$-operad.
	Then the $\infty$-category $\atom(n)$ is equivalent to the underlying $\infty$-category $\zO$ of $\zO^\o$.
\end{lm}

\begin{cor}
	\label{cor:atom(n)_equivalence_BO(1)}
	Let $\zO^\o$ be as above and assume moreover that its underlying $\infty$-category $\zO$ is an $\infty$-groupoid.
    Then the $\infty$-category $\atom(n)$ is equivalent to the classifying space $\mathrm{B}\zO(1)$ of the group of automorphisms of $\langle 1\rangle$ in $\zO$.
\end{cor}

\begin{proof}[Proof of lemma \ref{lm:atom(n)_equivalent_O}]
	First, fix an atomic morphism $\alpha\colon \langle n\rangle \to \langle n+1\rangle$ in $\Fin_*$ and consider the subcategory $\mathrm{Atom}_\zO^\alpha(n)$ of morphisms lying over $\alpha$. By definition, we have a cartesian square
	\[
	\begin{tikzcd}  
		\mathrm{Atom}_\zO^\alpha(n) \ar[r,""]\arrow[rd, "\lrcorner" very near start, phantom]\ar[d,"" left]	&	\atom(n) \ar[d,""] \\
		\{\alpha\} \ar[r,""]					& 	\mathrm{Atom}_{\Fin_*}(n).
	\end{tikzcd}
	\]
	Observe that the $\infty$-category $\mathrm{Atom}_{\Fin_*}(n)$ is the nerve of a $1$-category in which any two objects $j\colon \langle n\rangle \to S$ and $j'\colon \langle n \rangle \to S'$ are related by a unique morphism $S\to S'$ (namely the unique bijection that restricts to $j'\circ (j|^{\im(j)})\inv$ on the image of $j$). As a consequence, this ($\infty$-)category is terminal and we obtain a canonical equivalence of $\infty$-categories $\mathrm{Atom}_\zO^\alpha(n)\simeq \atom(n)$.

	Now we may decompose the atomic morphisms of $\mathrm{Atom}_\zO^\alpha(n)$ according to their arity using lemma \ref{lm:operad_decomposition}. The result is an equivalence of $\infty$-categories
	\[
	\mathrm{Atom}_\zO^\alpha(n) \simeq \left(\zO_{\langle 1\rangle /}\right)^n \times \zO_{\langle 0\rangle /}
	\]
	where the $\infty$-category $\zO_{\langle 0\rangle /}$ is a notation for the comma category 
	$$(\langle 0\rangle \in \zO^\o)\downarrow(\zO^\o \supset \zO).$$ 
	Since $\zO$ is assumed to be an $\infty$-groupoid, so is its slice $\zO_{\langle 1\rangle /}$; the latter has an initial object, it is therefore contractible. 
	We now turn to analysing the comma category  $\zO_{\langle 0\rangle /}$.
	By definition, it fits in a commutative diagram of cartesian squares
	\[
	\begin{tikzcd}  
		\zO_{\langle 0\rangle /} \ar[r,""] \arrow[rd, "\lrcorner" very near start, phantom]\ar[d,"" left]	&
		(\zO^\o)^{\langle 0\rangle /} \ar[d,""] \ar[r] \arrow[rd, "\lrcorner" very near start, phantom]
		& \Fun(\Delta^1,\zO^\o)\ar[d,"{(\mathrm{ev}_0,\mathrm{ev}_1)}"] \\
		\{\langle 0\rangle \} \times \zO \ar[r,hook]		& \{\langle 0\rangle \} \times 	\zO^\o \ar[r,hook] & \zO^\o \times \zO^\o .
	\end{tikzcd}
	\]
	Since the middle vertical map is a cocartesian fibration, so is the left vertical map. The fiber of this morphism at at object $X\in \zO$ is $\Map_{\zO^\o}(\langle 0\rangle, X)$, which is contractible by the assumption that $\zO^\o$ is unital. Therefore this cocartesian fibration is a trivial fibration $\zO_{\langle 0\rangle /} \stackrel{\simeq}{\to} \zO$, which completes the proof. 
    \todo{[DONE] solve the ambiguity of colors of $\zO^\o$ written as $\langle m \rangle$}
\end{proof}

\begin{proof}
    [Proof of Theorem \ref{thm:quotient_par_O(1)}] The homotopy cartesian square \eqref{eqn:homotopy_cartesian_square_thm_BO(1)} is obtained as the top left square in diagram \eqref{eqn:diag_def_variants_of_Ext(sigma)}, using the equivalences $\Ext^\triangle(\sigma)\simeq \Ext^\square(\sigma)\simeq \Ext(\sigma)$ of lemmas \ref{lm:equivalences_triangle_carres} and \ref{lm:ext_equivalent_Ext^square} and the equivalence $\mathrm{Atom}^\sigma_\zO(n)\simeq \mathrm{B}\zO(1)$ of corollary \ref{cor:atom(n)_equivalence_BO(1)}.
\end{proof}


\todo[inline]{change $\zO$ for $\zE$ when needed (?!)}
\todo[inline]{distinguer les endroits où l'on peut raccourcir la preuve, en l'ajoutant dans la thèse}

\section{Applications} 
\label{sec:applications}


As mentioned in the introduction, an important motivation for studying the brane action comes from string topology, as the $\bE_2$-algebra structure on free loop spaces arises from span diagrams given the brane action for the $\infty$-operad $\bE_2$.
More precisely, our work was motivated by the desire to use the formalism of brane actions to generalize string topology in the following two directions. 
\begin{itemize}
    \item On the one hand, we may consider analogs of free loop spaces $\Map(S^{n-1},X)$ based on higher dimensional spheres: this is the study of brane topology. Recall that the typical structure we expect in this setting is that of an $\bE_n^{\mathrm{fr}}$-algebra. Therefore, we would like to enhance the known $\bE_n$-structure of brane topology to take into account the action of groups of homeomorphisms of disks. The $\infty$-operads governing such structures are variants $\bE_n^G$ of the little disks $\infty$-operad $\bE_n$, that are given by semi-direct product of $\bE_n$ with a topological group $G$ endowed with a morphism to the $\infty$-group $\Top(n)$ of self-homeomorphisms of $\bR^n$.
    \item On the other hand, we want to obtain brane topology structures not just for the topological manifolds, but also for other geometric contexts, most notably for derived algebraic stacks. Note that this idea appears already in Toën's original work \cite{toen_branes}, where a sketch of proof of a higher formality theorem is given, that applies to quotients stacks of the form $[Y/G]$ for $G$ a smooth linear algebraic group over a field of characteristic zero acting on a quasiprojective derived scheme $Y$ of finite presentation.
\end{itemize}


The rest of this section is devoted to the study of $B$-framed little disks $\bE_B$, a family of $\infty$-operads that generalize the $\infty$-operads $\bE_n^G$ introduced above, a proof of their coherence and, as a consequence, a construction of new  operations on spaces of branes of perfect derived stacks. 

\subsection{Coherence of the \texorpdfstring{$\infty$}{infinity}-operad of little \texorpdfstring{$B$}{B}-framed disks} 
\label{sub:coherence_of_variants_of_the_little_disks_operad}

In this section, we define the $\infty$-operad of $B$-framed little disks and prove that it is coherent. This $\infty$-operad depends on the datum of a Kan fibration $B\to \mathrm{B}\Top(n)$ and recovers the variants $\bE_n^G$ of the little disks $\infty$-operad mentionned above when $B$ is the classifying space of a subgroup $G$ of $\Top(n)$.

%

\vspace{0.3cm}
We recall the definition of $\infty$-operad $\bE_B^\o$ introduced in \cite[Section 5.4.2]{ha}, following the presentation and  the notations thereof.

\begin{nota}
    Given two topological spaces $X$ and $Y$, we let $\Emb(X,Y)$ denote the topological space of open embeddings $X\to Y$, topologized as a subspace of the compact-open topology on the set $\Hom_{\Top}(X,Y)$. For $n\in \bN$, we let $\Top(n)$ denote the topological space of homeomorphisms of $\bR^n$, viewed as a subspace of $\Emb(\bR^n,\bR^n)$.
\end{nota}
\begin{rk}\label{rk:Kister_Mazur_thm}
    The Kister--Mazur theorem implies that the inclusion $\Top(n)\to \Emb(\bR^n,\bR^n)$ is a homotopy equivalence, for all $n\geq 0$ (see \cite[Theorem 5.4.1.5]{ha}).
\end{rk}

From now on, we fix a natural number $n$.

\begin{df}[{\cite[Definition 5.4.2.1]{ha}}]
    Let $^t\bE_{\mathrm{B}\Top(n)}^\o$ denote the topological  category whose objects are the finite pointed sets $\langle m \rangle \in \Fin_*$ and where mapping spaces are given by the formula
    \begin{equation}
        \Map_{^t\bE_{\mathrm{B}\Top(n)}^\o}(\langle m \rangle, \langle k \rangle)
        = \mathop{\coprod}\limits_{\alpha \in \Hom_{\Fin_*}(\langle m \rangle, \langle k \rangle)}
        \prod_{i = 1}^m
        \Emb(\bR^n\times \alpha\inv\{i\}, \bR^n).
    \end{equation}
    Let $\mathrm{B}\Top(n)^\o$ denote its homotopy coherent nerve, i.e. the $\infty$-category $N\left(^t\bE_{\mathrm{B}\Top(n)}^\o\right)$. 
    By \cite[Proposition 2.1.1.27]{ha}, $\mathrm{B}\Top(n)^\o$ forms an $\infty$-operad.
    We will denote by $\mathrm{B}\Top(n)$ its underlying $\infty$-category, which by remark \ref{rk:Kister_Mazur_thm} is a classifying space for the topological group $\Top(n)$.
\end{df}

Let us fix a Kan complex $B$ together with a Kan fibration $B \to \mathrm{B}\Top(n)$.

\begin{nota}
	Recall that given an $\infty$-category $\zC$, one can construct a cocartesian $\infty$-operad $\zC^\amalg$ whose spaces of multimorphisms are given by the formula
	$\Mul_{\zC^\amalg}(c_1,\dots, c_m;c) = \prod_{i=1}^m \Map_{\zC}(c_i,c)$ (see \cite[Section 2.4.3]{ha}).
\end{nota}

\begin{df}
	[{\cite[Definition 5.4.2.10]{ha}}]
	We let $\bE_B^\o$ denote the $\infty$-operad 
	\begin{equation}
		\bE_B^\o = \mathrm{B}\Top(n)^\o \mathop{\times}\limits_{\mathrm{B}\Top(n)^\amalg} B^\amalg
	\end{equation}
    and refer to it as the \emph{$\infty$-operad of $B$-framed little disks}.
\end{df}
Note that the underlying $\infty$-category of $\bE_B^\o$ is canonically equivalent to the Kan complex $B$. In particular,  one may identify the objects of $\bE_B^\o$ with those of $B$.

\begin{rk}[Examples]
    \begin{itemize}
        \item For $B$ a contractible Kan complex with a Kan fibration to $\mathrm{B}\Top(n)$, the associated $\infty$-operad $\bE_B^\o$ reduces to the ordinary $\infty$-operad $\bE_n^\o$  of little disks of dimension $n$.
        \item Consider a topological group together  with a map to $\Top(n)$. The induced morphism on classifying space can be represented up to equivalence by a Kan fibration $B := \mathrm{B}G\to \mathrm{B}\Top(n)$. Then the $\infty$-operad of $B$-framed little disks $\bE_B^\o$ is equivalent to a semi-direct $\bE_n^\o \rtimes G$ (in the sense of \cite{salvatore-wahl}). As a particular case, for $G = SO(n)$ we obtain the framed little disks $\infty$-operad $\bE_n^{\mathrm{fr}}$.
        \item  Let $M$ be a topological manifold of dimension $n$. Following \cite[Definition 5.4.5.1.]{ha}, let $\zC_M$ denote the topological category with two objects $M$ and $\bR^n$, whose mapping spaces are
            \begin{align*}
                \Map_{\zC_M}(\bR^n,\bR^n) &= \Emb(\bR^n,\bR^n)
                                          & 
                \Map_{\zC_M}(\bR^n,M)     &= \Emb(\bR^n,M)\\
                \Map_{\zC_M}(M,\bR^n) &= \emptyset
                                     & 
                \Map_{\zC_M}(M,M)    &= \{\id_M\}.
            \end{align*}
            Define $B_M$ as the Kan complex $\mathrm{B}\Top(n) \times_{\zN(\zC_M)} \zN(\zC_M)_{/M}$ and let $\bE_M^\o$ denote the $\infty$-operad $\bE_{B_M}^\o$. It is a variant on the $\infty$-operad $\bE_n^\o$ in which colors are open embedding $U\colon \bR^n\to M$ of disks of dimension $n$ into $M$ and operations of arity $k$ are diagrams of embeddings
            \begin{equation*}
                \begin{tikzcd}
                    \coprod_{i = 1}^k \bR^n \ar[rr]\ar[rd,"\coprod_i U_i"below left] && \bR^n \ar[ld,"U"]\\
                                            & M &
                \end{tikzcd}
            \end{equation*}
            together with an isotopy making the triangle commute.
            Note that $\bE_M$-algebras can be identified as locally constant factorization algebras on $M$, by theorem \cite[Theorem 5.4.5.9.]{ha}.
    \end{itemize}
\end{rk}


The $\infty$-operad of little disks operad $\bE_n$ is coherent, for all $n\in \bN$ by \cite[Theorem 5.1.1.1]{ha}. We generalize this result to the $B$-framed situation.
\begin{thm}
	\label{thm:B-framed_little_disks_coherent_operad}
	The $\infty$-operad of $B$-framed little disks $\bE_B^\o$ is coherent.
\end{thm}

The proof relies on theorem \ref{thm:quotient_par_O(1)} together with the following computation.
\begin{lm}
	\label{lm:computation_strict_version_extensions_E_B}
	Let $\sigma\colon (b_1,\dots,b_m)\to b$ be an active morphism in $\bE_B^\o$, with $b_1,\dots,b_m,b$ in $B$ and choose an additional color $b_{m+1}$ in $B$.
	Then Toën's model for the space of extensions of $\sigma$ in $\bE_B^\o$ is given by
	\begin{equation}
		\label{eqn:computation_toen_extension_E_B}
		\Mul_{\bE_B}(b_1,\dots,b_{m+1};b)\mathop{\times}\limits_{\Mul_{\bE_B}(b_1,\dots,b_m;b)}^\mathrm{h} \{\sigma\}
		\simeq 
		\begin{cases}
			 \Omega_b B\times \bigvee^m S^{n-1} & \text{ if } b_{m+1}\simeq b\text{ in } B,\\
			\emptyset & \text{ otherwise}.
		\end{cases}
	\end{equation}
\end{lm}

\begin{proof}
	By construction, the left hand side of equation \eqref{eqn:computation_toen_extension_E_B} is equivalent to the homotopy fiber at $\sigma$ of the map
	\begin{equation}
		\begin{tikzcd}
			\Emb(\bR^n\times \langle m+1\rangle^\circ,\bR^n) \mathop{\times}\limits_{\Emb(\bR^n,\bR^n)^{m+1}} \prod_{i=1}^{m+1} \Map_B(b_i,b)
			\ar[d] \\
			\Emb(\bR^n\times \langle m\rangle^\circ,\bR^n) \mathop{\times}\limits_{\Emb(\bR^n,\bR^n)^m} \prod_{i=1}^{m} \Map_B(b_i,b).
		\end{tikzcd}
	\end{equation}
	Commuting the fiber product with the homotopy fiber, we obtain the space
	\begin{equation}
		\left(
			\Emb(\bR^n\times\langle m+1\rangle^\circ, \bR^n)
		\mathop{\times}\limits_{
			\Emb(\bR^n\times\langle m  \rangle^\circ, \bR^n)
		}^\mathrm{h} \{\sigma\} \right)
		\mathop{\times}\limits_{
		\Emb(\bR^n,\bR^n)} \Map_B(b_{m+1},b)
	\end{equation}
	Since $B$ is a Kan complex, the factor $\Map_B(b_{m+1},b)$ is empty when $b_{m+1}$ and $b$ are in different connected components, and is equivalent to the based loop space $\Omega_b B$ otherwise.
    On the other hand, for any finite set $S$, the obvious map from the space $\Emb(\bR^n\times S, \bR^n)$  to the product $\Emb(\bR^n,\bR^n)^S\times \Conf(S,\bR^n)$ is an equivalence (see \cite[Proof of Proposition 5.4.2.8.]{ha}).
    As a consequence, we obtain an equivalence
	\begin{equation}
		\label{eqn:step_in_proof_lm_computation_strict_version_ext_E_B}
		\Emb(\bR^n\times\langle m+1\rangle^\circ, \bR^n)
		\mathop{\times}\limits_{
			\Emb(\bR^n\times\langle m  \rangle^\circ, \bR^n)
		}^\mathrm{h} \{\sigma\}
		\simeq 
	    \Emb(\bR^n,\bR^n)\times \Conf(S,\bR^n).
    \end{equation}
    Substituting this equivalence in \eqref{eqn:step_in_proof_lm_computation_strict_version_ext_E_B} and using that $\Conf(S,\bR^n)\simeq \bigvee^m S^{n-1}$,
 we obtain the result.
\end{proof}

\begin{lm}
	\label{lm:computation_extensions_E_B}
	Let $\sigma\colon (b_1,\dots,b_m)\to b$ be an active morphism in $\bE_B^\o$, with $b_1,\dots,b_m,b$ in $B$. 
	Then the space of extensions of $\sigma$ in $\bE_B^\o$ is equivalent to $\bigvee^m S^{n-1}$.
\end{lm}

\begin{proof}
    By theorem \ref{thm:quotient_par_O(1)}, any choice of a color $b_{m+1}$ yields a homotopy cartesian square
    \begin{equation}
        \label{eqn:diag_thm_quotient_applied_E_B}
        \begin{tikzcd}
            \Mul_{\bE_B}(b_1,\dots,b_{m+1};b)\mathop{\times}\limits_{\Mul_{\bE_B}(b_1,\dots,b_m;b)}^\mathrm{h} \{\sigma\}
            \arrow[rd, "\lrcorner" very near start, phantom]
            \ar[r]\ar[d] & 
            \Ext(\sigma) \ar[d]\\
            \{b_{m+1}\} \ar[r] & 
            B,
        \end{tikzcd}
    \end{equation}
    Upon taking base change of $\Ext(\sigma)\to B$ along the inclusion $B_{[b]}\to B$  of the connected component   of $b$ and using proposition \ref{prop:spaces_over_BG_are_G-quotients}, square \eqref{eqn:diag_thm_quotient_applied_E_B} endows the space of strict extensions (at the top left corner of equation \eqref{eqn:diag_thm_quotient_applied_E_B}) with an $\Omega_b B$-principal $\infty$-bundle structure over $\Ext(\sigma)$.

    If $b_{m+1}$ does not belong to the connected component of $b$ in $B$, then the corresponding fiber of $\Ext(\sigma)_{[b]}$ over $B_{[b]}$ is empty. In particular, the structural map $\Ext(\sigma)\to B$ factors through $B_{[b]}$.
    Now choose a point $b_{m+1}\in B_{[b]}$.
    Through the identification given by lemma \ref{lm:computation_strict_version_extensions_E_B}, the $\Omega_b B$-action on the space of strict extensions is the regular action on the first factor of $\Omega_b B\times \bigvee^m S^{n-1}$.
    Taking the quotient by this action, we see that the space $\Ext(\sigma)$ is equivalent to $\bigvee^m S^{n-1}$.
\end{proof}

\begin{proof}[Proof of theorem \ref{thm:B-framed_little_disks_coherent_operad}]
    First, it is clear that the $\infty$-operad $\bE_B^\o$ is unital. Moreover, its underlying $\infty$-category $ \mathrm{B}$ is a Kan complex by assumption.
    It remains to prove condition \ref{item:df_coherence_c} of definition \ref{df:coherence}. 
    By lemma \ref{lm:computation_extensions_E_B},  
    for a sequence of composable active morphisms $X \stackrel{\sigma}{\to} Y \stackrel{\tau}{\to} Z$ in  $\bE_B^\o$, with $X$, $Y$ and $Z$ of arity respectively $m$, $k$ and $1$,
    diagram \eqref{eqn:diag_df_coherence_rewritten} is equivalent  in the homotopy category of spaces to a  commutative square of the form
    \begin{equation}
        \begin{tikzcd}
            \label{eqn:diag_coherence_E_B}
            \mathop{\coprod}\limits^k S^{n-1} \ar[r]\ar[d] & 
            \mathop{\bigvee}^k S^{n-1} \ar[d]\\
            \mathop{\coprod}\limits_{i=1}^k \mathop{\bigvee}\limits_{p(\sigma)\inv\{i\}} S^{n-1} \ar[r] & 
            \mathop{\bigvee}\limits^m S^{n-1}
        \end{tikzcd}
    \end{equation}
    which is easily seen to be homotopy cocartesian (as in the case of the little disks $\infty$-operad $\bE_n$). This concludes the proof.
\end{proof}

\subsection{Operations on spaces of branes} 
\label{sub:string_topology_with_group_action}

As we just established,  the $\infty$-operad $\bE_B^\o$ is coherent (theorem \ref{thm:B-framed_little_disks_coherent_operad}); therefore, it admits a brane action by theorem \ref{thm:main_thm}. 
This yields the following result.
\begin{cor}
    \label{cor:E_B_brane_action_cospan}
    Using the same notations as above, there is a canonical morphism of $\infty$-operads $\bE_B^\o \to \Cospan(\zS)^{\amalg}$, sending a color $b$ to the space $\Ext(\id_b)\simeq S^{n-1}$ and an operation $\sigma\colon (b_1,\dots,b_m)\to b$ to a cospan diagram
    \begin{equation}
        \begin{tikzcd}
            \Ext(\id_b)^{\amalg m} \ar[r] &
            \Ext(\sigma) &
            \Ext(\id_b) \ar[l]
        \end{tikzcd}
    \end{equation}
    in which the middle space $\Ext(\sigma)$ is equivalent to a wedge of $m$ spheres $S^{n-1}$.
\end{cor}

The previous result can be applied to multiple geometric contexts to produce operations on spaces of branes, following \cite{toen_branes, robalo}. 
Let $\zX$ be an  $\infty$-topos. There is a  canonical functor $(-)_{\mathrm{cst}}\colon \zS\to \zX$ that sends every space $Z$ to the object $Z_{\mathrm{cst}}$ obtained as the colimit of the constant diagram with shape $Z$ and value the terminal object $*$. In other words, $Z_{\mathrm{cst}}$ is the locally constant stack with value $Z$, or equivalently the Betti shape of $Z$ in $\zX$.
Through this functor, we obtain from the $\bE_{B}$-algebra structure of Corollary \ref{cor:E_B_brane_action_cospan} a corresponding algebra in the $\infty$-category $\Cospan(\zX)$.

Given any object $X\in \zX$, we may now use internal Hom objects in $\zX$ and apply the functor $\Map(-,X)\colon \zX\op\to \zX$ to $\Ext(\sigma)_{\mathrm{cst}}$. This way, one obtains an $\bE_B$-algebra structure in $\Span(\zX)$ on the object $\Map((S^{n-1})_\mathrm{cst},X)$, called the \emph{space of $\bE_B$-branes of $X$} by Toën.
To summarize, we have the following result.

\begin{cor}
    \label{cor:E_B_brane_action_span_in_topos}
    For any object $X$ in $\zX$, the space $\Map((S^{n-1})_{\mathrm{cst}}, X)$ of $\bE_B$-branes internal to $\zX$ has a canonical $\bE_B$-algebra structure in $\Span(\zX)$, with structural morphism sending an operation $\sigma$ of arity $m$ to the span
    \begin{equation}
        \begin{tikzcd}
            \Map((S^{n-1})_{\mathrm{cst}},X)^{m}& 
            \Map(\Ext(\sigma)_{\mathrm{cst}},X)\ar[l]\ar[r] &
            \Map((S^{n-1})_{\mathrm{cst}},X). 
        \end{tikzcd}
    \end{equation}
\end{cor}

\begin{rk}
    The advantage of the above construction is its generality: one might say that the $\bE_B$-algebra structure on $\Map((S^{n-1})_{\mathrm{cst}},X)$ in $\Span(\zX)$  is motivic, in the sense that it  exists before taking any sort of linear invariant (chains, cohomology, quasi-coherent shaves, K-theory, etc.). 
    This is similar to the case of Gromov--Witten invariants \cite{robalo}, where the authors use the brane action to construct Gromov--Witten invariants at a purely geometric level and are then able to apply K-theory or ordinary cohomology functors to recover the enumerative invariants in their more classical form.
\end{rk}

In particular, specializing to the case $B = \mathrm{B}SO(n)$, the $\infty$-operad $\bE_B^\o$ recovers that of framed little disks $\bE_{n}^{\mathrm{fr}}$, so that the above corollary gives the following result. 
\begin{cor}
    \label{cor:partial_answer_sullivan_voronov_conj}
    Let $X$ be a homotopy type. Then the brane space $\Map(S^{n-1},X)$ carries an $\bE_{n}^{\mathrm{fr}}$-algebra structure in $\Span(\zS)$.
\end{cor}
To recover classical string and brane topology operations,  one would  need to pass from the previous $\bE_n^\mathrm{fr}$-algebra structure in spans to a corresponding structure in chain complexes by taking singular chains on the spaces. The rest of this section is devoted to this linearization process.

\subsubsection{Inverting spans}

In many applications, it is useful to "invert" the wrong-way morphisms appearing in the spans to obtain an algebra structure in a more tractable $\infty$-category, such as that of chain complexes or of spectra. 
To make this construction more precise, we rely on the universal property of the category of spans, as established by Stefanich in \cite{stefanich_higher_sheaf_theory} (see also \cite{gaitsgory-rozenblyum} for an earlier description of this universal property).

Given an $\infty$-category $\zC$ with pullbacks, we first need to consider an $(\infty,2)$-enhancement $\sspan(\zC)$ of the $\infty$-category $\Span(\zC)$ of spans in $\zC$ (see \cite{haugseng_spans} or \cite{stefanich_higher_sheaf_theory} for a precise construction).
Next, we define the Beck--Chevalley condition.
\begin{df}[Adjointable squares, {\cite[Definition 3.4.1]{stefanich_higher_sheaf_theory}}]
    \label{df:adjointable_squares}
    Let $\sD$ be an $(\infty,2)$-category. A commutative square
    \begin{equation}
        \begin{tikzcd}
            d' \arrow[d, "\beta'"] \arrow[r, "\alpha'"] & d \arrow[d, "\beta"] \\
            e' \arrow[r, "\alpha"]                      & e                   
        \end{tikzcd}
    \end{equation}
    in $\sD$ is called \emph{vertically right adjointable} if $\beta$ and $\beta'$ admit right adjoints $\beta^R$ and $\beta'^R$ and moreover the canonical $2$-morphism
    \begin{equation}
        \alpha'\beta'^R \to\beta^R\alpha
    \end{equation}
    constructed using the unit $\id_d\to \beta^R\beta$ and the counit $\beta'^R\beta'\to \id_{e'}$, is an isomorphism.
    
    The square is said to be \emph{horizontally right adjointable} if its transpose is vertically right adjointable. If it is both vertically and horizontally right adjointable, we simply say that the square is \emph{right adjointable}.
\end{df}

\begin{df}[Beck--Chevalley condition, {\cite[Definition 3.4.5]{stefanich_higher_sheaf_theory}}]
    \label{df:beck_chevalley}
    Let $\zC$ be an $\infty$-category with pullbacks and $\sD$ be an $(\infty,2)$-category. A functor $F\colon \zC\to \sD$ is said to satisfy the \emph{left Beck--Chevalley} condition if for every cospan $x\to s\leftarrow y$ in $\zC$, the induced commutative square in $\zD$
    \begin{equation*}
        \begin{tikzcd}
            F(x\times_s y) \arrow[d] \arrow[r] & F(y) \arrow[d] \\
            F(x) \arrow[r]                     & F(s)          
        \end{tikzcd}
    \end{equation*}
    is right adjointable.
\end{df}

Using these definitions, one can characterize the $2$-functors out of $(\infty,2)$-categories of spans.
\begin{thm}[2-categorical universal property of spans, {\cite[Theorem 3.4.18]{stefanich_higher_sheaf_theory}}]
    \label{thm:universal_property_span_beck_chevalley} 
    Let $\zC$ be an $\infty$-category with pullbacks and $\sD$ be an $(\infty,2)$-category. 
    Precomposition with the canonical functor $\zC\to \sspan(\zC)$ identifies the space of $2$-functors $\sspan(\zC)\to\sD$ with the subspace of $\Map_{\Cat_{\infty}}(\zC,\sD)$ consisting of those  functors $\zC\to \sD$ that satisfy the left Beck--Chevalley condition.
\end{thm}

\begin{rk}
    Given a functor $F\colon \zC\to \sD$ satisfying the left Beck--Chevalley condition, the associated $2$-functor $\sspan(\zC)\to \sD$ sends a span $x\stackrel{f}{\leftarrow} s \stackrel{g}{\to} y$ to the morphism $f^R\circ g\colon F(x)\to F(y)$. 
\end{rk}

We briefly discuss two geometric contexts in which to apply the above results: the case of derived stacks $\zX = \dst_k$ and that of spaces $\zX = \zS$.

\subsubsection{Algebro-geometric context}

Let $k$ be a field of characteristic $0$ and $\dst_k$ the $\infty$-category of derived étale stacks over $k$.

To invert the correspondences of derived stacks that arise from the brane action, we need to restrict our attention to a particularly well-behaved class of spaces, namely that of perfect stacks introduced by Ben-Zvi--Francis--Nadler in \cite{integral_transforms}.

\begin{df}[Perfect stacks, \cite{integral_transforms}]
    \label{df:perfect_stacks}
    A derived stack $X$ is said to be \emph{perfect} if its diagonal morphism is affine and if $\QCoh(X)$ is the ind-completion of its full subcategory of perfect complexes.
    We let $\zP$ denote the full subcategory of $\dst_k$ on perfect stacks.
\end{df}

\begin{ex}
    The class of perfect stacks contain many examples of interest. For instance, every quotient $Y/G$ of a quasi-projective derived scheme $Y$ by a linear action of an affine group $G$ is perfect. Perfect stacks are moreover stable under fiber products and if $X\in \zP$, so is $\Map((K)_\mathrm{cst},X)$ for every finite simplicial set $K$.
\end{ex}

Consider the $(\infty,2)$-category $\dgCat_k^\mathrm{L}$ of (possibly large) $k$-linear presentable dg-categories, with functors preserving small colimits as morphisms.
Let  $\QCoh \colon \dst_k\to \dgCat_k^{\mathrm{L}}$ denote the functor that assigns  to every derived stack its derived $\infty$-category of quasi-coherent sheaves.

By \cite[Proposition 3.10]{integral_transforms}, the restriction of $\QCoh$ to $\zP$ satisfies the left Beck--Chevalley property and therefore extends to a $2$-functor 
$$\QCoh\colon \sspan(\zP)\to \dgCat_k^{\mathrm{L}},$$ using Theorem \ref{thm:universal_property_span_beck_chevalley}. 
Moreover, we can upgrade this $2$-functor $\QCoh$ to a  symmetric monoidal one, using \cite[Corollary 1.2.2]{stefanich_higher_sheaf_theory}.
Together with Corollary \ref{cor:E_B_brane_action_span_in_topos}, this  shows the following result.

\begin{cor}
    \label{cor:E_B_algebra_inverting_spans}
    Let $X$ be a perfect stack. Then the $\infty$-category of quasicoherent sheaves on its space of branes $\Map((S^{n-1})_{\mathrm{cst}},X)$ carries a canonical $\bE_B$-algebra structure in $\dgCat_k^{\mathrm{L}}$.
\end{cor}

This results extends results of Toën \cite[Corollary 5.1]{toen_branes} and of Ben-Zvi--Francis--Nadler \cite{integral_transforms}, which correspond to the particular case of the $\bE_n$-operad (that is $\bE_B$ for $B\simeq *$).

\subsubsection{Topological context}

We want to adapt the above construction to the case of topological spaces, in order to recover the classical string topology operations, and more generally to obtain the conjectural $\bE_n^{\mathrm{fr}}$-algebra structure on chains on spaces of branes. 

However, one immediately runs into the problem of defining \emph{functorial} umkehr (or Gysin) maps \emph{at the chain level}. In particular, this requires to handle all the coherence data needed to produce a functor from a suitable subcategory $\zP$ of $\zS$ to that of chain complexes (or suitable variants of such). One difficulty is that the subcategory $\zP$ would have to contain the free loop spaces $\LX$, which are infinite-dimensional manifolds, and the sought functor to chain complexes would have to specialize to the classical Thom--Pontryagin construction of umkehr maps upon taking homology.
It seems that no such construction has appeared in literature yet. The author plans to come back to this problem in a future work.

	\appendix

\section{Appendix}

The appendix is divided into two parts. In the first, we give definitions and results on marked simplicial sets that are used in this paper. In the second part, we gather some standard facts on group actions in the higher categorical setting, with no claim to originality. 

\subsection{Some results on marked simplicial sets} 
\label{sec:marked_simplicial_sets}

In this section, we collect various facts about marked simplicial sets that are used in the proof of theorem \ref{thm:comparaison_bo_ext_intro}.
These results are standard and well-known to specialists, with perhaps the exception of proposition \ref{prop:right_simplification_marked_anodyne}, stating that anodyne morphisms satisfy a  weak form of the right simplification property, which seems to have not appeared in the literature. This last result might be of independent interest.


\subsubsection{Some properties of marked simplicial sets}
\label{sub:some_properties_marked_simplicial_sets}

\begin{df}
    A \emph{marked simplicial set} is a pair $(X,mX)$ where $X$ is a simplicial set and $mX$ is a subset of $X_1$ that contains all degenerate edges. 
    A morphism of marked simplicial sets $(X,mX)\to (Y,mY)$ is a morphism of simplicial sets $f\colon X\to Y$ such that $f(mX) \subseteq mY$.

    The category of marked simplicial sets is denoted $\sSet^+$.
\end{df}

\begin{nota}
    Given a simplicial set $X$, one can associated three marked simplicial sets:
    \begin{itemize}
        \item the minimal marking $X^\flat$, consisting only of degenerate edges,
        \item the cartesian marking $X^\natural$ in which an edge is marked if and only if it is an equivalence,
        \item the maximal marking $X^\sharp$, containing all edges.
    \end{itemize}
    Given an edge $e$ in a marked simplicial set $Y$, we let $Y[e]$ denote the marked simplicial set obtained by further marking $e$. 
    In other words, $Y[e]$ is the initial marked simplicial set whose underlying simplicial set is $Y$ and such that both canonical maps $Y\to Y[e]$ and $e\colon \Delta^{1,\sharp}\to Y[e]$ are morphisms of marked simplicial sets.
\end{nota}

Given a category $\zC$ and a class $S$ of morphisms in $\zC$, we say that $S$ is \emph{weakly saturated} if 
it contains all isomorphisms and is closed under cobase change, transfinite composition, coproducts and retracts. 
The smallest weakly saturated class containing $S$ is denoted $\bar{S}$ and called the \emph{saturation} of $S$.

We introduce several classes of morphisms in $\sSet$ and $\sSet^+$:
\begin{itemize}
    \item the class $\mathrm{Cell} = \{\partial\Delta^n\subset \Delta^n \:\mid \:   n\in \bN\}$,
    \item the class $\mathrm{InnHorn} = \{\Lambda_k^n\subset \Delta^n \:\mid \: 0 < k < n, \: n\geq 2 \}$,
    \item the class $\mathrm{Cell}^\flat = \{\Lambda_k^{n,\flat}\subset \Delta^{n,\flat} \:\mid \: 0 \leq k < n, \: n\geq 1 \}$,
    \item the class $\mathrm{InnHorn}^\flat = \{\Lambda_k^{n,\flat}\subset \Delta^{n,\flat} \:\mid \: 0 < k < n, \: n\geq 2 \}$,
    \item the class $\mathrm{LHorn}^\sharp = \mathrm{InnHorn}^\flat \cup \{\Lambda_0^{n,\flat}[0\to 1]\subset \Delta^{n,\flat}[0\to 1]\:\mid \: n\in \bN^* \}$,
    \item the class $\mathrm{RHorn}^\sharp = \mathrm{InnHorn}^\flat \cup \{\Lambda_n^{n,\flat}[n-1\to n]\subset \Delta^{n,\flat}[n-1\to n] \:\mid \: n\in \bN^* \}$.
\end{itemize}
The saturations $\bar{\mathrm{Cell}}$ and $\bar{\mathrm{InnHorn}}$ are respectively the class of monomorphisms and that of inner anodyne morphisms.
We now introduce a notion of anodyne morphisms for marked morphisms that is suitable for our computations of section \ref{sec:comparison_with_lurie_s_spaces_of_extensions}.

\begin{df}[Marked anodyne morphisms]
    \label{df:marked_anodyne_morphisms}
    The class $\Mark$ of marked anodyne morphisms is defined as the saturation of the union of $\mathrm{LHorn}^\sharp$ and $\mathrm{RHorn}^\sharp$ together with the map
    $$ \Lambda_1^{2,\sharp} \mathop{\coprod}_{\Lambda_1^{2,\flat}} \Delta^{2,\flat} \too \Delta^{2,\sharp}$$
    as well as the maps $K^\flat \to K^\sharp$ for all Kan complexes $K$.
\end{df}

\begin{rk}[Difference with Lurie's definition]
    Beware that the previous definition differs from that \cite[Definition 3.1.1.1]{htt} in that our definition is symmetric, whereas Lurie's include $\mathrm{RHorn}^\sharp$ but not $\mathrm{LHorn}^\sharp$.
    The conceptual reason for this discrepancy is the following: Lurie's marked anodyne morphisms are examples of trivial cofibration in the cartesian model structure on $\sSet^+$, while our marked anodyne morphisms should be trivial cofibrations in an appropriate model structure of bifibrations on $\sSet^+$. 
    However, for the purpose of this work, we shall not need the full power of such a model structure.
\end{rk}
\begin{df}
    \label{df:right_anodyne_morphisms}
Morphisms satisfying definition  \cite[Definition 3.1.1.1]{htt} will be called \emph{marked right anodyne} in this thesis. The obvious dual definition gives the class of  \emph{marked left anodyne} morphisms.
\end{df}


\begin{lm}
    \label{lm:left_lifting_property_marked_anodyne}
    Every marked anodyne morphism has the left lifting property against all morphisms of the form $X^{\natural}\to *$ for $X$ an $\infty$-category.
\end{lm}
\begin{proof}
    By \cite[Proposition 3.1.1.6]{htt}, marked \emph{right} anodyne morphisms have the desired lifting property.
    By symmetry of the argument, so do marked left anodyne morphisms. Since $\Mark$ is the saturation of the class given as the union of these two types of anodyne morphisms, we deduce the result.
\end{proof}

In section \ref{sec:comparison_with_lurie_s_spaces_of_extensions}, we will use that the class of marked anodyne morphisms satisfy the following weak form of right cancellation property.
\begin{prop}[Right cancellation  property for marked anodyne morphisms]
    \label{prop:right_simplification_marked_anodyne}
    Let $i\colon A\to B$ and $j\colon B\to C$ be monomorphisms of marked simplicial sets. 
    Assume that $i$ and $j\circ i$ are marked anodyne morphisms and that $j$ is bijective on $0$-simplices.
    Then $j$ has the left lifting property with respect to all morphisms of the form $X^{\natural}\to *$ for $X$ an $\infty$-category.
\end{prop}
\begin{proof}
    Our proof is merely an adaptation to the marked simplicial setting of the argument of \cite[Theorem 1.5]{stevenson_stability_inner_fibrations} which states that the class of inner anodyne maps has the right cancellation property. 
    We give details here for completeness.
    
    Let $X$ be an $\infty$-category.
    We will show that $j$ has the left lifting property against $X^\natural\to *$. 
    Suppose we are given a map $u\colon B\to X$ of marked simplicial sets.
    By lemma \ref{lm:left_lifting_property_marked_anodyne}, that $j\circ i$ is marked anodyne allows to pick a morphism $\varphi\colon C\to X$ satisfying $\varphi \circ j \circ i = u\circ i$.
    This implies that $u$ and $\varphi \circ j$ are in the same fiber of the map 
    $i^*\colon \Map^\sharp(B,X^\natural) \too \Map^\sharp(A,X^\natural)$.
    By \cite[Proposition 3.1.3.3 and the following remark]{htt} (or more precisely a generalization thereof to arbitrary marked anodyne maps in our sense), the map $i^*$ is a trivial Kan fibration.
    We may then take a homotopy between $u$ and $\varphi\circ j$ over their common image by $i^*$. This homotopy takes the form of a morphism $h\colon \Delta^{1,\sharp}\times B\to X^\natural$
    with the following properties:
    \begin{equation*}
        h|_{\{0\}\times B} = \varphi\circ j \qquad
        h|_{\{1\}\times B} = u\qquad 
        h\circ (\id_A\times i) = u\circ i \circ\mathrm{proj}_A.
    \end{equation*}
    Consequently, $h$  and $\varphi$ induce a map $w=(\varphi,h)\colon \{0\}\times C \cup \Delta^{1,\sharp}\times B \to X^\natural$.
    The problem therefore reduces to finding a lift $d\colon \Delta^{1,\sharp}\times C\to X^{\natural}$ in the diagram
    \begin{equation*}
        \begin{tikzcd}
            \{0\}\times C \cup \Delta^{1,\sharp}\times B \ar[d]\ar[r,"w"] & 
            X^\natural\\
            \Delta^{1,\sharp}\times C \ar[ru,bend right=20,dashed,"d"] & 
        \end{tikzcd}
    \end{equation*}
    for then $d|_{\{1\}\times C}$ will provide the desired lift of $u$ along $j$.
    Using the skeleton filtration on $C$, write $C(n) = B\cup \mathrm{sk_n(C)}$. 
    Note that we have the equality $B = C(0)$, since $j$ is bijective on objects.
    Working inductively, it therefore suffices to prove the existence of a lift in the following diagram:
    \begin{equation*}
        \begin{tikzcd}
            (\{0\}\times C(n+1)) \cup (\Delta^{1,\sharp}\times C(n)) \ar[d]\ar[r] & 
            X^\natural\\
            \Delta^{1,\sharp}\times C(n+1) \ar[ru,bend right=20,dashed], & 
        \end{tikzcd}
    \end{equation*}
    for every $n\in \bN$.
Since every monomorphism of marked simplicial sets   is obtained by cell attachments and edge markings,  the proof reduces to the case where the inclusion $C(n)\to C(n+1)$ is either  $\partial\Delta^{n+1,\flat}\to \Delta^{n+1,\flat}$ or $\Delta^{1,\flat}\to\Delta^{1,\sharp}$.
    Using \cite[Corollary 3.1.1.7]{htt}, one easily sees that $(\Delta^{1,\sharp}\times\Delta^{1,\flat})\cup (\{0\}\times\Delta^{1,\sharp})\to \Delta^{1,\sharp}\times \Delta^{1,\sharp}$ is marked anodyne, so the result follows for the case of the latter inclusion.
For the former inclusion, one can adapt the argument of \cite[Lemma 2.4]{stevenson_stability_inner_fibrations}:  decompose the inclusion
    $$(\{0\}\times \Delta^{n+1,\flat}) \cup (\Delta^{1,\sharp}\times\partial\Delta^{n,\flat})
    \too 
    \Delta^{1,\sharp}\times \Delta^{n+1,\flat}$$
    as a sequence of inner horn inclusions that successively add the different top-dimensional simplices, composed with the inclusion of a left marked horn $\Lambda_0^{n+1}[0\to 1]$ into $\Delta^{n+1}[0\to 1]$.
\end{proof}

\subsubsection{Calculus of pushout-joins}
\label{sub:calculus_pushout_joins}

Given maps $i\colon A\to B$ and $j\colon K\to L$ of simplicial sets, define the \emph{pushout-join} $i\boxast j$ as the map
\begin{equation}
    i\boxast j \colon A\ast L \mathop{\coprod}\limits_{A\ast K} B\ast K
    \xrightarrow[]{(i\ast \id_L, \id_B\ast j)}
    B\ast L.
\end{equation}
If $i$ and $j$ are instead maps of marked simplicial sets, then $i\boxast j$ also defines a map of marked simplicial sets.

\begin{lm}
    \label{lm:useful_property_pushoutjoin}
    Let $S$ and $T$ be two classes of morphisms, either both in $\sSet$ or in $\sSet^+$. Then $\bar{S}\boxast \bar{T} \subseteq \bar{S\boxast T}$.
\end{lm}
\begin{proof}
    For the case of $\sSet$, this is \cite[Proposition 30.12]{rezk_notes}. The case of marked simplicial sets is a straightforward adaption of the argument thereof.    
\end{proof}

\begin{lm}
		\label{lm:calcul_box_join}
		We have the following inclusions of classes of morphisms in $\sSet$ and $\sSet^+$: 
        $$\overline{\mathrm{RHorn}}\boxast\overline{\mathrm{Cell}}\subseteq \overline{\mathrm{InnHorn}}\quad  \text{ and }\quad  
		\overline{\mathrm{Cell}}\boxast\overline{\mathrm{LHorn}}\subseteq \overline{\mathrm{InnHorn}},$$
        $$\overline{\mathrm{Cell}^\flat}\boxast\overline{\mathrm{RHorn}^\sharp}\subseteq \overline{\mathrm{RHorn}^\sharp} \quad \text{ and }\quad 
	\overline{\mathrm{LHorn}^\sharp}\boxast\overline{\mathrm{Cell}^\flat}\subseteq \overline{\mathrm{LHorn}^\sharp}.$$ 
\end{lm}
\begin{proof}
    The results follow from lemma \ref{lm:useful_property_pushoutjoin} together with the following computation: for $j,k,n \in \bN$ with $0\leq j \leq n$, there are canonical isomorphisms
    \begin{eqnarray*}
        (\Lambda_j^n\subset \Delta^n)\boxast(\partial\Delta^k\subset \Delta^k)
        & \cong & (\Lambda_{j}^{n+1+k}\subset \Delta^{n+1+k}),\\
        (\partial \Delta^k\subset \Delta^k)\boxast (\Lambda_j^n\subset \Delta^n)
        & \cong & (\Lambda_{k+1+j}^{n+1+k} \subset \Delta^n).
    \end{eqnarray*}
\end{proof}

The following computation will be essential  in section \ref{sec:comparison_with_lurie_s_spaces_of_extensions}.

\begin{lm}
	\label{lm:combinatorial_inclusion}
	Let $I$ and $J$ be finite linear orders and $J_0\subset J$ a suborder. Let $i$ and $j$ denote respectively the inclusions $\emptyset\subseteq I$ and $J_0\subset J$.
	Then
	\begin{enumerate}
		\item the inclusion $$\id_I \ast j\colon (\Delta^I)^\flat\ast(\Delta^{J_0})^\sharp \to (\Delta^I)^\flat \ast (\Delta^J)^\sharp$$
		\item and the inclusion 
		$$i\boxast j \colon 
		(\Delta^J)^\sharp
		\mathop{\cup}\limits_{(\Delta^{J_0})^\sharp}
		((\Delta^{I})^\flat\ast (\Delta^{J_0})^\sharp)
		\to (\Delta^{I})^\flat\ast(\Delta^J)^\sharp$$
	\end{enumerate}
	are both marked anodyne.
\end{lm}

\begin{proof}
	For both assertions, it suffices to show the result for $J = J_0\cup \{y\}$.

	\begin{enumerate}
		\item We consider the inclusion $\id_I\ast j$. We proceed by induction on the cardinality of $I$. 
	Suppose first that $I = \emptyset$. Assuming that $y$ is not an extremum in $J$, let $y_-$ (respectively $y_+$) denote the maximum (resp. minimum) of the elements $x\in J$ such that $x<y$ (resp. $x>y$). Consider the spine inclusion $\mathrm{Sp}^J\to \Delta^J$, which is inner anodyne. Note that this map factors through the simplicial set $T = \Delta^{J_0}\cup \Delta^{y_- y y_+}$ as a inner anodyne inclusion $\mathrm{Sp}^J\to T$. As $\Delta^{J_0}\to \Delta^J$ also factors through $T$, it is enough to show that $(\Delta^J_0)^\sharp\to T^\sharp$ is marked anodyne: this follows from the two inclusions $(\Delta^{y_-y_+})^\sharp \to (\Lambda^{y_-yy_+}_{y_+})^\sharp \to (\Delta^{y_-yy_+})^\sharp$ being marked anodyne. The case where $y$ is the maximum (resp. minimum) of $J$ is analogous, replacing $T$ by $\Delta^{J_0}\cup \Delta^{y_-y}$ (resp. $\Delta^{J_0}\cup\Delta^{yy_+}$).

	For $I\neq \emptyset$, assume the result for finite linear orders of cardinality less than $I$. Let $x$ be the minimum of $I$ and let $I_0 = I\setminus \{x\}$.
	It suffices to show that the inclusion
	\[
	(\Delta^{I_0} \subset \Delta^I)\boxast (\Delta^{J_0}\subset \Delta^J) \colon 
	(\Delta^{I_0} 	\ast \Delta^J)\mathop{\cup}\limits_{\Delta^{I_0}\ast\Delta^{J_0}} 
	(\Delta^{I} 	\ast \Delta^{J_0}) \to (\Delta^{I}\ast\Delta^J)
	\]
	is inner anodyne. This follows from lemma \ref{lm:calcul_box_join}, since $\Delta^{I_0}\subset \Delta^I$ is right anodyne and $\Delta^{I_0}\subset\Delta^I$ is a monomorphism.

	\item We now turn to the inclusion $i\boxast j$. If $y > J_0$, then $j$ is marked left anodyne; using lemma \ref{lm:calcul_box_join} we obtain that $i\boxast j$ is inner anodyne, hence a marked equivalence. Otherwise, we can partition $J_0$ as $J_0^-\coprod J_0^+$ such that $J_0^- < y < J_0^+$ and $J_0^+$ is non-empty. Then $j$ factors as the composite
	\[
	\Delta^{J_0} {\too}
	\Delta^{J_0} \mathop{\coprod}_{\Delta^{J_0^+} } 
	\left(\Delta^{y}\ast \Delta^{J_0^+} \right)
	\too 
	\Delta^J
	\]
	where the first map is induced by the inclusion $e\colon {\Delta^{J_0^+}}\to \Delta^y\ast\Delta^{J_0^+}$ and the second map is $(\emptyset \subseteq \Delta^{J_0^-})\boxast e$. Since $e$ is marked right anodyne, using lemma \ref{lm:calcul_box_join}, we deduce that so is $j$ and therefore also $i\boxast j$, as desired.
	\end{enumerate}
\end{proof}





\subsection{Recollections} 
\label{sec:recollections}

\subsubsection{A note on cartesian fibrations and spaces of lifts} 
\label{sec:cartesian_fibrations}

In this section, we recall a useful characterization of cartesianity of a functor in terms of contractibility of a certain space of lifts, which is extracted from \cite{ha}.
Let $p \colon X\to S$ be an inner fibration of $\infty$-categories.

\begin{df}[$p$-cartesian edges]
A morphism $f \colon x\to y$ in $X$ is said to be \emph{cartesian} if the canonical map
\[
q_f \colon X_{/f} \too X_{/y}\times_{S_{/py}} S_{/pf}
\]
is a trivial fibration.
\end{df}

\begin{nota}
\label{not:lift}
	Given a morphism $f\colon x\to y$ and an object $z$ in $X$, base-changing $q_f$ along $z\colon *\to X$ yields a functor
	\[
	q_z \colon X_{/f} \times_{X} \{z\}\too \zD_z,
	\]
	where $\zD_z$ denotes the $\infty$-category $\left( X_{/y}\times_{S_{/py}} S_{/pf}\right) \times_{X} \{z\}$.
	The fiber of $q_z$ at an object $u$ will be denoted $\lift$ and refered to as the \emph{space of lifts} of $u$ along $f$, leaving the dependance on $(f,z,u)$ implicit in the notation. The situation is summarized in the following commutative diagram of $\infty$-categories
	\begin{equation*}	
		\begin{tikzcd}
		\zL \ar[r]\ar[d,"q_{u}" left]\arrow[dr, phantom, "\lrcorner", very near start]
		& X_{/f} \times_{X} \{z\}\ar[r]\ar[d,"q_z"]	
		\arrow[dr, phantom, "\lrcorner", very near start]
		& X_{/f} \ar[d,"q_f"] \\
		* \ar[r,"u"] & 
		\zD_z
		\ar[r]\ar[d]
		\arrow[dr, phantom, "\lrcorner", very near start]
		& X_{/y} \times_{S_{/py}} S_{/pf}\ar[d]\\
		& * \ar[r,"z"] & X.
		\end{tikzcd}
	\end{equation*}
	in which all squares are cartesian.
\end{nota}
 
We will use the following equivalent description of cartesian edges, which is essentially a rewording of Proposition 2.4.4.3 in \cite{htt} and its proof.

\begin{lm}
\label{lm:equivalent_description_cartesian_edge}
Let $f \colon x\to y$ be a morphism in $X$. 
Then $f$ is $p$-cartesian if and only if for all $z\in X$, every fiber $\lift$ of $q_z$ is contractible.
\end{lm}

\begin{proof}
	By Proposition 2.1.2.1 in \cite{htt}, the morphism $q_f$ is a right fibration, hence so are $q_u$ and $q_z$. Now note that every fiber of $q_f$ is of the form $\lift$ for some choice of objects $z$ and $u$.
	Since a right fibration is trivial if and only if each of its fibers is contractible, we get the result.
\end{proof}

\begin{rk}
\label{rk:lift_is_Kan}
The proof also shows that $\lift$ is a Kan complex, since $q_u$ is a right fibration whose codomain is a Kan complex. This justifies the use of the terminology \emph{space} of lifts for $\lift$.
\end{rk}

\begin{df}
[Cartesian fibrations]
The functor $p\colon X\to S$ is a \emph{cartesian fibration} if for all $y\in X$, every morphism $\bar{x}\to p(y)$ in $S$ admits a lift $x\to y$ along $p$ which is $p$-cartesian.
\end{df}


\subsubsection{Principal \texorpdfstring{$\infty$}{infinity}-bundles} 
\label{sec:principal_infinity-bundles}

In this section, we recall some definitions and basic properties of groupoid objects, $\infty$-groups and principal bundles in higher category theory. We essentially follow the exposition of 
\cite{nikolaus_schreiber_stevenson_2014_principal_bundles_general_theory} and make no claim to originality.\\

Let $\zT$ be an $\infty$-topos.
\begin{df}[{\cite[Definition 6.1.2.7]{htt}}]
	A \emph{groupoid object} in $\zT$ is a simplicial object $G_\bullet \colon \Delta\op \to \zT$ such that for every $n\in \bN$ and every partition $[k]\cup [k'] = [n]$ with $[k]\cap[k'] = \{*\}$, the induced diagram
	\begin{equation}
		\begin{tikzcd}
			G_n\ar[r] \ar[d] & G_k\ar[d]\\
			G_{k'}\ar[r] & G_0
		\end{tikzcd}
	\end{equation}
	is a pullback in $\zT$.
	The full subcategory of $\Fun(\Delta\op,\zT)$ on the groupoid objects is denoted $\Grpd(\zT)$.
\end{df}

\begin{df}
	For $f \colon X\to Y$ a morphism in $\zT$, one has a associated groupoid object $\check{C}(X\to Y)$ in $\zT$ called the \v{C}ech nerve of $f$
given in degree $n$ by the $(n+1)$-fold fiber product
$$
\check{C}(X\to Y)_n = X\times_Y X \times_Y \dots \times_Y X.
$$
We say that $f$ is an \emph{effective epimorphism} if it is the colimiting cocone of its \v{C}ech nerve, i.e. if we may write
$$
f \colon X\too \check{C}(f).
$$
Let $\mathrm{Eff}(\zT)$ denote the full subcategory of $\Fun(\Delta^1,\zT)$ on the effective epimorphisms.
\end{df}

\begin{prop}[{\cite[Corollary 6.2.3.5]{htt}}] The \v{C}ech nerve construction provides an equivalence of $\infty$-categories
\begin{equation}
	\label{eqn:equivalence_grpd(T)_effective_epi}
	\check{C}\colon \mathrm{Eff}(\zT)\simeq
	\Grpd(\zT) 
\end{equation}
whose inverse sends a groupoid $G_\bullet$ to the colimiting cocone $G_0\to \colim G_\bullet$.
\end{prop}

\begin{df}
	An \emph{$\infty$-group in $\zT$} is a groupoid $G_\bullet$ in $\zT$ such that $G_0\simeq *$.
	The corresponding full subcategory of $\Grpd(\zT)$ is denoted $\Grp(\zT)$.
	We usually write $G$ for the space $G_1$ and will often abuse notation by refering to $G$ as the $\infty$-group, leaving the rest of the simplicial structure $G_\bullet$ implicit.
\end{df}

We now recall the delooping equivalence.
\begin{prop}[{\cite[Lemma 7.2.2.11]{htt}}]
	\label{prop:B_Omega_equivalence}
	The loop space functor $\Omega$ canonically extends to a functor from the $\infty$-category $\zT_*$ of pointed objects in $\zT$  to $\Grp(\zT)$.
	Its restriction to connected pointed objects yields an equivalence of $\infty$-categories
	\begin{equation*}
		\Omega\colon (\zT_*)_{\geq 1} \simeq \Grp(\zT)\colon \mathrm{B}.
	\end{equation*}
\end{prop}
The functor $\mathrm{B}$ inverse to $\Omega$ is called the \emph{delooping}, or \emph{classifying space} functor.
The effective epimorphism $\ast\to \mathrm{B}G$ associated to an $\infty$-group $G$ is the colimiting cocone
\begin{equation*}
	\label{eqn:colimit_BG}
	\begin{tikzcd}
		\dots \ar[r, shift right=3]\ar[r, shift left=3] \ar[r, shift right]\ar[r, shift left]& 
		G\times G \ar[r]\ar[r, shift right=2]\ar[r, shift left=2] &
		G \ar[r, shift left]\ar[r, shift right] & \ast \ar[rr] && \mathrm{B}G
	\end{tikzcd}
\end{equation*}
induced by the simplicial object $G\colon \Delta\op\to \zT$.

\begin{df}[{\cite[Definition 3.1]{nikolaus_schreiber_stevenson_2014_principal_bundles_general_theory}}]
	Let $G_\bullet\in \Grp(\zT)$ be a group object and $X$ an object in $\zT$. A \emph{$G$-action on $X$} is a groupoid object $(X//G)_\bullet$ in $\zT$ of the form
	\begin{equation*}
		\label{eqn:df_G-action}
		\begin{tikzcd}
			\dots \ar[r, shift right=3]\ar[r, shift left=3] \ar[r, shift right]\ar[r, shift left]& 
			X\times G\times G \ar[r]\ar[r, shift right=2]\ar[r, shift left=2] &
			X\times G \ar[r, shift left]\ar[r, shift right,"\mathrm{proj}"below] & X
		\end{tikzcd}
	\end{equation*}
	such that the degreewise projection maps $X\times G^n\to G^n$ yield a morphism of groupoid objects $(X//G)_\bullet \to G_\bullet$.
	The \emph{$\infty$-quotient} of the action is the colimit object $X//G := \colim (X//G)$ in $\zT$.\\
	The $\infty$-category $G\mathrm{Action}(\zT)$ of $G$-actions in $\zT$ is the full subcategory of $\Grpd(\zT)_{/G_\bullet}$ on $G$-actions.
\end{df}

\begin{df}[{\cite[Definition 3.4]{nikolaus_schreiber_stevenson_2014_principal_bundles_general_theory}}]
	\label{df:G-principal_bundle}
	Let $G_\bullet\in \Grp(\zT)$ be a group object and $X$ an object in $\zT$. A \emph{$G$-principal $\infty$-bundle} over $X$ is a morphism $Y\to X$ in $\zT$ together with a $G$-action on $Y$, such that $Y\to X$ exhibits $X$ as the quotient $Y//G$.\\
	The $\infty$-category $G\mathrm{Bun}(X)$ of $G$-principal $\infty$-bundles over $X$ is the homotopy fiber at $X$ of the quotient functor
	\[
	G\mathrm{Action}(\zT) \subseteq \Grpd(\zT)_{/G_\bullet} \too \Grpd(\zT) \stackrel{\colim}{\too} \zT. 
	\]
\end{df}

The following result will be useful in section \ref{sec:homotopy_type_of_spaces_of_extensions}.
\begin{prop}
	[{
		\cite[Proposition 3.8]{nikolaus_schreiber_stevenson_2014_principal_bundles_general_theory}
	}]
	\label{prop:spaces_over_BG_are_G-quotients}
	If $G$ is an $\infty$-group and $X\to \mathrm{B}G$ a morphism in $\zT$,  then its homotopy fiber $Y\to X$ at the distinguished point of $\mathrm{B}G$ carries a canonical structure of a $G$-principal $\infty$-bundle over $X$. 
\end{prop}

The $G$-principal $\infty$-bundle structure is obtained by considering the following morphism of augmented simplicial objects
\begin{equation*}
	\begin{tikzcd}
		\dots \ar[r, shift right=3]\ar[r, shift left=3] \ar[r, shift right]\ar[r, shift left]& 
		Y\times G\times G \ar[r]\ar[r, shift right=2]\ar[r, shift left=2] \ar[d]&
		Y\times G \ar[r, shift left]\ar[r, shift right,"\mathrm{proj}"below] \ar[d]& Y \ar[r] \ar[d]& X\ar[d]\\
		\dots \ar[r, shift right=3]\ar[r, shift left=3] \ar[r, shift right]\ar[r, shift left]& 
		G\times G \ar[r]\ar[r, shift right=2]\ar[r, shift left=2] &
		G \ar[r, shift left]\ar[r, shift right] & \ast \ar[r] & \mathrm{B}G\\
	\end{tikzcd}
\end{equation*}
in which all rectangle are cartesian squares.

Moreover, all principal $\infty$-bundles are obtained through this construction, as stated in the next result.

\begin{thm}[Classification of \texorpdfstring{$G$}{G}-principal \texorpdfstring{$\infty$}{infinity}-bundles, {{\cite[Theorem 3.17]{nikolaus_schreiber_stevenson_2014_principal_bundles_general_theory}}}]
	\label{thm:classification_principal_bundles}
	For all $\infty$-groups $G\in \Grp(\zT)$ and all objects $X\in \zT$, there is a natural equivalence of $\infty$-groupoids
	\[
	G\mathrm{Bun}(X) \simeq \Map_{\zT}(X,\mathrm{B}G)
	\]
	given on objects by the construction $(p\colon X\to \mathrm{B}G)\mapsto (\mathrm{hofib}(p)\to X)$ of proposition \ref{prop:spaces_over_BG_are_G-quotients}.
\end{thm}

\ifdraft{
\subsection{todo} 
\label{sec:todo}

\todo[inline]{
	Notation à harmoniser, à regrouper au début et flèches à mettre\\
	Tout à la fin : plus de rappels catégories (structures de modèles, structure cartésienne faible) en appendice\\
}


}{} 

\makeatletter
\providecommand\@dotsep{5}
\makeatother
\listoftodos\relax

	\bibliographystyle{alpha}
	\bibliography{ms}

\end{document}